\documentclass{article}
\usepackage[applemac]{inputenc}
\usepackage[T1]{fontenc}
\usepackage{lmodern}
\usepackage[a4paper]{geometry}
\usepackage[english]{babel}
\usepackage{graphicx}
\usepackage{onlyamsmath}
\usepackage{amsfonts}
\usepackage{amssymb}
\usepackage{ntheorem}

\newcommand{\e}{\epsilon}
\newcommand{\se}{\sqrt{\epsilon}}
\newcommand{\ga}{\gamma}
\newcommand{\gaa}{\lvert\gamma\rvert}
\newcommand{\sga}{\sqrt{\gamma}}
\newcommand{\sgap}{\sqrt{\gamma'}}
\newcommand{\sgaa}{\sqrt{\lvert\gamma\rvert}}

\newcommand{\E}{\mathbb{E}}
\newcommand{\esp}{\mathcal{H}^\omega}

\newcommand{\ko}[1]{k^2 (\omega_{#1})}
\newcommand{\N}[1]{N (\omega{#1})}

\newcommand{\Bh}[2]{\beta_{#1} (\omega{#2})}

\newcommand{\A}[2]{\widehat{a}_{#1}^\epsilon (\omega,#2)}

\newcommand{\B}[2]{\widehat{b}_{#1}^\epsilon (\omega,#2)}

\newcommand{\hf}{\widehat{g}(\omega)}
\newcommand{\1}{\textbf{1}}
\newcommand{\dz}{\frac{d}{dz}}

\newtheorem{thm}{Theorem}[section]
\newtheorem{defi}{Definition}[section]
\newtheorem{prop}{Proposition}[section]

\newtheorem{lem}{Lemma}[section]
\newtheorem{assumption}{Assumption}
\theoremstyle{nonumberplain}
\theorembodyfont{\upshape}
\newtheorem{preuve}{Proof}

\title{Wave Propagation in Underwater Acoustic Waveguides with Rough Boundaries.}

\author{Christophe Gomez\thanks{Laboratoire d'Analyse, Topologie, Probabilit\'es, UMR 7353, Aix-Marseille Universit\'e, Marseille, France, christophe.gomez@latp.univ-mrs.fr.}}

\date{} 

\begin{document}

\maketitle

\begin{abstract}

In underwater acoustic waveguides a pressure field can be decomposed over three kinds of modes: the propagating modes, the radiating modes and the evanescent modes. In this paper, we analyze the effects produced by a randomly perturbed free surface and an uneven bottom topography on the coupling mechanism between these three kinds of modes. Using an asymptotic analysis based on a separation of scales technique we derive the asymptotic form of the distribution of the forward mode amplitudes. We show that the surface and bottom fluctuations affect the propagating-mode amplitudes mainly in the same way. We observe an effective amplitude attenuation which is mainly due to the coupling between the propagating modes themselves. However, for the highest propagating modes this mechanism is stronger and due to an efficient coupling with the radiating modes.

\end{abstract}

\begin{flushleft}
\end{flushleft}

\begin{flushleft}
\end{flushleft}

\section*{Introduction} 

Acoustic wave propagation in waveguides has been studied for a long time because of its numerous domains of applications. One of the most important application is submarine detection with active or passive sonars, but it can also be used in underwater communication, mines or archaeological artifacts detection, and to study the ocean's structure or ocean biology. Underwater acoustic waveguide are used to model acoustic wave propagation media such as a continental shelves. These environments are very complex media because of indices of refraction with spatial and time dependences. However, the sound speed in water, which is about $1500$m/s, is sufficiently large with respect to the motions of water masses so that we can consider this medium as being time independent. Moreover, the presence of small spatial inhomogeneities in the water, heave or the ocean bottom roughness can induce significant effects over large propagation distances. 

\begin{figure}\begin{center}
\begin{tabular}{cc}
\includegraphics*[scale=0.2]{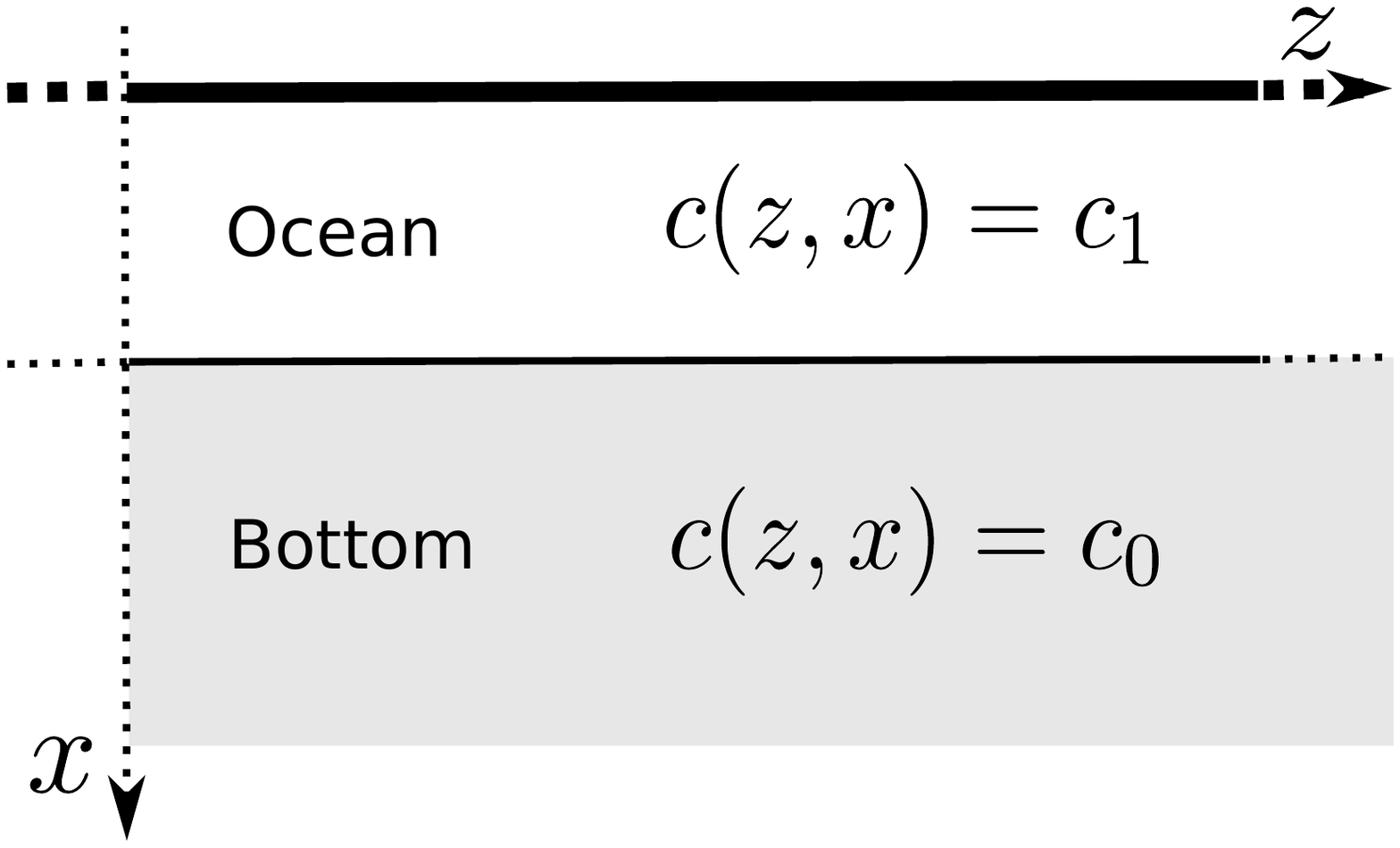} & \includegraphics*[scale=0.2]{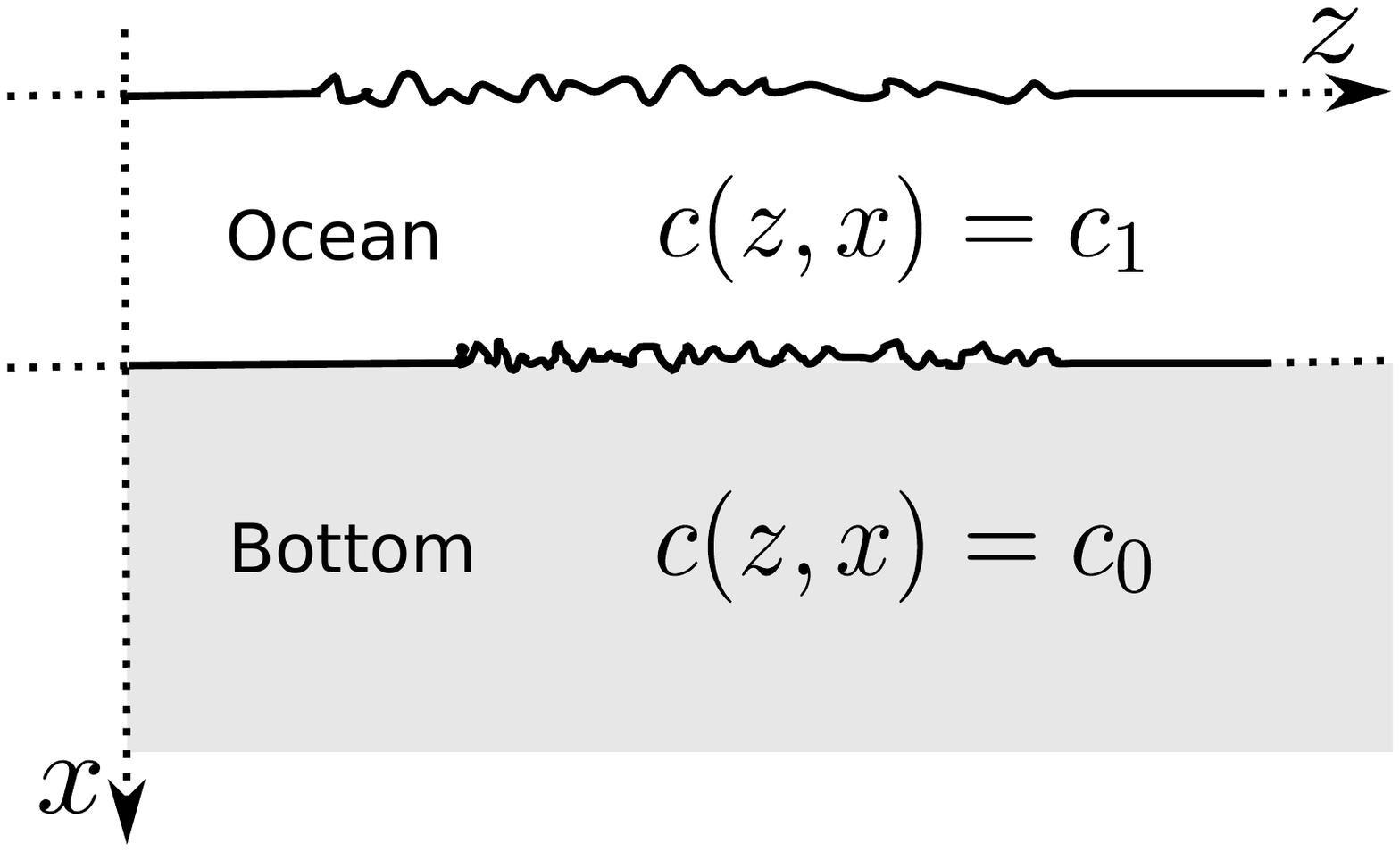}\\
$(a)$ & $(b)$
\end{tabular}
\end{center}
\caption{\label{wavegmod} Illustration of two semi-infinite waveguides with a Pekeris profile. In $(a)$ an unperturbed waveguide with the rigid-lid assumption is considered. In $(b)$ a perturbed waveguide with a free surface and an uneven bottom topography is considered. In the two cases the waveguide are homogeneous in the interior with propagation speed profile $c(z,x)$ equal to $c_1$ in the ocean section of the waveguide, and $c_0$ in its bottom.}
\end{figure}

The effects of random boundaries has been studied in many physical setups. For instance, for fluid flows in a medium with random boundaries has been studied in \cite{basson}, for water wave propagation with a free surface or a random depth has been studied in \cite{dya,garniernach}, and also in wave propagation in underwater acoustics with a perturbed sea surface \cite{dozier, kuperman, kuperman2, kuttler}. Mathematical studies regarding acoustic wave propagation in randomly perturbed waveguides have been carried out in many papers \cite{book,garniereva, garnier2,papa, gomez2,gomez3,papanicolaou}, but only under the rigid-lid assumption at the waveguide boundaries, and with random perturbations inside the waveguide through variations of the index of refraction. Recently, waveguides with a bounded cross-section and randomly perturbed boundaries have been considered in \cite{alonso}. Using a change of coordinates the authors show that the scattering effects differ from the ones produced by internal perturbations (see \cite{book}). In this paper, we consider a two-dimensional acoustic waveguide model (unbounded cross-section) with a Pekeris profile, with a randomly perturbed free surface and uneven bottom topography (see Figure \ref{wavegmod}). Because of the unbounded cross-section of our waveguide model, we consider a conformal transformation approach. In fact, the conformal transformation is a convenient tool allowing the use of the modal decomposition of the unperturbed waveguide (see Figure \ref{wavegmod} $(a)$). In our model a propagating field can be decomposed over three kinds of modes: the propagating modes which propagate over long distances, the evanescent modes which decrease exponentially with the propagation distance, and the radiating modes representing modes which can penetrate under the ocean bottom. Using an asymptotic analysis based on a separation of scale technique precisely described in Section \ref{sect1P2}, we show in this paper that the random perturbations of the waveguide geometry induce a mode coupling between the three kinds of modes. Theorem \ref{thasymptP21}, Theorem \ref{thasymptP12}, and Theorem \ref{thasympP2} describe the coupling mechanisms between the propagating and the radiating mode amplitudes in term of diffusion models but taking into account the effects produced by the evanescent modes. These mechanisms are similar to the ones observed in \cite{gomez2} and produce an effective attenuation of the mode amplitudes (Theorem \ref{thmain} and see Section \ref{meanmode}). It turns out that the surface and bottom fluctuations affect the propagating-mode amplitudes mainly in the same way. However, for almost all the propagating modes the attenuation mechanism is mainly due to the coupling between the propagating mode themselves. Nevertheless, for the highest propagating modes the attenuation mechanism is due to the coupling with the radiating modes, and it is significantly stronger than for the other modes.

The organization of this paper is as follows. In Section \ref{sect1P2} we present the waveguide model and give a summary of the tools. We also introduce some simplifying assumptions used in this paper before introduce the main result. The remaining of the paper consists in introducing more precisely the tools, and techniques used to prove the results. In Section \ref{conformal} we introduce the conformal transformation. In Section \ref{sect3P2} we study the mode coupling mechanism when the three kinds of modes are taken into account and we derive the coupled mode equations. In Section \ref{coupledprocP2}, under the forward scattering approximation, we study the asymptotic form of the joint distribution of the propagating and radiating mode amplitudes. Finally, we describe the attenuation of the propagating-mode amplitudes.

\section{Waveguide Model and Main Result}\label{sect1P2}

\begin{figure}\begin{center}
\begin{tabular}{cc}
\includegraphics*[scale=0.2]{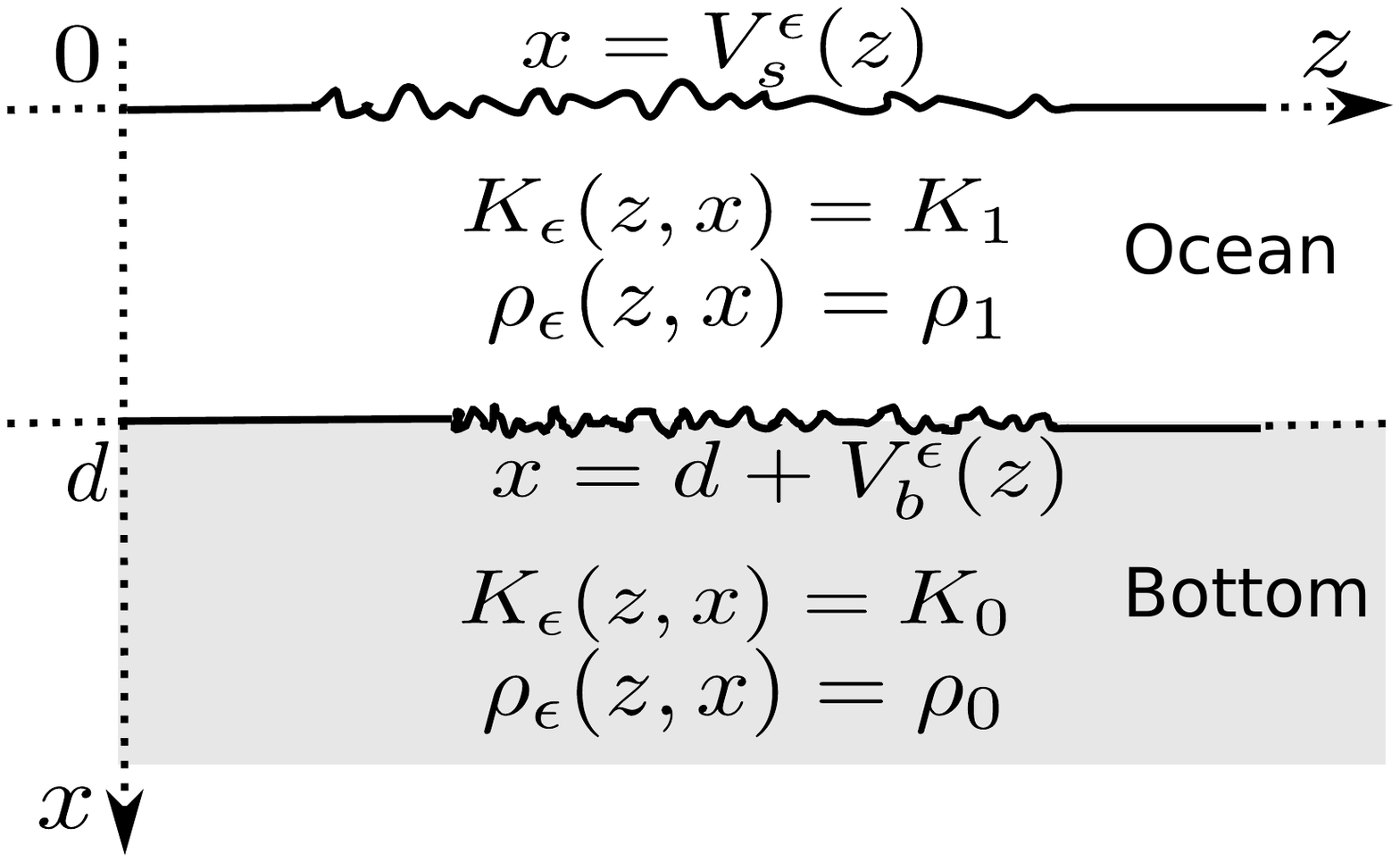} & \includegraphics*[scale=0.2]{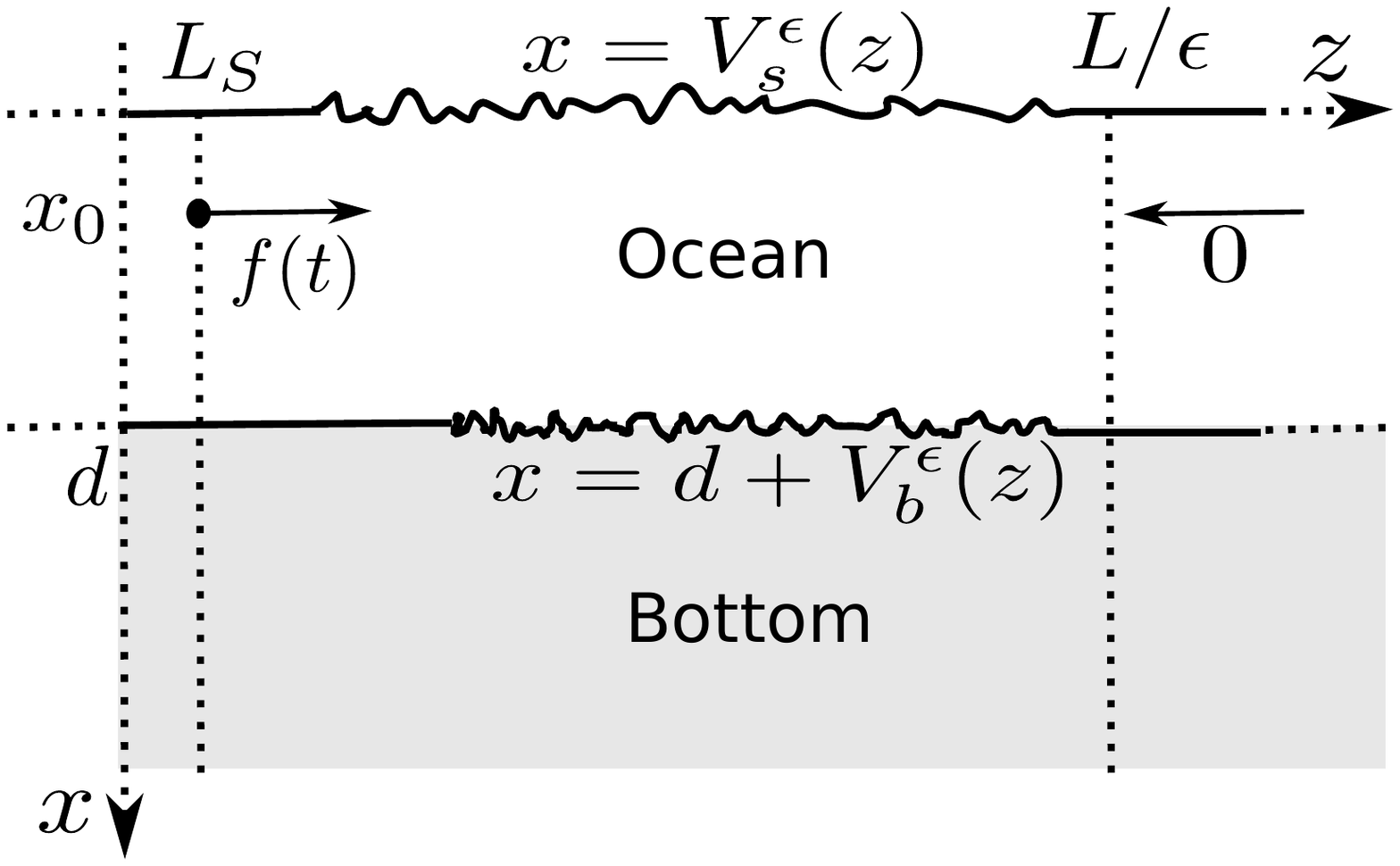}\\
$(a)$ & $(b)$
\end{tabular}
\end{center}
\caption{\label{figureP2} Illustration of the waveguide model. In both figures $d$ represents the mean ocean depth, and the random fluctuations are given by the graphs of $V^\e_s$ at the free surface,  and $d+V^\e_b$ at the bottom of the waveguide model. In $(a)$ we illustrate the medium parameters of the waveguide model and its random fluctuations at the free surface and at the bottom. In this figure, the bulk modulus $K_\e(z,x)$ (resp. $\rho_\e(z,x)$) is equal to  $K_1$ (resp. $\rho_1$) in the ocean section of the waveguide and $K_0$ (resp. $\rho_0$) in its bottom. In $(b)$ we illustrate the source term $\textbf{F} (t,z,x)$, which is approximately a point-like source at $(L_S,x_0)$ emitting a wave $f(t)$ is the $z$-direction. $L/\e$ characterizes the size of the section in which the random fluctuations are included (see \eqref{randompert}). The right arrow (with a $0$) pointing to the left at $z=L/\e$ indicates that no wave is coming from the right at the end of the random section. }
\end{figure}

Throughout this paper, we consider a two-dimensional linear acoustic wave model. In acoustic waveguides sound propagation can be described by Euler's equations \cite[Appendix A]{kat}, so that the linearized conservation equations of mass and momentum are given by \cite[Chapter 2]{book}
\begin{equation}\label{conservationP2}
\begin{split}
\frac{1}{K_\e(z,x)} \partial_t p+ \nabla . \textbf{u} &=0,  \\
\rho_\e ( z,x)\partial_t \textbf{u} + \nabla p &=\textbf{F},
\end{split}
\end{equation}
where $p$ is the acoustic pressure, $\textbf{u}$ is the acoustic velocity, $\rho_\e$ is the density of the medium, and $K_\e$ is the bulk modulus. The coordinate $z$ represents the propagation axis along the waveguide, and the coordinate $x$ represents the transverse section of the waveguide (see Figure \ref{wavegmod} or Figure \ref{figureP2}). The forcing term $\textbf{F} (t,z,x)$ is given by
\begin{equation}\label{sourcetermP2}  \textbf{F} (t,z,x)=\Psi(t,x)\delta(z-L_S)\textbf{e}_{z},\end{equation}
where $\textbf{e}_{z}$ is the unit vector pointing in the $z$-direction. The source profile $\Psi(t,x)$ is defined by \eqref{profsourceP2}, so that $\Psi(t,x)\simeq f(t)\delta(x-x_0)$.
Consequently, the forcing term 
\begin{equation}\label{approxsource}\textbf{F} (t,z,x)\simeq f(t)\delta(x-x_0)\delta(z-L_S)\textbf{e}_{z}\end{equation}
 models a source which is closed to a point-like source located at $(L_S,x_0)$, and emitting a signal $f(t)$ in the $z$-direction (see Figure \ref{figureP2} $(b)$).
Let $d>0$ be the average ocean depth, we assume that the medium parameters are given by:

\begin{equation}\label{parametre}\begin{split}
K_\e(z,x) & =  \left\{ \begin{array}{ccl} 
                                            K_1& \text{ if }  & x\in(V^\e_s(z),d+ V^\e_b (z)), \\
                                             K_0 & \text{ if }  &x\in [d+ V^\e_b (z),+\infty),\\
                                          \end{array} \right. z\in (-\infty,+\infty), \\
\rho_\e(z,x)&=   \left\{ \begin{array}{ccl} 
                                            \rho_1& \text{ if }  & x\in(V^\e_s(z),d+ V^\e_b (z)), \\
                                             \rho_0 & \text{ if }  &x\in[d+ V^\e_b (z),+\infty),\\
                                          \end{array} \right. z\in (-\infty,+\infty),
\end{split}\end{equation}
where $V^\e_s$ and $V^\e_b$ model respectively the free surface and the bottom topography (see Figure \ref{figureP2} $(a)$). In our context, the free surface can model heave produced by the speed of the wind, and the uneven bottom topography can model a sandy bottom with variations produced by water currents \cite{kat}.   

From the continuity of the pressure and the vertical component of the particle velocity across the transition from the ocean to the bottom, we have
\[\lim_{x\to (d+ V^\e_b (z))^+}p(t,x,z)=\lim_{x\to (d+ V^\e_b (z))^-}p(t,x,z)\]
and
\[\lim_{x\to (d+ V^\e_b (z))^+}\partial_x p(t,x,z)=\lim_{x\to (d+ V^\e_b (z))^-}\partial_x p(t,x,z),\]
so that we obtain from equations \eqref{conservationP2} the wave equation for the pressure field 
\begin{equation}\label{waveeq}
\Delta p - \frac{1}{ c^2_\e(z,x)  }\frac{\partial ^2  p }{\partial t^2} = \nabla .\textbf{F},
\end{equation}
with velocity field
\begin{equation}\label{propspeed}
  c_\e(z,x)= \sqrt{K_\e(z,x)/\rho_\e(z,x)} =\left\{ \begin{array}{ccl} 
                                            c_1 & \text{ if }  & x\in(V^\e_s(z),d+ V^\e_b (z)), \\
                                             c_0 & \text{ if }  &x\in [d+ V^\e_b (z),+\infty),\\
                                          \end{array} \right. z\in (-\infty,+\infty),
\end{equation}
and $\Delta = \partial ^2 _x + \partial ^2 _{z}$ (see Figure \ref{wavegmod} $(b)$ for an illustration).

Moreover, from the continuity of the pressure field at the free surface $x=V^\e_s(z)$, which is tantamount to neglect the surface tension, the wave equation \eqref{waveeq} is complemented by the following boundary conditions 
\[p(t,z,V^\e_s(z))=p_{surface}\quad \forall(t,z)\in[0,+\infty)\times(-\infty,+\infty),\] 
where $p_{surface}$ is the atmospheric pressure. However, one can assume without loss of generality that 
\[p_{surface}=0,\] 
and then consider Dirichlet boundary conditions at the free surface, which corresponds to a pressure-release condition.  

 In this paper we consider the Pekeris waveguide model \eqref{propspeed}. This kind of model has been extensively studied for half a century \cite{pekeris}, and has been widely used to model an ocean with a constant propagation speed profile (see Figure \ref{wavegmod}). Such conditions can be found during the winter in Earth's mid latitudes and in water shallower than about $30$ meters \cite{kuperman3}. The Pekeris profile leads us to simplified algebra but it underestimates the complexity of the medium. However, the analysis carried out in this paper can be extended to more general propagation speed profiles. This model could also be used to study the propagation of electromagnetic waves in a dielectric slab or an optical fiber with randomly perturbed boundaries. In fact, in the electromagnetic setting the bulk modulus and the density of the medium become the permeability and conductivity \cite{magnanini,marcuse,rowe,wilcox2}.

In the definition of the medium parameters \eqref{parametre} and \eqref{propspeed}, $V^\e_s$ and $V^\e_b$ are given by
\begin{equation}\label{randompert}
 V^\e_s(z)= \se V_s( z/l_c) f_s(\e z), \quad\text{and}\quad V^\e_b(z)= \se V_b(z/l_c) f_b(\e z), \quad\text{with}\quad \e\ll1. 
\end{equation} 
We assume that $V_s$ and $V_b$ are two independent mean zero stationary bounded stochastic processes, and $f_s$ and $f_b$ are two smooth functions with support included in $(\e\eta, L)$ (with $\eta>0$). Therefore, $f_s(\e z)$ and $f_b(\e z)$ represent the locations of the random fluctuations of the free surface and the bottom topography for $z\in(\eta,L/\e)$. The scaling of $f_b$ and $f_s$ as been chosen according to the size of the random fluctuations ($\se\ll1$). In fact, we will see that we have to wait for long propagation distances (of order $1/\e$) to observe significant cumulative stochastic effects on the pressure wave (see \eqref{pressureinit} and Theorem \ref{thmain}). For instance, wind speeds at about $5$m/s induce a standard deviation of the surface roughness at about $0.1$m, and the same order of magnitude can be considered for the roughness of sandy bottoms (see \cite[Chapter 2]{kat} and \cite{kuperman2}). Moreover, the typical propagation distance is about $10km$. The correlation lengths $l_c$ of the surface and bottom roughness in \eqref{randompert} are considered to be about $10$m, that is large compared to the roughness magnitude \cite[Chapter 2]{kat}. In fact, in most ocean-surface conditions the ratio between the standard deviation and the correlation length is small, and regarding the bottom roughness the correlation length is larger than its standard deviation (see \cite{kuperman2}). Let us note that the two random processes $ V_s$ and $V_b$ in \eqref{randompert} are considered to be independent. In our context, the reason is that the surface and bottom standard deviations are small compared to the ocean depth ($d\sim30$m), so that the bottom topography does not produce any ripples on the free surface. In this paper we consider the power spectral densities for $V_s$ and $V_b$ (i.e the Fourier transform of their autocorrelation functions) given by (see \cite[Chapter 2, Section 2.4 and Section 2.9]{kat})
\begin{equation}\label{formuleI}
I_s(u)=\int R_s(z)e^{iuz}dz= \frac{C_{1,s}}{\vert u\vert ^{\alpha_{I_s}}}e^{-\frac{C_{2,s}}{\vert u\vert^{\mu_{I_s}}}}\quad\text{and}\quad I_b(u)=\int R_b(z)e^{iuz}dz=\frac{C_b}{4\pi(v^2+1)^{\alpha_{I_b}/2}},
\end{equation}
with $C_{1,s},C_{2,s}, C_b>0$, and 
\[R_s(z)=\E[V_s(z+z_0)V_s(z_0)] \qquad \text{and}\qquad R_b(z)=\E[V_b(z+z_0)V_b(z_0)].\]
In \cite[Chapter 2, Section 2.4 and Section 2.9]{kat} the authors refer to $\alpha_{I_b}=4$, $\alpha_{I_s}=6$ and $\mu_{I_s}=2$ (Pierson-Neuman spectra), or $\alpha_{I_s}=5$ and $\mu_{I_s}=4$ (Pierson-Moscovitz spectra). Consequently, the trajectories of $V^\e_s$ and $V^\e_b$ are at least of class $\mathcal{C}^2$  with bounded derivatives, and therefore the surface fluctuations given by the graph $V^\e_s$ do not have breaking waves. Finally, we have a free surface and a bottom topography model describing a smooth ocean surface and bottom with small fluctuations.

Regarding the mixing properties of $V^\e_s$ and $V^\e_b$, we assume that $V_s$ and $V_b$ in \eqref{randompert} are $\phi$-mixing processes \cite{kushner}, that is considering
\begin{equation}\label{filtration}\mathcal{F}_{z}=\mathcal{F}_{0,z}=\sigma(V_s(u), V_b(u), \quad 0\leq u\leq z )\quad \text{and}\quad \mathcal{F}_{z,+\infty}=\sigma(V_s(u), V_b(u), \quad z\leq u ),\end{equation}
the $\sigma$-algebras generated by $(V_s(u), V_b(u))_{0\leq u\leq z}$ and $(V_s(u), V_b(u))_{z\leq u}$, we have 
\begin{equation}\label{mixcond}
\sup_{\substack{z\geq 0\\A\in \mathcal{F}_{z+u,+\infty}\\B\in\mathcal{F}_{0,z}}}\lvert \mathbb{P}(A\vert B)-\mathbb{P}(A)\rvert\leq \phi(u).
\end{equation}
This mixing property \eqref{mixcond} describes the decorrelating behavior of the random free surface $V_s$ and bottom topography $V_b$, where $\phi\in L^1(\mathbb{R})\cap L^{1/2}(\mathbb{R})$ is a nonnegative function characterizing the decorrelation speed. In other words, the mixing condition \eqref{mixcond} describes the decreasing dependence between $\mathcal{F}_{z}$ and $ \mathcal{F}_{z+Z,+\infty}$ as $Z\to+\infty$, that is for $Z$ large enough  $\mathcal{F}_{z}$ and $ \mathcal{F}_{z+Z,+\infty}$ are almost independent. 

In this paper we propose an asymptotic analysis based on a separation of scale technique. The important scales of the problem considered in this paper are : the typical wavelength $\lambda$ of the pulse $f(t)$ in \eqref{approxsource}, the propagation distance $L_\e=L/\e$, the correlation length $l_c$, and the standard deviation $\sigma_\e=\se$ of the surface and bottom roughness. These parameters are chosen so that we are in the well-known weakly heterogeneous regime \cite[Chapter 5]{book}: 
\begin{equation}\label{regime} \lambda^{-1} l_c\sim 1,\quad \lambda^{-1} L_\e \gg 1, \quad\text{and}\quad \sigma_\e\ll1. \end{equation}
In \eqref{regime}, the first relation means that the typical wavelength and the correlation length of the perturbations are of the same order, so that the acoustic pressure wave produced by the source can interact with all the scales of variation of the medium perturbations. The second one means that we consider a high-frequency regime, while the last one means that the amplitude of the perturbations are small. For example, the typical value of the propagation distance $L_\e$ is about $10$km, the typical wavelength in the ocean $\lambda$ and the correlation length $l_c$  are about $10$m, and the amplitude of the fluctuations are about $10$cm \cite{kuperman2,kuperman3}. In fact, as already noted, we have to wait for long propagation distances to observe significant cumulative effects on the acoustic pressure wave when the fluctuations of the free surface and the bottom topography are small. In this paper we set $l_c=1$, so that $\lambda\sim 1$. Moreover, we also have 
\[ \lambda^{-1} d \sim 1. \]
meaning that the ocean depth and the typical wavelength are of the same order. This condition allows us to have a fixed number of propagating modes ($\N{}\sim  \lambda^{-1}d$) in \eqref{modedec}. 

\begin{figure}
\begin{center}
\includegraphics*[scale=0.25]{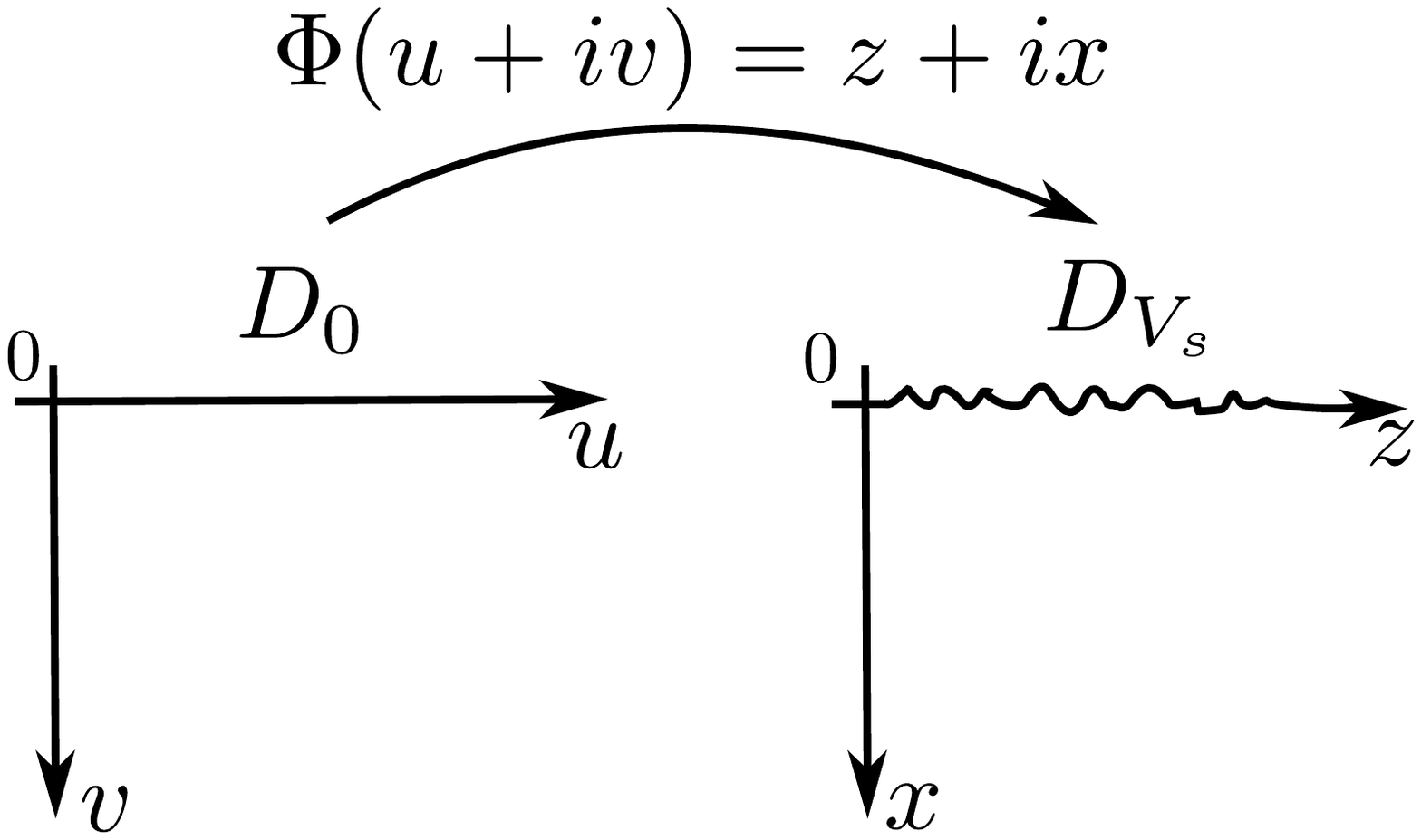} 
\end{center}
\caption{ \label{conftrans1} Schematic representation of the change of coordinate from the unperturbed domain $D_0$ onto the perturbed domain $D_{V_s}$.}
\end{figure}

To study the asymptotic behavior of the acoustic pressure wave in the perturbed waveguide domain, we use a change of coordinates through a conformal transformation \cite{rudin} described in Section \ref{conformal}, from
\[D_0=\big\{u+iv\in \mathbb{C}, \quad\text{such that }u>0 \big\}\]
onto
\[D_{V_s}=\big\{z+ix\in \mathbb{C}, \quad\text{such that } x>V^\e_s(z) \big\},\]
and allowing us to study the pressure wave $p(t,z,x)$ in a waveguide domain with a flat surface (see Figure \ref{conftrans1}). 
As a result, one can focus our attention on $p_0(t,u,v)$ defined by
\[p_0(t,u,v)=p(t,z(u,v),x(u,v)),\]
which represents the pressure wave $p(t,z,x)$ solution of \eqref{waveeq} in a waveguide domaine with a flat surface, and where $z(u,v)+ix(u,v)=\Phi(u+iv)$ is the change of coordinates  from $D_0$ onto $D_{V_s}$. Thanks to this conformal transformation the acoustic pressure wave $p_0$ satisfies now Dirichlet boundaries conditions at the flat surface of the waveguide:
\[p_0(t,u,0)=0, \quad \forall (t,u)\in[0,+\infty)\times(-\infty,+\infty).\]
Moreover, the original pressure wave $p(t,z,x)$ can be recovered from $p_0(t,u,v)$ by inverting the change of coordinates:
\[ p(t,z,x) = p_0(t,Re(\Phi^{-1}(z+ix)), Im(\Phi^{-1}(z+ix)) ).\]
We show in Section \ref{conformal} that this change of coordinates transforms the random fluctuations of the free surface to random fluctuations in the interior and at the bottom of the waveguide (see Figure \ref{figuredef}). 

Throughout this paper the propagation model is a linear model, so that the pressure wave $p_0(t,u,v)$ can be expressed as the superposition of monochromatic waves by taking its Fourier transform. Here, the Fourier transform and the inverse Fourier transform, with respect to time, are defined by
\begin{equation}\label{fourier} \hf =\int g(t)e^{i \omega t} dt, \quad g(t)=\frac{1}{2 \pi} \int \hf e^{-i \omega t} d\omega.\end{equation}
In the new system of coordinates $D_0$, and in the frequency domain one can decompose the monochromatic pressure field $\widehat{p}_0(\omega,u,v)$ as follows,
\begin{equation}\label{modedec}\widehat{p}_0(\omega,u,v)=\underbrace{\sum_{j=1}^{\N{}} \widehat{p}_j (\omega,u )\phi_j(\omega,v)}_{\text{propagating modes}}+\underbrace{\int_{0}^{\ko{}}\widehat{p}_\ga (\omega,u)\phi_\ga (\omega,v)d \ga}_{\text{radiating modes}}+\underbrace{\int_{-\infty}^{0}\widehat{p}_\ga (\omega,u)\phi_\ga (\omega,v)d \ga}_{\text{evanescent modes}},\end{equation}
where the modal decomposition corresponds to the one of the unperturbed Pekeris propagation model recalled in Section \ref{spectralP2} (see Figure \ref{wavegmod} $(a)$). Here $k(\omega)=\omega/c_0$ is the wavenumber, $c_0$ is the propagation speed in the bottom of the waveguide model, and $\N{}$ is the number of propagating modes defined by \eqref{number}. The coefficients $\widehat{p}_j (\omega,u )$ ($j\in \big\{1,\dots,\N{}\big\}$) and $\widehat{p}_\ga (\omega,u)$ ($\forall \ga\in(-\infty,\ko{})$), defined by \eqref{defmode1} and \eqref{defmode2}, represent the amplitudes of each kind of modes (propagating, radiating, and evanescent modes), Moreover, $\phi_j(\omega,u)$ ($j\in \big\{1,\dots,\N{}\big\}$) and $\phi_\ga(\omega,u)$ ($\forall \ga\in(-\infty,\ko{})$), defined by \eqref{phij} and \eqref{phiga}, are the eigenelements of the Pekeris operator defined by \eqref{pekerisop}. Thanks to this decomposition, one can study the effects of the random fluctuations of the waveguide on the pressure field $\widehat{p}_0(\omega,u,v)$ through the coupled mode equations \eqref{eqdiff2P2} introduced in Section \ref{couplemodeeq}. In \eqref{modedec}, the propagating modes can propagate over long distances, while the amplitudes of the evanescent modes decrease exponentially with the propagation distance. Finally, the radiating modes represent modes which can penetrate under the ocean bottom. 

Before introducing the main results of this paper, we need to consider some simplifications. First, the change of coordinate we consider for the analysis (Assumption \ref{hyp1} in Section \ref{conformal}) is an approximation of the conformal transformation introduced above but does not change the overall results of this paper. For technical reason, we restrict the continuous spectrum in \eqref{modedec} to $(-1/\xi,-\xi)\cup(\xi,\ko{})$ (Assumption \ref{hyp2} in Section \ref{couplemodeeq}). Here, we have introduced a new parameter $\xi$ such that
\[\e\ll\xi\ll1,\]
which means that both of these parameters will go to $0$ in the asymptotic analysis carried out in this paper, but $\e$ first, then $\xi$. Moreover, we neglect  the mode coupling mechanism between the radiating and the evanescent modes (Assumption \ref{hyp3} in Section \ref{thasyptoticP2}). This assumption is not restrictive since for random perturbations taking place in the interior of the waveguide the coupling mechanism between these modes do not play any significant role (see \cite{gomez2}). We consider a radiation condition for the evanescent modes $(\ga\in(-1/\xi,-\xi))$ (Assumption \ref{hyp6} in Section \ref{iemprP2}) meaning that the energy carried by these modes decay as the propagation distance becomes large. We also consider the forward scattering approximation (Assumption \ref{hyp4} in Section \ref{fsaP2}) which consists in neglecting the backscattering effects during the propagation.  Finally, we consider a nondegeneracy condition on the modal wavenumbers $\beta_j(\omega)$ defined by \eqref{spectrumRP2} in Section \ref{spectralP2} (Assumption \ref{hyp5} in Section \ref{thasyptoticP2}). Under these assumption and because of the strong attenuation of the evanescent modes, the pressure field $\widehat{p}_0$ at the end of the perturbed section ($u=L/\e$) can be expressed
\begin{equation}\label{pressureinit} \widehat{p}_0\Big(\omega,\frac{L}{\e},v\Big)\simeq \sum_{j=1}^{\N{}}\frac{\widehat{a}^{\e,\xi}_j(\omega, L)}{\sqrt{\beta_j(\omega)}}e^{i\beta_j(\omega)L/\e}\phi_j(\omega,v)+\int_\xi^{\ko{}} \frac{\widehat{a}^{\e,\xi}_\ga(\omega, L)}{\ga^{1/4}}e^{i\sqrt{\ga} L/\e}\phi_\ga(\omega,v)d\ga,\end{equation}
where $\widehat{a}^{\e,\xi}(\omega, u)$ are the forward mode amplitudes satisfying \eqref{eqfora}. In \eqref{pressureinit} we have introduced the modal wavenumbers $(\beta_1(\omega),\dots,\beta_{\N{}}(\omega))$ of the propagating modes and $(\sga, \ga\in(\xi,\ko{}))$ for the radiating modes, which are defined in Section \ref{spectralP2}. According to Theorem \ref{thasymptP21}, Theorem \ref{thasymptP12}, and Theorem \ref{thasympP2}, one can observe that in the asymptotic $\e \to 0$ the propagating mode amplitudes behave diffusively while the radiating mode amplitudes remain constant. The main result of this paper regarding the asymptotic behavior of the propagating-mode amplitudes is given by the following theorem.
\begin{thm}\label{thmain}
Let us assume that $f_s=f_b=\mathbf{1}_{(0,L)}$ in \eqref{randompert}. Under Assumptions \ref{hyp1}-\ref{hyp5}, and the mixing conditions \eqref{mixcond}, we have for all $j\in\{1,\dots,\N{}\}$
\[\lim_{\xi\to 0}\lim_{\e\to 0} \vert\E[ \widehat{a}^{\e,\xi}_j(\omega, L) ] \vert=e^{-C_j(\omega)L}\vert \widehat{a}_j(\omega, 0)\vert.\]
Moreover, in the limit of a large number of propagation modes $k(\omega)\gg 1$, we have :
\begin{itemize}
\item for $j\ll \N{}^{1/2}$ 
\[C_j(\omega) \underset{k(\omega)\gg1}{\sim} - C_1\, k^{3/2}(\omega)\frac{j^2}{\N{}^2}  \int_{0}^{+\infty}\sqrt{v}(I_s+I_b)(v)dv,\]
\item for $j\sim \N{}^{1/2}$
\[C_j(\omega)\underset{k(\omega)\gg1}{\sim} -  C_2  \, k^{3/2}(\omega) \frac{j^2}{\N{}^2}  \int_{0}^{+\infty}\sqrt{v}(I_s+I_b)(v-j^2\pi\theta/(2\N{}d))dv,\]
\item for $\N{}^{1/2}\ll j\lesssim \nu \N{}$, with $\nu\in(0,1)$
\[C_j(\omega)\underset{k(\omega)\gg1}{\sim} -  C_3 \, k^2(\omega)\frac{j^3}{\N{}^3} \int_{0}^{+\infty}(I_s+I_b)(v)dv,\]
\item for $j\sim \N{}$
\[ C_4 \, k(\omega)^{5/2} \leq C_j(\omega)\leq C_5 \, k(\omega)^{5/2}\] 
\end{itemize}
when $k(\omega)\gg 1$, where $C_1,C_2,C_3,C_4,C_5>0$, and $\theta$ is defined by \eqref{deftheta}.
\end{thm}
In Theorem \ref{thmain} we give the order of magnitude of the mode-dependent and frequency-dependent decay rates $C_j(\omega)$ responsible for the effective attenuation of the forward propagating-mode amplitudes. The first remark is that the surface and bottom fluctuations affect the mode amplitudes in the same way. The second remark is that the higher the propagating mode is the stronger is the attenuation. In fact, the higher the mode is the more it bounces on the random boundaries, and then the more it is scattered. However, as described more precisely in Proposition \ref{coefdec} and Proposition \ref{coefdec2}, while $j$ is not of order $\N{}$ the phenomena is mainly due to the mode coupling between the propagating modes themselveses. However, for the highest order modes $j\sim \N{}$, this attenuation is mainly produced by the coupling with the radiating modes. The reason is that only the highest order modes can couple efficiently with the radiating modes (see \eqref{coefatt}, \eqref{coefatt2}, and \eqref{spectrumRP2}) and therefore produced a stronger attenuation since the radiating modes induce strong loss in the bottom of the waveguide (see \cite{gomez2}). 

The remaining of this paper consists in introducing the tools used in this paper and the simplifying assumptions (Assumptions \ref{hyp1}-\ref{hyp5}) to prove Theorem \ref{thmain} through Theorem \ref{thasymptP21}, Theorem \ref{thasymptP12}, Theorem \ref{thasympP2}, Proposition  \ref{coefdec}, and Proposition  \ref{coefdec2}.

\section{Conformal Transformation of the Waveguide}\label{conformal}

The free surface of our waveguide model (see Figure \ref{figureP2}) is not convenient for the forthcoming asymptotic analysis based on a spectral decomposition of the acoustic  pressure wave. In this paper, we consider a conformal transformation allowing us to transform our waveguide model with a free surface to a waveguide with an unperturbed surface (see Figure \ref{conftrans1}). This strategy has already been used in several contexts \cite{dozier, dya,garniernach, kuttler}. We will see in Section \ref{sect3P2} that this transformation by flattening the surface induces deformations in the structure of the waveguide, that is deformations of the bottom and also induces variations of the index of refraction of the waveguide (see Figure \ref{figuredef}). 

To introduce the conformal transformation we consider the waveguide domain with a free surface as a subdomain of the complex plan $\mathbb{C}$:
\[D_{V_s}=\Big\{z+ix \in\mathbb{C}\,:\quad x >V^\e_s(z ) \Big\},\] 
where $V^\e_s$ is defined by \eqref{randompert} (see Figure \ref{conftrans1}), and the domain corresponding to the waveguide with a flat surface is then given by
\[D_0=\Big\{u+iv \in\mathbb{C}\,:\quad v>0 \Big\}.\] 
We refer to \cite{rudin} for the basic properties of conformal transformations. In order to provide an explicit expression of this transformation (from $D_0$ onto $D_{V_s}$) we need to consider boundary conditions at infinity. We have assumed that the random perturbations of the waveguide geometry are included in the interval $z\in(0,L/\e)$, then the difference between the conformal map and the identity map goes to $0$ as $Z=u+iv$ goes to $+\infty$. The following proposition, which is a consequence of \cite[Theorem 14.8]{rudin}, proves the existence and gives explicit formula of a conformal transformation from $D_0$ onto $D_{V_s}$.  

\begin{prop}\label{propconfmap}
There exists a bijective conformal map $\Phi$  from $D_0$ onto $D_{V_s}$ defined by
\[\Phi(u+iv)=z(u,v)+ix(u,v),\]
with
\begin{equation}\label{conformmap}z(u,v)=u-\frac{1}{\pi}\int_{-\infty}^{+\infty}d\tilde{u}\frac{(\tilde{u}-u)V^\e_s(z(\tilde{u},0))}{(\tilde{u}-u)^2+v^2}
\quad \text{and}\quad x(u,v)=v-\frac{v}{\pi}\int_{-\infty}^{+\infty}d\tilde{u}\frac{V^\e_s(z(\tilde{u},0))}{(\tilde{u}-u)^2+v^2},
\end{equation}
for all $(u,v)$ such that $u+iv\in D_0$. 
\end{prop}
Let us note that we can extend the conformal transformation $\Phi$ as a homeomorphism from the closure of $D_0$ onto the closure of $D_{V_s}$ thanks to \cite[Theorem 14.19]{rudin}, so that the derivation of \eqref{conformmap} is a classical result about harmonic functions which can be found in \cite{axler,nehari} for instance.

Consequently, one can consider the pressure wave in the new system of coordinates:
\[p_0(t,u,v)=p(t,z(u,v),x(u,v)),\quad \forall (t,u,v)\in[0,+\infty)\times(-\infty,+\infty)\times(0,+\infty),\]
where $p$ is the solution of \eqref{waveeq}.
Therefore, $p_0(t,u,v)$ evolves in the unperturbed domain $D_0$ and satisfies Dirichlet boundary conditions at the flat surface of the waveguide:
\[p_0(t,u,0)=0, \quad \forall (t,u)\in[0,+\infty)\times(-\infty,+\infty).\]
Moreover, the original pressure wave $p(t,z,x)$ can be recovered from $p_0(t,u,v)$ by inverting the change of coordinates as follows,
\[ p(t,z,x) = p_0(t,Re(\Phi^{-1}(z+ix)), Im(\Phi^{-1}(z+ix)) ).\]
However, \eqref{conformmap} involves the term $z(u,0)$ which is not convenient for the forthcoming analysis. In fact, at the boundary $v=0$ it is difficult to use the inductive formula $z(u,0)=u-\mathcal{U}(z(.,0))(u)$, where $\mathcal{U}$ is the Hilbert transform defined in Proposition \ref{propconfmap2}. In the following proposition we give an approximation of $z(u,0)$ and then an approximation of $x(u,0)$ to motivate Assumption \ref{hyp1} introduced below.
\begin{prop}\label{propconfmap2}
We have for every $p\in(2,+\infty)$
\begin{equation}\label{zo}
z(u,0)=u-\mathcal{U}(V^{\e}_s)(u)+\e^{1-1/p} A^\e_1(u),
\end{equation}
and then
\begin{equation}\label{xo}
x(u,0)=V^\e_s(z(u,0))=V^\e_s(u)- \mathcal{U}(V^{\e}_s)(u)\frac{d}{du}V^\e_s(u) +\e^{3/2-1/p} A^\e_2(u),
\end{equation}
with $\sup_\e \| A^\e_1\|_{L^p(\mathbb{R})}+\| A^\e_2\|_{L^p(\mathbb{R})}\leq C$ almost surely, and where $C$ is a deterministic constant. Moreover, in \eqref{zo} and \eqref{xo}, $\mathcal{U}$ stands for the Hilbert transform defined by
\[\mathcal{U}(f)(u)=\frac{1}{\pi}\emph{\text{p.v.}} \int \frac{f(\tilde{u})}{u-\tilde{u}}d\tilde{u},\]
where p.v. denotes the Cauchy principal value.
\end{prop}
We detail how to obtain \eqref{zo} in Section \ref{proofpropconfmap}, and \eqref{xo} comes from \eqref{zo} in addition to a Taylor expansion of $V^\e_s$ defined by \eqref{randompert}. The term $\e^{3/2-1/p}$ in \eqref{xo} comes from the product between $dV^\e_s/du$ and $\e^{1-1/p} A^\e_1(u)$, and all the remaining terms in the Taylor expansion have a smaller order of magnitude. Let us remark that we prove this approximation only in $L^p(\mathbb{R})$ with $p<\infty$ for the reason that the Hilbert transform is not a bounded operator in $L^{\infty}(\mathbb{R})$, and we chose $p>2$ so that $A^\e_1(u)$ in \eqref{zo} and $ A^\e_2(u)$ in \eqref{xo} do not play any significant role in the forthcoming analysis. In particular we need to have $1-1/p>1/2$ ($p>2$) in \eqref{zo} and $3/2-1/p>1$ in \eqref{xo}.    

Consequently, for the sake of simplicity in the proof of Theorem \ref{thasymptP21} (Section \ref{proofth}) let us consider the following assumption. 
\begin{assumption}\label{hyp1}
\begin{equation}\label{zo2}z(u,0)=u-\mathcal{U}(V^{\e}_s)(u),\end{equation}
and
\begin{equation}\label{xo2}x(u,0)=V^\e_s(u)-\mathcal{U}(V^{\e}_s)(u)\frac{d}{du}V^\e_s(u).\end{equation}
\end{assumption}
This assumption give us an approximation of the conformal transformation $\Phi$ introduced in Proposition \ref{propconfmap} by neglecting the contributions of $\e^{1-1/p} A^\e_1(u)$ and $\e^{3/2-1/p} A^\e_2(u)$. However, this assumption will not change the overall results of this paper because the contributions of the correctors $\e^{1-1/p} A^\e_1(u)$ (resp. $\e^{3/2-1/p} A^\e_2(u)$) are too small with respect to $\se$ (resp. $\e$) to produce any significant effects. In Assumption \ref{hyp1} we keep only the terms which produce significant effects.   

We can remark that \eqref{zo2} illustrates the fact that $z(u,0)$ is closed to the identity map (since $V^\e_s$ is defined by \eqref{randompert}) but where the small correction comes from the conformal transform. The same remark holds for \eqref{xo2}, at the boundary of the domain $D_{0}$, $x(u,0)$ is closed to the original free surface $V^\e_s$ but also with a small correction.  

Finally, we describe in the following proposition the asymptotic behavior of the real and imaginary part of the conformal transformation given by \eqref{conformmap} under Assumption \ref{hyp1}.
\begin{prop}\label{asymptconf}
For all $p>2$, we have 
\[\sup_{u, \,v\geq 0}\lvert z(u,v)-u \rvert + \lvert x(u,v)-v \rvert\leq C\e^{1/2-1/p}\quad \text{almost surely},\]
and
\[\sup_{u}\int_0^{+\infty}\E\Big[  \lvert z(u,v)-u \rvert^2 + \lvert x(u,v)-v \rvert^2 \Big] dv \leq C\e^{1/2-1/p}, \]
where $C$ is a deterministic constant independent of $\e$.
\end{prop}
Proposition \ref{asymptconf} means, as expected, that the conformal transformation is closed to the identity map since the amplitude of the surface fluctuations are small. The proof of Proposition \ref{asymptconf} is given in Section \ref{proofasymptconf}.

\section{Mode Coupling in Random Waveguides} \label{sect3P2}

In this paper, the stochastic effects produced on the wave propagation is studied using the spectral decomposition of the unperturbed Pekeris operator defined by \eqref{pekerisop} and introduced in Section \ref{spectralP2}. In this section we introduce the mode coupling mechanism induces by the random free surface and the uneven bottom topography.

The propagation model \eqref{waveeq} considered in this paper is a linear model, so that the pressure wave $p(t,z,x)$ can be expressed as the superposition of monochromatic waves by taking its Fourier transform \eqref{fourier}, and then satisfies the time-harmonic wave equation
\[  \partial^2_z \widehat{p}(\omega,z,x)+\partial^2_x \widehat{p}(\omega,z,x)+\frac{\omega^2}{c^2_\e(z,x)}\widehat{p}(\omega,z,x)=\nabla\cdot \widehat{F}(\omega,z,x),\]
where $c_\e$ is defined by \eqref{propspeed}, and $\widehat{F}(\omega,z,x)$ is the Fourier transform of the source term $F(t,z,x)$ defined by \eqref{sourcetermP2}.

Let us first describe the mode decomposition of the pressure field entering the perturbed section of the waveguide $(0,L/\e)$. According to the spectral decomposition introduced in Section \ref{spectralP2} and described in Section \ref{sect1P2}, the monochromatic pressure field can be expanded as follows in the unperturbed section $[L_S,0]$,   
\[\widehat{p}(\omega,z,x)=\underbrace{\sum_{j=1}^{\N{}} \widehat{p}_j (\omega,z)\phi_j(\omega,x)}_{\text{propagating modes}}+\underbrace{\int_{0}^{\ko{}}\widehat{p}_\ga (\omega,z)\phi_\ga (\omega,x)d \ga}_{\text{radiating modes}}+\underbrace{\int_{-\infty}^{0}\widehat{p}_\ga (\omega,z)\phi_\ga (\omega,x)d \ga}_{\text{evanescent modes}},\]
where the mode amplitudes are defined by \eqref{defmode1} and \eqref{defmode2}, and where $\phi_j(\omega,u)$ ($j\in \big\{1,\dots,\N{}\big\}$) and $\phi_\ga(\omega,u)$ ($\forall \ga\in(-\infty,\ko{})$), defined by \eqref{phij} and \eqref{phiga}, are the eigenelements of the Pekeris operator. In this paper we consider a source profile $\Psi$ in \eqref{sourcetermP2} admitting in the Fourier domain the following decomposition,
\begin{equation}\label{profsourceP2}
 \widehat{\Psi}(\omega,x)=\widehat{f}(\omega)\left[\sum_{j=1}^{\N{}}\phi_j(\omega,x_0)\phi_j(\omega,x)+\int_{(-1/\xi,-\xi)\cup(\xi,\ko{})}\phi_\ga(\omega,x_0)\phi_\ga(\omega,x)d \ga\right],
 \end{equation}
where $x_0\in (0,d)$. 
The meaning of the bound $\xi$ is that the source does not excite the high evanescent modes nor the low radiating and evanescent modes. However, $\xi$ can be arbitrarily small, so that $\Psi(t,x)\simeq f(t)\delta(x-x_0)$ models a source profile closed to a point source emitting a signal $f(t)$. Let us note that $\xi$ is introduced for technical reasons as explain in Section \ref{sect1P2} and Section \ref{couplemodeeq}. The source term \eqref {sourcetermP2} implies the following jump conditions for the pressure field across the plane $z=L_S$ \cite{book,gomez2}, 
\begin{equation}\label{jumpscondP2}\begin{array}{ccl}
\widehat{p}(\omega,x,L_S^+)-\widehat{p}(\omega,x,L_S^-)&=& \widehat{\Psi}(\omega,x),\\
\partial_z \widehat{p}(\omega,x,L_S^+)-\partial_z \widehat{p}(\omega,x,L_S^-)&=&0,
\end{array}
\end{equation}
so that the pressure field $\widehat{p}(\omega,z,x)$ for $z\in(L_S,0]$ is given by
\begin{equation}\label{coefP2}
\begin{split}
\widehat{p}(\omega,z,x) &=\sum_{j=1}^{\N{}}  \frac{\widehat{a}^\xi_{j,0} (\omega)}{\sqrt{\Bh{j}{}}}e^{i\Bh{j}{}z} \phi_j(\omega,x)+\int_{\xi}^{\ko{}}  \frac{\widehat{a}^\xi_{\ga,0} (\omega)}{\ga^{1/4}} e^{i \sga z}\phi_\ga(\omega,x)d \ga\\
&\hspace{1.5cm}+\int_{-1/\xi}^{-\xi}\frac{\widehat{c}^\xi_{\ga,0} (\omega)}{\gaa^{1/4}}  e^{-\sgaa z}\phi_\ga(\omega,x)d\ga, \quad \forall x>0,
\end{split}
\end{equation}
and where
\begin{equation}\label{condinith}
\begin{split}
\widehat{a}^\xi_{j,0} (\omega)&=\frac{\sqrt{\Bh{j}{}}}{2}\widehat{f}(\omega)\phi_j(\omega,x_0)e^{-i\Bh{j}{}L_S}\quad j\in\{1,\dots,\N{}\},\\
\widehat{a}^\xi_{\ga,0} (\omega)&=\frac{\ga^{1/4}}{2}\widehat{f}(\omega)\phi_\ga(\omega,x_0)e^{-i\sga L_S}\quad \ga\in(\xi,\ko{}),\\
\widehat{c}^\xi_{\ga,0} (\omega)&=\frac{\gaa^{1/4}}{2}\widehat{f}(\omega)\phi_\ga(\omega,x_0)e^{\sgaa L_S}\quad \ga\in(-1/\xi,-\xi).
\end{split}\end{equation}
Now, to study the evolution of the pressure field in the perturbed section $(0,L/\e)$ through a mode coupling mechanism, we  consider the new system of coordinates introduced in Section \ref{conformal}. This transformation allows us to flatten the free surface of the waveguide domain, and then to decompose the pressure field with respect to the propagating, radiating, and evanescent modes. 

\begin{figure}\begin{center}
\begin{tabular}{cc}
\includegraphics*[scale=0.3]{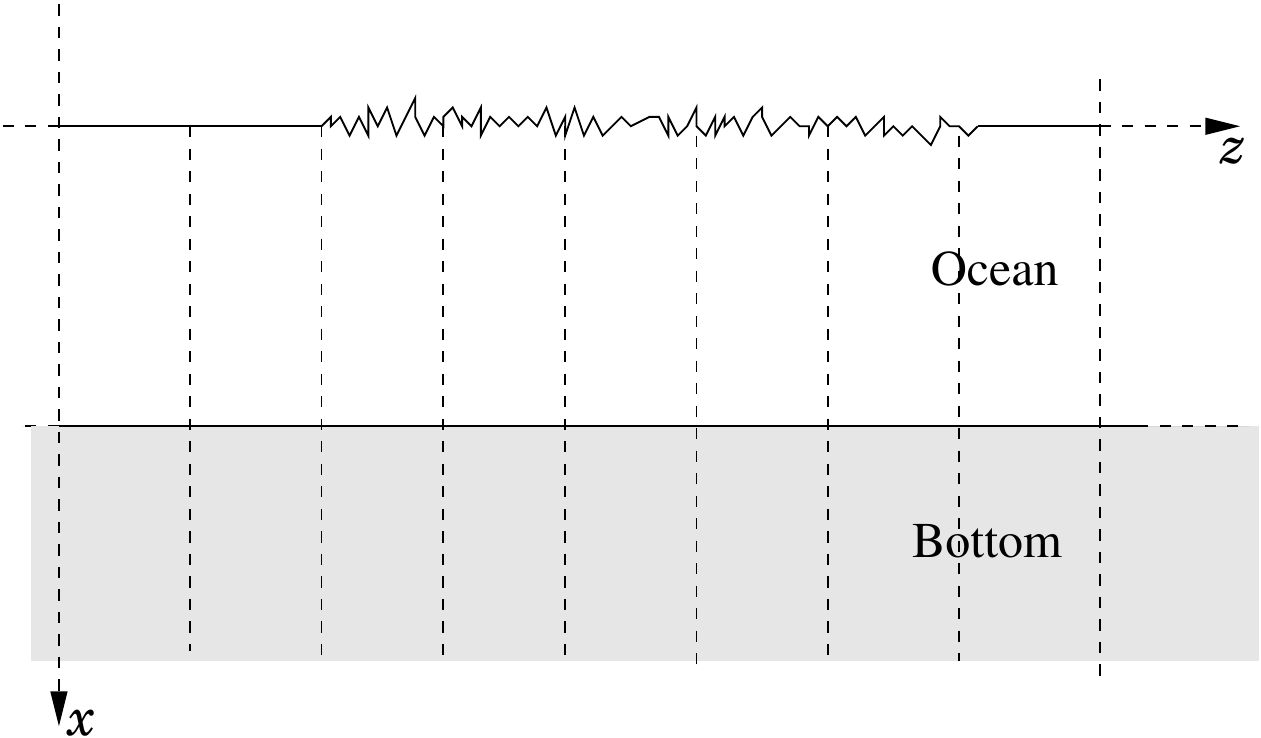} & \includegraphics*[scale=0.3]{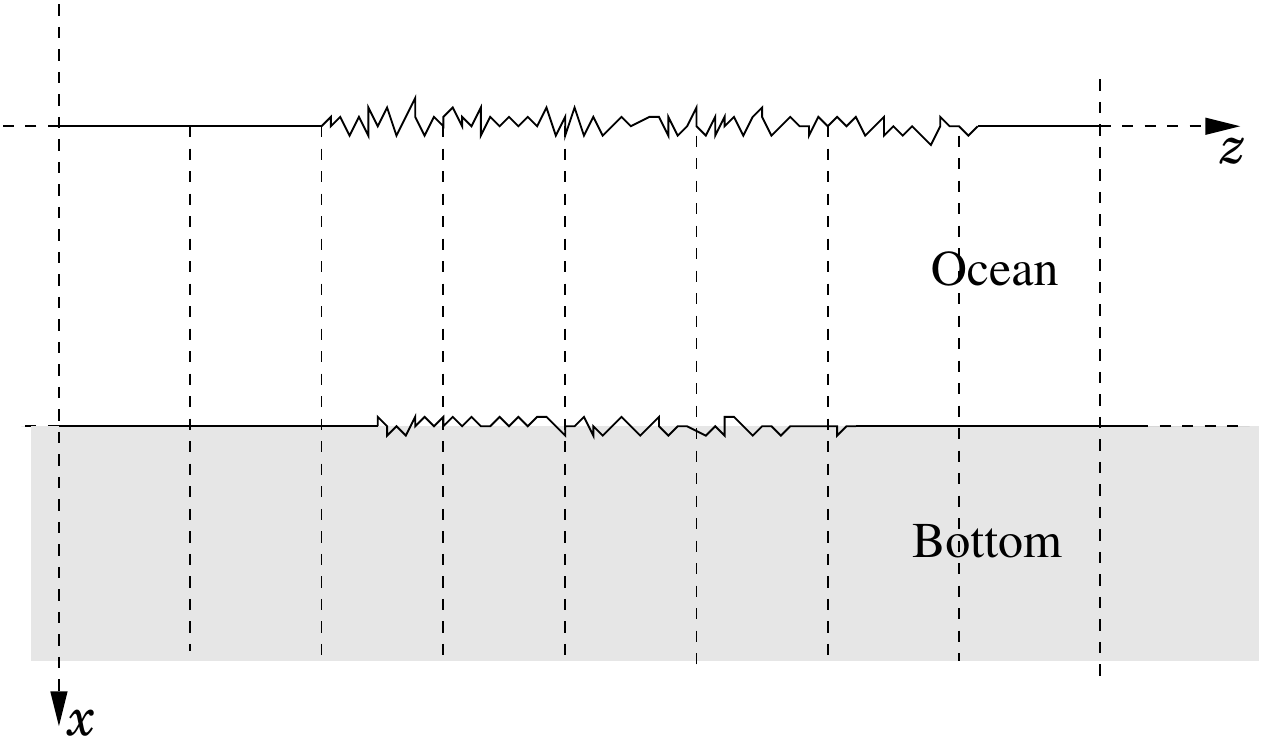}\\
$(a)$ & $(b)$\\
\includegraphics*[scale=0.3]{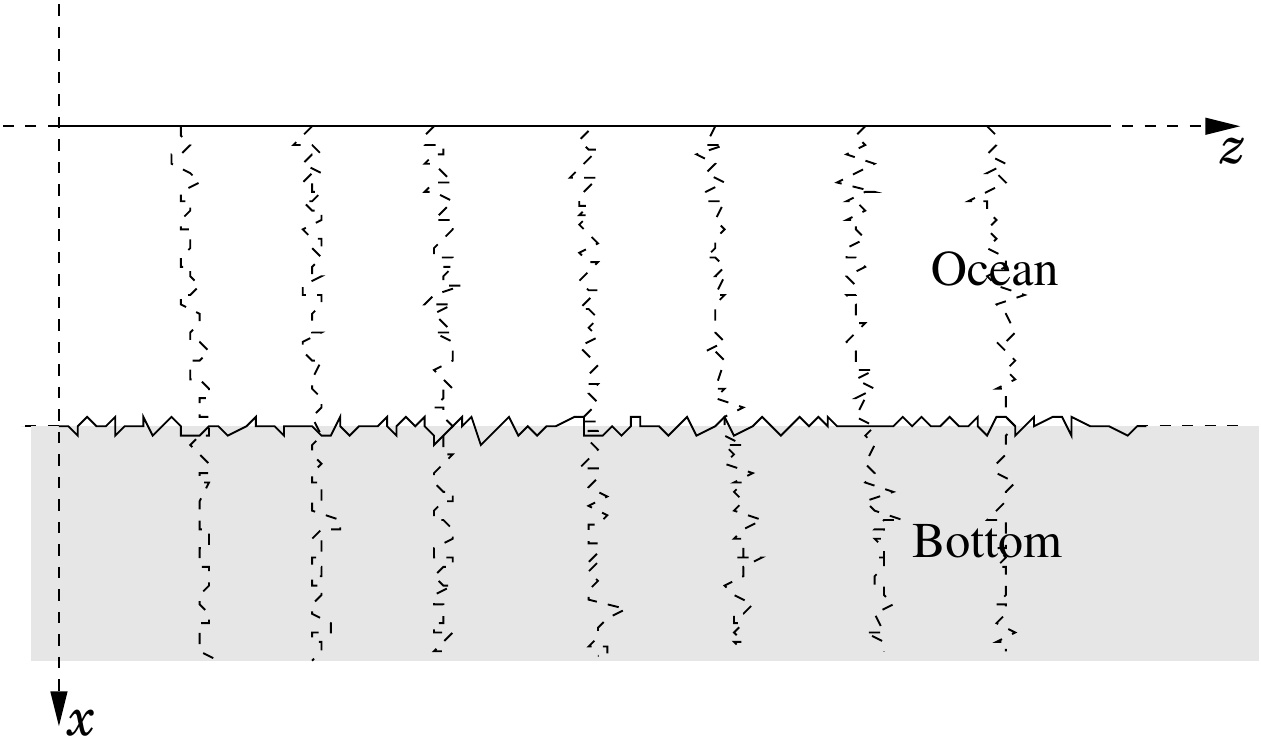} & \\
$(c)$ & 
\end{tabular}
\end{center}
\caption{\label{figuredef}Illustration of the conformal transformation on the waveguide with a free surface and an uneven bottom topography. Figure $(a)$ and Figure $(b)$ represent waveguides with a free surface. In these figures the vertical dashed lines mean that the medium parameters are not perturbed. In $(a)$ the bottom is not perturbed, that is $V^\e _b=0$, and in $(b)$ the bottom of the waveguide is perturbed. Figure $(c)$ represents the waveguide after the conformal transformation.  In this figure the perturbed vertical dashed lines mean that the medium parameters are perturbed. Moreover, if the bottom was already perturbed before the conformal transformation (Figure $(b)$), it is now also is perturbed by the waveguide transformation.} \end{figure}

The pressure field $\widehat{p}(\omega,z,x)$ in the new system of coordinates is given by
\begin{equation}\label{p0}\widehat{p}_0(\omega,u,v)=\widehat{p}(\omega,z(u,v),x(u,v)),\quad \forall (u,v)\in(-\infty,+\infty)\times(0,+\infty),\end{equation}
where $z(u,v)$ and $x(u,v)$ are defined by \eqref{conformmap} under Assumption \ref{hyp1}. As a result, thanks to Proposition \ref{asymptconf}, the monochromatic pressure field $\widehat{p}_0(\omega,u,v)$ satisfies the following time-harmonic wave equation
\begin{equation}\label{hwaveeq}  \partial^2_u \widehat{p}_0(\omega,u,v)+\partial^2_v \widehat{p}_0(\omega,u,v) +k^2(\omega)J_\e(u,v)n^2_\e(u,v) \widehat{p}_0(\omega,u,v)=0,\end{equation}
for all $(u,v)\in [0,L/\e]\times [0,+\infty)$, with Dirichlet boundary conditions
\begin{equation}\label{dbc}\widehat{p}_0(\omega,u,0)=0, \quad \forall u\in (-\infty,+\infty).\end{equation}
Moreover, in \eqref{hwaveeq}, the Jacobian $J_\e$ and the index of refraction $n_\e$ are given by
\begin{equation}\label{jacodef}J_\e(u,v)=\Big(\partial_u x(u,v)\Big)^2+\Big(\partial_v x(u,v)\Big)^2, \end{equation}
and for all $(u,v)\in (-\infty,+\infty)\times(0,+\infty)$
\begin{equation}\label{indexrefraction} n_\e(u,v)=c_0/c_\e(z(u,v),x(u,v))=\left\{ \begin{array}{ccl} 
                                            n_1=c_0/c_1>1 & \text{ if }  & v\in(0,x^{-1}(u,d+ V^\e_b (z(u,v)))), \\
                                             1 & \text{ if }  & v \in [x^{-1}(u,d+ V^\e_b (z(u,v))),+\infty).\\
                                          \end{array} \right.
\end{equation}
According to \eqref{hwaveeq}, let us remark that the waveguide transformation transfers the surface perturbations to the interior of the waveguide through the Jacobian $J_\e$, and to the bottom of the waveguide through the index of refraction $n_\e(u,v)$ defined by \eqref{indexrefraction} (see Figure \ref{figuredef}). Thanks to the Dirichlet boundary conditions \eqref{dbc}, the effects of the random perturbations can be studied using the modal decomposition corresponding to an unperturbed waveguide (see Figure \ref{wavegmod} $(a)$). 

\subsection{Coupled Mode Equations}\label{couplemodeeq}

Let us begin this section with two remarks. First, $J_\e(u,v)n^2_\e(u,v) \widehat{p}_0(\omega,u,v)$ does not necessarily belong to $H=L^2(0,+\infty)$ since the Jacobian $J_\e$ blows up as $v$ goes to $0$. To correct this problem we need to introduce the following assumption on the modal decomposition \eqref{modedec} of the pressure field $\widehat{p}_0(\omega,u,v)$ defined by \eqref{p0}. 
\begin{assumption} \label{hyp2}
\begin{equation}\label{modechyp}\widehat{p}_0(\omega,u,v)=\sum_{j=1}^{\N{}} \widehat{p}_j (\omega,u )\phi_j(\omega, v)+\int_{(-1/\xi,-\xi)\cup(\xi,\ko{})}\widehat{p}_\ga (\omega,u)\phi_\ga (\omega,v)d \ga,\end{equation}
\end{assumption}
where the modes amplitude are defined by \eqref{defmode1} and \eqref{defmode2}. The cut-off induced by $\xi$ in the previous decomposition allows to obtain rigorous derivation of the following coupled mode equations more easily, but do not change the overall results obtained in Section \ref{thasyptoticP2}. Consequently, according to Assumption \ref{hyp2} the coupling mechanism of the spectral components of the pressure field will be studied in the space
\[ \mathbb{C}^{\N{}}\times L^2\big((-1/\xi,-\xi)\cup(\xi,\ko{})\big).\]
Finally, we assume that $\e\ll\xi\ll1$ so that we have two distinct scales. In this paper, we consider first the asymptotic $\e$ goes to $0$ and in a second time the asymptotic $\xi$ goes to $0$.

Consequently, the mode amplitudes in \eqref{modechyp} satisfy the following coupled mode equations for $u\in[0,L/\e]$:
\begin{equation}\label{eqdiff2P2}\begin{split}
\frac{d^2}{du^2}\widehat{p}_j(\omega,u)+\beta^2 _j(\omega) \widehat{p}_j(\omega,u)&+\ko{}\sum_{l=1}^{\N{}}C^\e_{jl}(\omega,u)\widehat{p}_l(\omega,u)\\
            &+\ko{} \int_{(-1/\xi,-\xi)\cup(\xi,\ko{})}C^\e_{j\ga'}(\omega,u)\widehat{p}_{\ga'}(\omega,u)d\ga'=0,\\
\frac{d^2}{du^2}\widehat{p}_{\ga}(\omega,u)+\ga \,\,\widehat{p}_{\ga}(\omega,u)&+\ko{}\sum_{l=1}^{\N{}}C^\e_{\ga l}(\omega, u)\widehat{p}_{l}(\omega,u)\\
            &+\ko{} \int_{(-1/\xi,-\xi)\cup(\xi,\ko{})} C^\e_{\ga \ga'}(\omega, u)\widehat{p}_{\ga'}(\omega, u)d\ga'=0,          
\end{split}\end{equation}
where
\begin{equation}\label{randomcoef}
C^\e_{rs}(\omega,u)=\int_0^{+\infty} (J_\e(u,v)n^2_\e(u,v)-n^2(v))  \phi_r(\omega,v)\phi_s(\omega,v)dv,
\end{equation}
for $(r,s)\in(\{1,\dots,\N{}\}\times(-1/\xi,-\xi)\cup(\xi,\ko{}))^2$, and
\begin{equation}\label{indexrefraction2} n(v)=\left\{ \begin{array}{ccl} 
                                            n_1=c_0/c_1>1 & \text{ if }  & v\in(0,d), \\
                                             1 & \text{ if }  & v \in [d,+\infty).\\
                                          \end{array} \right.\end{equation}
Under  Assumption \ref{hyp2}, the coupling coefficients $C^\e(\omega,u)$ and then the previous coupled mode equations in $\mathbb{C}^{\N{}}\times L^2\big((-1/\xi,-\xi)\cup(\xi,\ko{})\big)$ are well defined since $\exists K_{1,\e}>0$ such that $\sup_u \lvert J_\e(u,v)-1 \rvert\leq K_{1,\e}/v^2$, and 
\begin{equation}\label{est}\phi_r(\omega,v)\underset{v\to0}{\sim}{K_2}v\quad\forall r\in\{1,\dots,\N{}\}\times(-1/\xi,-\xi)\cup(\xi,\ko{}).\end{equation}
Moreover, $\sup_u\lvert n^2_\e(u,v)-n^2(v)\rvert$ is bounded and has a bounded support with respect to the transverse variable $v$.

Let us remark that 
\[J_\e(u,v)n^2_\e(u,v)-n^2(v)=n^2_\e(u,v)(J_\e(u,v)-1)+(n^2_\e(u,v)-n^2(v)),\]
so that the first part of the decomposition induces random perturbations in the core of the waveguide, while the second part which is the difference of the index of refractions induces random perturbations at the bottom of the waveguide (see Figure \ref{figuredef}). 
 
 The two asymptotic results in Section \ref{coupledprocP2}  (Theorem \ref{thasymptP21} and Theorem \ref{thasymptP12}) are based on a diffusion-approximation result for the solution of ordinary differential equations with random coefficients in a Hilbert space. However, according to the formula of the conformal transformation \eqref{conformmap} and the coefficients \eqref{randomcoef}, we cannot directly apply the asymptotic results obtained in \cite{book,gomez2,gomez3} since the mode amplitudes are not adapted to the $\phi$-mixing filtration \eqref{filtration}. In fact, in \eqref{conformmap} all the trajectory of $V^\e_s$ is involved, which leads us to technical difficulties to manage the random coefficients \eqref{randomcoef}, and to use the mixing properties of the random perturbations \eqref{mixcond}. To avoid this technical problem we consider the following assumption. 
\begin{assumption}\label{hyp3}
\[C^\e_{\ga\ga'}(\omega,u)=0, \quad \forall (\ga,\ga')\in (-1/\xi,-\xi)\cup(\xi,k^2(\omega)).\]
\end{assumption}
This assumption means that we neglect the coupling mechanism between the radiating and evanescent modes, they do not interact with each others. Let us remark that in \cite{gomez2}, in which the random perturbations take place in the interior of the waveguide, the mode coupling mechanism between these modes do not play any significant role.

Next, we introduce the amplitudes of the  generalized right- and left-going modes 
$\widehat{a}(\omega,z)$ and $\widehat{b}(\omega,z)$, which are given by
\[\begin{split}
\widehat{p}_j(u)&=\frac{1}{\sqrt{\beta_j}}\Big( \widehat{a}_j(u) e^{i\beta_j u} +\widehat{b}_j(u)e^{-i\beta_j u} \Big), \quad \frac{ d }{du}\widehat{p}_j(u) = i \sqrt{\beta_j} \Big( \widehat{a}_j(u) e^{i\beta_j u} - \widehat{b}_j(u)e^{-i \beta_j u} \Big),\\
\widehat{p}_\ga(u)&= \frac{1}{\ga^{1/4}}\Big( \widehat{a}_\ga(u)e^{i\sga u} +\widehat{b}_\ga(u)e^{-i\sga u} \Big),\quad \frac{ d }{du}\widehat{p}_\ga(u) = i \ga^{1/4} \Big( \widehat{a}_\ga(u) e^{i\sga u} - \widehat{b}_\ga(u)e^{-i\sga u} \Big),
\end{split}\]
so that 
\begin{equation}\label{ab}\begin{split}
\widehat{a} _{j}(u)&=\frac{i\beta_j \widehat{p}_j(u)+\frac{d}{dz} \widehat{p}_j(u)}{2i \sqrt{\beta_j}}e^{-i \beta_j u},\quad \widehat{a} _{\ga}( u)=\frac{i\sga \widehat{p}_\ga(u)+\frac{d}{dz} \widehat{p}_\ga(u)}{2i \ga^{1/4}}e^{-i\sga u},\\
\widehat{b} _{j}(u)&=\frac{i\beta_j \widehat{p}_j(u)-\frac{d}{dz} \widehat{p}_j(u)}{2i \sqrt{\beta_j}}e^{i\beta_j u},\quad \widehat{b} _{\ga}(u)=\frac{i\sga \widehat{p}_\ga(u)-\frac{d}{dz} \widehat{p}_\ga(u)}{2i \ga^{1/4}}e^{i\sga u},
\end{split}\end{equation}
for all $j \in \big\{1,\dots,\N{}\big\}$ and almost every $\ga \in (\xi,\ko{})$. 
From now on, let us denote
\[\esp_\xi=\mathbb{C}^{\N{}}\times L^2(\xi,\ko{}).\]
From \eqref{eqdiff2P2} and \eqref{ab}, we obtain the coupled mode equation in $\esp_\xi \times \esp_\xi\times L^2(-1/\xi,-\xi)$ for the mode amplitudes $\big(\widehat{a}(\omega,u),\widehat{b}(\omega,u),\widehat{p}(\omega,u)\big)$:
\begin{equation}\label{mcP21} \begin{split} 
\frac{ d }{du} \widehat{a} _{j}(u)&=\textbf{H}^{aa}_{\e,\xi,j}(u)(\widehat{a}(u))+\textbf{H}^{ab}_{\e,\xi,j}(u)(\widehat{b}(u))+\frac{ik^2}{2}\int_{-1/\xi}^{-\xi} \frac{C^\e_{j\ga'}(u)}{\sqrt{\beta_j}}\widehat{p}_{\ga'}(u)d\ga' e^{-i \beta_j  u},\\
\frac{ d }{du} \widehat{a} _{\ga}(u)&=\textbf{H}^{aa}_{\e,\xi,\ga}(u)(\widehat{a}(u))+\textbf{H}^{ab}_{\e,\xi,\ga}(u)(\widehat{b}(u)),
\end{split}\end{equation}
\begin{equation}\label{mcP24}\begin{split}
\frac{ d }{du} \widehat{b} _{j}(u)&=\overline{\textbf{H}^{ab}_{\e,\xi,j}(u)(\overline{\widehat{a}(u)})}+\overline{\textbf{H}^{aa}_{\e,\xi,j}(u)(\overline{\widehat{b}(u)})}-\frac{ik^2}{2}\int_{-1/\xi}^{-\xi} \frac{C^\e_{j\ga'}(u)}{\sqrt{\beta_j}}\widehat{p}_{\ga'}(u)d\ga' e^{i \beta_j  u},\\
\frac{ d }{du} \widehat{b} _{\ga}(u)&=\overline{\textbf{H}^{ab}_{\e,\xi,\ga}(u)(\overline{\widehat{a}(u)})}+\overline{\textbf{H}^{aa}_{\e,\xi,\ga}(u)(\overline{\widehat{b}(u)})},
\end{split}\end{equation}
\begin{equation}\label{evaeqdiff}
\frac{d^2}{du^2}\widehat{p}_{\ga}(\omega,u)+\ga \,\,\widehat{p}_{\ga}(\omega,u)+\se g_\gamma (\omega, u)=0,           
\end{equation}
where
\begin{equation}\label{g}\begin{split} g_\gamma ( u)&= k^2\sum_{l=1}^{N}\frac{C^\e_{\ga l}(u)}{\sqrt{\beta_l}}\left(   \widehat{a} _{l}(u)e^{i \beta_l u} + \widehat{b} _{l}(u) e^{-i \beta_l u} \right)\\
&+k^2 \int_{\xi}^{k^2} \frac{C^\e_{\ga \ga'}(u)}{{\ga'}^{1/4}}\left(   \widehat{a} _{\ga'}(u)e^{i \sgap u} + \widehat{b} _{\ga'}(u) e^{-i \sgap u} \right)d\ga'+k^2\int_{-1/\xi}^{-\xi}C^\e_{\ga \ga'}(u)\widehat{p}_{\ga'}(u)d\ga'.\end{split}\end{equation}
In \eqref{mcP21}-\eqref{mcP24} we have
\begin{equation}\label{haajP2}
\textbf{H}^{aa}_{\e,\xi,j}(u)(y)=\frac{ik^2}{2}\Big[\sum_{l=1}^{N} \frac{C^\e_{jl}(u)}{\sqrt{\beta_j \beta_l}} e^{i(\beta_l -\beta_j)u}y_l+ \int_{\xi}^{k^2}\frac{C^\e_{j\ga'}(u)}{\sqrt{\beta_j \sgap}}e^{i(\sgap -\beta_j)u}y_{\ga'}d\ga' \Big],
\end{equation}
\begin{equation}\label{haagaP2}
\textbf{H}^{aa}_{\e,\xi,\ga}(u)(y)=\frac{ik^2}{2}\Big[\sum_{l=1}^{N} \frac{C^\e_{\ga l}(u)}{\sqrt{\sga \beta_l}} e^{i(\beta_l -\sga)u}y_l+ \int_{\xi}^{k^2}\frac{C^\e_{\ga \ga'}(u)}{\ga^{1/4}{\ga'}^{1/4}}e^{i(\sgap -\sga)u} y_{\ga'} d\ga' \Big],
\end{equation}
\begin{equation}\label{habjP2}
\textbf{H}^{ab}_{\e,\xi,j}(u)(y)=\frac{ik^2}{2}\Big[\sum_{l=1}^{N} \frac{C^\e_{jl}(u)}{\sqrt{\beta_j \beta_l}} e^{-i(\beta_l +\beta_j)u}  y_l+ \int_{\xi}^{k^2}\frac{C^\e_{j\ga'}(u)}{\sqrt{\beta_j \sgap}}e^{-i(\sgap + \beta_j)u} y_{\ga'} d\ga' \Big],
\end{equation}
\begin{equation}\label{habgaP2}
\textbf{H}^{ab}_{\e,\xi,\ga}(u)(y)\frac{ik^2}{2}\Big[\sum_{l=1}^{N} \frac{C^\e_{\ga l}(u)}{\sqrt{\sga \beta_l}} e^{-i(\beta_l +\sga)u}y_l+ \int_{\xi}^{k^2}\frac{C^\e_{\ga \ga'}(u)}{\ga^{1/4}{\ga'}^{1/4}}e^{-i(\sgap + \sga)u}y_{\ga'} d\ga' \Big],
\end{equation}
The operators $\textbf{H}^{aa}_{\e,\xi}(\omega,u)$ and $\textbf{H}^{ab}_{\e,\xi}(\omega,u)$ represent the coupling between the propagating and the radiating modes with themselves. Moreover, $\textbf{H}^{aa}_{\e,\xi}(\omega,u)$ describes the coupling between the forward-going modes, while $\textbf{H}^{ab}_{\e,\xi}(\omega,u)$ describes the coupling between the forward- and backward-going modes. Let us note that in absence of random perturbations, the mode amplitudes $\widehat{a}(\omega,u)$ and $\widehat{b}(\omega,u)$ are constant. Moreover, to complement this system, we need boundary conditions. First, let us recall that the pressure wave entering the random section at $z=0$ is given by \eqref{coefP2}. Second, we assumed that no wave is coming from the right at $z=L/\e$ (see Figure \ref{figureP2} $(b)$). Consequently, we can complement the coupled mode system \eqref{mcP21}-\eqref{mcP24} with the boundary conditions
\begin{equation} \label{modeinit}\widehat{a}(\omega,0)=\widehat{a}^{\e,\xi}_{0}(\omega) \quad \text{ and } \quad\widehat{b}\left(\omega,\frac{L}{\e}\right)=\widehat{b}^{\e,\xi}_{L}(\omega) \quad\text{in }\esp_\xi,\end{equation}
where $\widehat{a}^{\e,\xi}_{0}(\omega)$ and $\widehat{b}^{\e,\xi}_{L}(\omega)$ are defined by \eqref{ab}, \eqref{defmode1}, and \eqref{defmode2}, for respectively $u=0$ and $u=L/\e$. In \eqref{modeinit}, $\widehat{a}^{\e,\xi}_{0}(\omega)$ represents the initial amplitudes of the right-going propagating and radiating modes at $u=0$, while $\widehat{b}^{\e,\xi}_{L}(\omega)$ represents the initial amplitudes of the left-going propagating and radiating modes at $u=L/\e$. The following proposition shows that $\widehat{a}^{\e,\xi}_{0}(\omega)$ is closed to the mode amplitudes of the pressure field \eqref{coefP2} coming from the right homogeneous part of the waveguide and entering the random section at $z=0$. As a result, $\widehat{a}^{\e,\xi}_{0}(\omega)$ is an approximation of $\widehat{a}^{\xi}_{0}(\omega)$  defined by \eqref{condinith} in the new system of coordinates introduced in Section \ref{conformal}. The convergence of $\widehat{b}^{\e,\xi}_{L}(\omega)$ toward $0$ in the following proposition is due to the fact that we assumed that no wave is coming from the right at $z=L/\e$ (see Figure \ref{figureP2} $(b)$). 
\begin{prop}\label{propsource} We have
\[\lim_{\e\to 0}\E\Big[\|\widehat{a}^{\e,\xi}_{0}(\omega)-\widehat{a}^\xi_{0}(\omega) \|^2_{\esp_\xi} +\|\widehat{b}^{\e,\xi}_{L}(\omega) \|^2_{\esp_\xi} \Big]=0. \]
\end{prop}
Proposition \ref{propsource} gives us an approximation of the coefficient $\widehat{a}^{\e,\xi}_{0}(\omega)$ in $\esp_\xi$, which is used in Theorem \ref{thasymptP21} and Theorem \ref{thasymptP12} to determine the initial conditions of the diffusion processes. This proposition is a direct consequence of Proposition \ref{asymptconf} and the Parseval's equality associated to the spectral decomposition introduced in Section \ref{spectralP2} since we have
\[\begin{split} 
 \int_0^{+\infty} \E\Big[ \lvert \widehat{p}_0(\omega, u_\ast,v)- \widehat{p}(\omega, u_\ast,v))  \rvert^2&+ \lvert \partial_u\widehat{p}_0(\omega, u_\ast,v)- \partial_u\widehat{p}(\omega, u_\ast,v))  \rvert^2 \Big]dv\\
&\leq C  \int_0^{+\infty} \E\Big[ \lvert z( u_\ast,v)- u_\ast  \rvert^2+ \lvert x( u_\ast,v)- v \rvert^2 \Big]dv,
 \end{split}\] 
where $u_\ast $ is either equal to $0$ or $L/\e$. Let us remark that this inequality holds because $u_\ast $ is out of the randomly perturbed section of the waveguide (see Figure \ref{figureP2}), so that we can use the variable $v\in[0,+\infty)$ of the unperturbed transverse section also for the original monochromatic pressure field $\widehat{p}(\omega, z,x)$.

\subsection{Influence of the Evanescent Modes on the Propagating and Radiating Modes}\label{iemprP2}

In this section we describe the influence of the evanescent modes on the coupling mechanism between the propagating and the radiating modes. The main goal is to obtain a more convenient coupled mode equation involving only the propagating and the radiating modes, but taking into account the contribution of the evanescent modes. However, to obtain the following proposition we introduce the following radiation condition for the evanescent modes meaning that the energy carried by the evanescent modes decay as the propagation distance becomes large.
\begin{assumption}\label{hyp6}
\[\lim_{u\to+\infty}\int_0^{+\infty}\Big \lvert  \int_{-1/\xi}^{-\xi} \widehat{p}_\ga (\omega,u)\phi_\ga(\omega,v)d\ga\Big\rvert^2 dv =0.\]
\end{assumption}
Consequently, thanks to Assumption \ref{hyp6} and following the proof of \cite[Section 4.3]{gomez2}, we can rewrite \eqref{mcP21}-\eqref{mcP24} in a closed  form in $\esp_\xi \times \esp_\xi$.
\begin{prop}\label{propresidual} We have
\begin{equation}\label{fullcoupledeq1}\begin{split}
\frac{d}{du} \widehat{a}(\omega, u)&= \textbf{\emph{H}}^{aa}_{\e,\xi}(\omega,u)\big(\widehat{a}(\omega,u)\big)+\textbf{\emph{H}}^{ab}_{\e,\xi}(\omega,u)\big(\widehat{b}(\omega, u)\big)+\textbf{\emph{R}}^{a,L_S}_{\e,\xi}(\omega,u)\\
&+\textbf{\emph{G}}^{aa}_{\e,\xi}(\omega,u)\big(\widehat{a}(\omega, u)\big)+\textbf{\emph{G}}^{ab}_{\e,\xi}(\omega,u)\big(\widehat{b}(\omega, u)\big)+\tilde{\textbf{\emph{R}}}^{a,L_S}_{\e,\xi}(\omega,u)+A^\e(u),
\end{split}\end{equation}
\begin{equation}\label{fullcoupledeq2}\begin{split}
\frac{d}{du} \widehat{b}(\omega, u)&=\overline{\textbf{\emph{H}}^{ab}_{\e,\xi}(\omega,u)\big(\overline{\widehat{a}(\omega, u)}\big)}+\overline{\textbf{\emph{H}}^{aa}_{\e,\xi}(\omega,z)\big(\overline{\widehat{b}(\omega, u)}\big)}+\textbf{\emph{R}}^{b,L_S}_{\e,\xi}(\omega,u)\\
&+\overline{\textbf{\emph{G}}^{ab}_{\e,\xi}(\omega,u)\big(\overline{\widehat{a}(\omega, u)}\big)}+\overline{\textbf{\emph{G}}^{aa}_{\e,\xi}(\omega,u)\big(\overline{\widehat{b}(\omega, u)}\big)}+\tilde{\textbf{\emph{R}}}^{b,L_S}_{\e,\xi}(\omega,u)+B^\e(u),
\end{split}\end{equation}
for all $ u\in[0,L/\e]$, where 
\begin{equation}\label{residual}
\sup_{u\in[0,L/\e]}\|A^\e(u)\|_{\esp_\xi}+\|B^\e(u)\|_{\esp_\xi}\leq C \frac{\e^{3/2-1/p}}{\xi}\sup_{u\in[0,L/\e]}\Big(\|\widehat{a}(\omega,u)\|_{\esp_\xi}+\|\widehat{b}(\omega,u)\|_{\esp_\xi}\Big) 
\end{equation}
for every $p>2$.
\end{prop}
In \eqref{fullcoupledeq1} and \eqref{fullcoupledeq2}, the operator $\textbf{H}^{aa}_{\e,\xi}$ and $\textbf{H}^{ab}_{\e,\xi}$ are defined by \eqref{haajP2}-\eqref{habgaP2}, and according to Assumption \ref{hyp3}
\begin{equation}\label{gaajP2}
\textbf{G}^{aa}_{\e,\xi,j}(u)(y)= \frac{ik^4}{4}\sum_{l=1}^{N}  \int_{-1/\xi}^{-\xi}\int   \frac{C^\e_{j\ga'}(u)C^\e_{\ga' l}(\tilde{u}+u)}{\sqrt{\beta_j\lvert \ga'\rvert \beta_l}}e^{-\sqrt{\lvert \ga'\rvert}\lvert \tilde{u}\rvert -i\beta_j u+i\beta_l (\tilde{u}+u)}y_l
\end{equation}
\[\textbf{G}^{ab}_{\e,\xi,j}(u)(y)= \frac{ik^4}{4}\sum_{l=1}^{N} \int_{-1/\xi}^{-\xi}\int   \frac{C^\e_{j\ga'}(u)C^\e_{\ga' l}(\omega,\tilde{u}+u)}{\sqrt{\beta_j(\omega)\lvert \ga'\rvert \beta_l }}e^{-\sqrt{\lvert \ga'\rvert}\lvert \tilde{u}\rvert -i\beta_j u -i\beta_l (\tilde{u}+u)}y_l
\]
\begin{equation}\label{gaagaP2}
\textbf{G}^{aa}_{\e,\xi,\ga}(u)(y)=\textbf{G}^{ab}_{\e,\xi,\ga}(u)(y)=0,\quad\text{for almost all } \ga\in (\xi,k^2),
\end{equation}
 with $y\in \esp_\xi$ and where $C^\e(\omega,u)=0$ if $u\not\in[0,L/\e]$. Here, the operators $\textbf{G}^{aa}_{\e,\xi}$ and $\textbf{G}^{ab}_{\e,\xi}$ represent the coupling between the evanescent modes with the propagating and the radiating modes. Moreover, $\textbf{R}^{a,L_S}_{\e,\xi}$, $\tilde{\textbf{R}}^{a,L_S}_{\e,\xi}$, $\textbf{R}^{b,L_S}_{\e,\xi}$, and $\tilde{\textbf{R}}^{b,L_S}_{\e,\xi}$ represent the excitation of the evanescent modes produced by the source term on the propagating and the radiating modes through the random perturbations. These terms are defined by 
\begin{equation}\label{Rj}
\textbf{R}^{a,L_S}_{\e,\xi,j} (u)=\overline{\textbf{R}^{b,L_S}_{\e,j} (u)}=\frac{i k^2}{2}\int_{-1/\xi}^{-\xi}\frac{C^\e_{j\ga'}(u)}{\sqrt{\beta_j}}\phi_{\ga'}( x_0)e^{-\sqrt{\lvert \ga' \rvert}(u-L_S)}d\ga' e^{-i \beta_j u},
\end{equation}
\begin{equation}\label{Rtj}\begin{split}
\tilde{\textbf{R}}^{a,L_S}_{\e,\xi,j} (u)=\overline{\tilde{\textbf{R}}^{b,L_S}_{\e,j} (u)}&=\frac{ik^{4}}{4}\int_{-1/\xi}^{-\xi}\int_{-1/\xi}^{-\xi} \int  \frac{C^\e_{j\ga'}(u)C^\e_{\ga' \ga''}(\tilde{u}+u)}{\sqrt{\beta_j \lvert\ga' \rvert}}e^{-\sqrt{\lvert \ga' \rvert}\lvert \tilde{u}\rvert-i\beta_j u}d\tilde{u} \\ 
&\times \phi_{\ga''}(x_0)e^{-\sqrt{\lvert \ga'' \rvert}(\tilde{u}+u-L_S)}d\ga''d\ga'\,,
\end{split}\end{equation}
\begin{equation}\label{Rga}
\textbf{R}^{a,L_S}_{\e,\xi,\ga} (u)=\textbf{R}^{b,L_S}_{\e,\ga} (u)=\textbf{R}^{a,L_S}_{\e,\xi,\ga} (u)=\tilde{\textbf{R}}^{b,L_S}_{\e,\ga} (u)=0,\quad\text{for almost all } \ga\in (\xi,k^2).
\end{equation}
Let us make two remarks. The first one concerns the fact that we take $p>2$ in order to have $3/2-1/p>1$ in \eqref{residual}, so that the correctors $A^\e$ and $B^\e$ are small enough to do not play any significant role in the forthcoming analysis. Second, the proof of this proposition given in Section \ref{proofpropresidual} follows closely the one of \cite[Section 4.3]{gomez2}. The only difference comes from error estimates near the boundary $v=0$ and the use of several truncated expansions in \eqref{randomcoef}.

Equations \eqref{fullcoupledeq1} and \eqref{fullcoupledeq2} describe the coupling process between the propagating and the radiating mode  including the influence of the evanescent modes and the evanescent part of the source term.  However, it is difficult to get good a priori estimates on 
\[\sup_{u\in[0,L/\e]}\Big(\|\widehat{a}(\omega,u)\|_{\esp_\xi}+\|\widehat{b}(\omega,u)\|_{\esp_\xi}\Big)\]
to be sure that the residuals $A^\e$ and $B^\e$ in Proposition \ref{propresidual} are small even after large propagation distance $L/\e$. One way consists in introducing the following stopping "time"
\[L^\e =\inf\Big(L>0, \quad \sup_{u\in[0,L/\e]}\|\widehat{a}(\omega,u) \|_{\esp_\xi}+\|\widehat{b}(\omega,u) \|_{\esp_\xi} \geq \frac{1}{\e^\alpha}\Big),\]
where $\alpha$ is such that $3/2-1/p-\alpha>1$, in order to limit the size of the random section and control the size of the residuals $A^\e$ and $B^\e$. However, will see in Section \ref{fsaP2} that under the forward scattering approximation, we have
\[\lim_{\e \to 0}\mathbb{P}(L^\e \leq L)=0\qquad \forall L >0,\]
which means that with high probability the residual terms are small.

\subsection{Forward Scattering Approximation}\label{fsaP2}

The forward scattering approximation is widely used in the literature, and consists in assuming that the coupling between forward- and backward-propagating modes is negligible compared to the coupling between the forward-propagating modes. In this case, the amplitudes of the left-going modes $\widehat{b}$ should converge to $0$ as $\e \to 0$ since we have assumed that no wave is coming from the right (see Figure \ref{figureP2} $(b)$). A rigorous derivation of the forward scattering approximation is technically complex. The main technical problem is that the mode amplitudes $(\widehat{a},\widehat{b})$ are not uniformly bounded on $[0,L/\e]$. To correct this issue we could think that it suffices to use in a first step the stopping times corresponding to the first exit times of closed balls with respect to the norm $\|.\|_{\esp_{\xi}}$ as described in \cite[Chapter 11]{stroock}. Unfortunately there are two problems. The first problem is that it is not possible to show a limit theorem on $\mathcal{C}([0,L],(\esp_{\xi},\|.\|_{\esp_{\xi}}))$. In fact, if the convergence holds on $\mathcal{C}([0,L],(\esp_{\xi},\|.\|_{\esp_{\xi}}))$ the energy conservation property \eqref{conservrelation} should be also valid for the asymptotic diffusion process according to the Portmanteau's theorem \cite{biling}. However, this conservation property contradict the energy dissipation of the limit process\footnote{We refer to \cite[Theorem 6.1]{gomez2} for the statement of this result in the case of perturbations in the interior of the waveguide. The same theorem holds in our context but the precise study of the asymptotic mean mode powers will be addressed in a later work}. As a result, it seems to be more appropriate to obtain a limit theorem on $\mathcal{C}([0,L],\esp_{\xi,w})$ where $\esp_{\xi,w}$ stands for $\esp_\xi$ equipped with the weak topology. However, the second problem is that these first exit times are lower semicontinuous with respect to the topology of $\mathcal{C}([0,L],(\esp_{\xi},\|.\|_{\esp_{\xi}}))$, but not with respect to the one of $\mathcal{C}([0,L],\esp_{\xi,w})$. Therefore, the classical technique of \cite[Chapter 11]{stroock} cannot be applied.

However, one can formally show using the diffusion approximation theorem proved in \cite{book} that the coupling between right-going propagating modes and left-going propagating modes involves coefficients of the form
\[ \int_{0}^{+\infty}R_s(u)\cos\big((\Bh{l}{}+\Bh{j}{})u\big)du\quad\text{and}\quad \int_{0}^{+\infty}R_b(u)\cos\big((\Bh{l}{}+\Bh{j}{})u\big)du,\] 
while the coupling between two right-going propagating modes or two left-going propagating modes involves  coefficients of the form
\[ \int_{0}^{+\infty}R_s(u)\cos\big((\Bh{l}{}-\Bh{j}{})u\big)du\quad\text{and}\quad \int_{0}^{+\infty}R_b(u)\cos\big((\Bh{l}{}-\Bh{j}{})u\big)du,\] 
for all $(j,l) \in \big\{1,\dots,\N{}\big\}^2$, and where $R_s,R_b\in L^1(\mathbb{R})$ are the autocorrelations function of $V_s$ and $V_b$ (see \eqref{formuleI}). Consequently, in our context if we forget these technical problems, according to \cite{book, papa} the forward scattering approximation should be valid in the asymptotic $\e\to 0$ under the following assumption, meaning that the Fourier transform of its $u$-autocorrelation function possesses a cut-off wavenumber. 
\begin{assumption}\label{hyp4} 
For all $(j,l) \in \big\{1,\dots,\N{}\big\}^2$,
\[ \int_{0}^{+\infty}R_s(u)\cos\big((\Bh{l}{}+\Bh{j}{})u\big)du=0\quad\text{and}\quad  \int_{0}^{+\infty}R_b(u)\cos\big((\Bh{l}{}+\Bh{j}{})u\big)du=0.\]
\end{assumption}
As a result, under these assumptions, there is no effective mode coupling between the right-going and left-going propagating modes, but there is still coupling between the right-going propagating modes which is described in Section \ref{coupledprocP2}. Let us remark that the continuous part $(0,\ko{})$ of the spectrum corresponding to the radiating modes play no role in the previous assumption. The reason is that the radiating modes and the evanescent modes play no role in the coupling mechanism between the right- and left-going modes as it has been shown in \cite{gomez2,gomez3}.

Consequently, according to the forward scattering approximation we only consider the simplified version of the forward coupled mode equation \eqref{fullcoupledeq1} and \eqref{fullcoupledeq2}:
\begin{equation}\label{eqdiffa}
\dz \widehat{a}(\omega, u)=\textbf{H}^{aa}_{\e,\xi}(\omega,u)\left(\widehat{a}(\omega, u)\right)+ \textbf{G}^{aa}_{\e,\xi}(\omega,u)\left(\widehat{a}(\omega, u)\right)+\textbf{R}^{a,L_S}_{\e,\xi}(\omega,u)+\tilde{\textbf{R}}^{a,L_S}_{\e,\xi}(\omega,u)+A^\e(u)
\end{equation}
in $\esp_\xi$, for $u\in[0,L/\e]$, where $\textbf{H}^{aa}_{\e,\xi}$, $\textbf{G}^{aa}_{\e,\xi}$, $\textbf{R}^{a,L_S}_{\e,\xi}$ and $\tilde{\textbf{R}}^{a,L_S}_{\e,\xi}$ are defined by \eqref{haajP2}-\eqref{haagaP2}, \eqref{gaajP2}-\eqref{gaagaP2}, and \eqref{Rj}-\eqref{Rga}, and with
\[
\sup_{u\in[0,L/\e]}\|A^\e(u)\|_{\esp_\xi}\leq C \frac{\e^{3/2-1/p}}{\xi}\sup_{u\in[0,L/\e]}\|\widehat{a}(\omega,u)\|_{\esp_\xi},\qquad \forall p>2.
\]
The evolution equation \eqref{eqdiffa} is complemented with initial condition $\widehat{a}(\omega, 0)=\widehat{a}^{\e,\xi}_0(\omega)$, where $\widehat{a}^{\e,\xi}_0(\omega)$ is introduced in Section \ref{couplemodeeq}. As already discussed, under the forward scattering approximation, we have the following proposition.
\begin{prop}\label{TA} We have
\[\lim_{\e \to 0}\mathbb{P}(L^\e \leq L)=0, \quad\text{where}\quad L^\e =\inf\Big(L>0, \quad \sup_{u\in[0,L/\e]}\|\widehat{a}(\omega,u) \|_{\esp_\xi} \geq \frac{1}{\e^\alpha}\Big),\] 
\end{prop}
The proof of Proposition \ref{TA} follows exactly the one of \cite[Proposition 4.1]{gomez2bis}. Consequently, under the forward scattering approximation, one can adapt the results obtained in \cite[Section 4.2]{gomez2} to the system \eqref{eqdiffa} in order to derive the local energy flux conservation for the propagating and the radiating modes, that is for all $\eta>0$
\begin{equation}\label{conservrelation} \lim_{\e\to0}\mathbb{P}\Big(\sup_{u\in[0,L/\e]}\left\lvert   \|\widehat{a}(\omega,u)\|^2_{\esp_\xi}- \|\widehat{a}^\xi_0(\omega)\|^2_{\esp_\xi} \right\rvert>\eta  \Big)=0,\end{equation}
meaning that the amplitude $\widehat{a}(\omega, u)$ is asymptotically uniformly bounded on $[0,L/\e]$ in the limit $\e\to0$. The derivation of \eqref{conservrelation} is based on Lemma \ref{approxcoef}, Proposition \ref{propsource}, and Proposition \ref{TA}. Consequently, using the Gronwall's lemma, Proposition \ref{propsource}, and Proposition \ref{TA}, we have the following proposition.
\begin{prop}\label{propcme}
We have for all $\eta>0$
\[ \lim_{\e\to0}\mathbb{P}\Big(\sup_{u\in[0,L/\e]}  \|\widehat{a}(\omega,u)-\widehat{a}_1(\omega,u)\|^2_{\esp_\xi}>\eta  \Big)=0,\]
where $\widehat{a}_1(\omega,u)$ is the unique solution of the differential equation
\begin{equation}\label{approxcme}
\frac{d}{du} \widehat{a}_1(\omega, u)= \textbf{\emph{H}}^{aa}_{\e,\xi} (\omega,u)(\widehat{a}_1(\omega, u))+ \textbf{\emph{G}}^{aa}_{\e,\xi}(\omega,u)(\widehat{a}_1(\omega, u))\quad\text{with}\quad   \widehat{a}_1(\omega,0)=\widehat{a}^\xi_0(\omega).
\end{equation}
\end{prop} 
As a result, thanks to Proposition \ref{propcme} and \cite[Theorem 3.1]{biling}, we can focus our attention in what follows on the process $\widehat{a}_1(\omega,u)$, which is studied in detail in Section \ref{coupledprocP2}. Moreover, this proposition implies that the information transmitted to the propagating and radiating modes about the evanescent part of the source ($\textbf{R}^{a,L_S}_{\e,\xi}$ and $\tilde{\textbf{R}}^{a,L_S}_{\e,\xi}$) is lost during the propagation over large distances.

\section{Coupled Mode Processes}\label{coupledprocP2}

In this section, we study in a first time the asymptotic behavior, as $\e\to 0$ first and then $\xi \to 0$, of the statistical properties of the mode coupling mechanism \eqref{approxcme} in terms of a diffusion process. In a second time we study the attenuation properties of the propagating-mode amplitudes.

First of all, let us define the rescaled process according to the size of the random section $[0,L/\e]$,
\[\widehat{a} ^{\e,\xi} (\omega ,u)=\widehat{a}_1\left(\omega,\frac{u}{\e}\right)\quad \forall u\in[0,L].\]
 This process is the unique solution on $[0,L]$ of the rescaled forward coupled mode equations 
\begin{equation}\label{eqfora}
\frac{d}{du} \widehat{a}^{\e,\xi}(\omega, u)=\frac{1}{\e} \textbf{H}^{aa}_{\e,\xi} \Big(\omega,\frac{u}{\e}\Big)(\widehat{a}^{\e,\xi}(\omega, u))+\frac{1}{\e} \textbf{G}^{aa}_{\e,\xi} \Big(\omega,\frac{u}{\e}\Big)(\widehat{a}^{\e,\xi}(\omega, u))
\end{equation}
in $\esp_\xi$, with initial conditions $\widehat{a}^{\e,\xi}(\omega, 0)=\widehat{a}^\xi_0(\omega)$ defined by \eqref{condinith}, and where $\textbf{H}^{aa}_{\e,\xi}$ and $\textbf{G}^{aa}_{\e,\xi}$ are defined by \eqref{haajP2}-\eqref{haagaP2} and \eqref{gaajP2}-\eqref{gaagaP2}.

\subsection{Limit Theorem}\label{thasyptoticP2}

In this section, to simplify the presentation we first focus our attention only on the random perturbations at the surface of the waveguide, that is we assume that $V^\e_b=0$. Afterward, we take into account also the randomly perturbed bottom topography of the waveguide. To prove Theorem \ref{thasymptP21} and Theorem \ref{thasymptP12} we consider the classical nondegeneracy condition \cite{book,papa, papanicolaou}.
\begin{assumption}\label{hyp5}
The modal wavenumbers $(\beta_1(\omega),\dots,\beta_{\N{}}(\omega))$ defined in Section \ref{spectralP2} are distinct along with their sums and differences. 
\end{assumption}  
Let us note that this condition is not necessarily satisfied for nonplanar waveguide model. However, in this case, Assumption \ref{hyp5} can be still valid if the  statistical properties of the random fluctuations are rotationally invariant \cite{papanicolaou}.

\subsubsection{The Random Surface}\label{randomsurfsec}

This section is devoted to the asymptotic analysis of the forward mode amplitude $\widehat{a}^{\e,\xi}(\omega, u)$ solution of \eqref{eqfora}, and where the random perturbations of the waveguide are only due to the random surface, that is $V^\e_b=0$.
\begin{thm}\label{thasymptP21}
In the case $V^\e_b=0$, under Assumptions \ref{hyp1}--\ref{hyp5} and the mixing conditions \eqref{mixcond}, for all $L>0$ the family $\big(\widehat{a}^{\e,\xi}(\omega, \cdot)\big)_{\e\in(0,1)}$, unique solution of \eqref{eqfora}, converges in distribution on 
$\mathcal{C}([0,L], \esp_{\xi,w})$ as $\e \to 0$ to a limit denoted by $\widehat{a}^\xi(\omega, \cdot)$. Here $\esp_{\xi,w}$ stands for the Hilbert space $\esp_{\xi}$ equipped with the weak topology. This limit is the unique diffusion process on $\esp_{\xi}$, starting from $\widehat{a}^\xi_0(\omega)$ defined by \eqref{condinith}, and associated to the infinitesimal generator
\begin{equation}\label{Ls}\mathcal{L}^{\omega,s}_\xi(u)=\mathcal{L}^{\omega,s}_1(u)+\mathcal{L}^{\omega,s}_{2,\xi}(u)+\mathcal{L}^{\omega,s}_{3,\xi}(u),\end{equation}
where
\[\begin{split}
\mathcal{L}^{\omega,s}_{1}(u)&=\frac{f_s^2(u)}{2}\sum_{\substack{j,l =1\\ j\not=l}}^{\N{}}\Gamma^{c,s}_{jl}(\omega)\left(T_j\overline{T_j}\partial_{T_l}\partial_{\overline{T_l}}+T_l\overline{T_l}\partial_{T_j}\partial_{\overline{T_j}}-T_jT_l\partial_{T_j}\partial_{T_l}-\overline{T_j}\overline{T_l}\partial_{\overline{T_j}}\partial_{\overline{T_l}}\right)\\
&+\frac{f_s^2(u)}{2}\sum_{j=1}^{\N{}} \left( \Gamma^{c,s}_{jj}(\omega) -\Gamma^{1,s}_{jj}(\omega)\right)\left(T_j \partial_{T_j}+\overline{T_j}\partial_{\overline{T_j}}\right)+i \frac{ f_s^2(u)}{2}\sum_{j=1}^{\N{}} \Gamma^{s,s}_{jj}(\omega) \left(T_j \partial_{T_j}-\overline{T_j}\partial_{\overline{T_j}}\right),
\end{split}\]
and
\begin{equation}\label{driftls}\begin{split}
\mathcal{L}^{\omega,s}_{2,\xi}&=-\frac{f_s^2(u)}{2}\sum_{j=1}^{\N{}} \left(\Lambda^{c,s,\xi}_j(\omega) +i\Lambda^{s,s,\xi}_j(\omega)\right)    T_j \partial_{T_j}+ \left(\Lambda^{c,s,\xi}_j(\omega)  - i\Lambda^{s,s,\xi}_j(\omega)\right)\overline{T_j}\partial_{\overline{T_j}},\\
\mathcal{L}^{\omega,s}_{3,\xi}&=i f_s^2(u) \sum_{j=1}^{\N{}} \kappa^{s,\xi}_j(\omega) \left(T_j \partial_{T_j}-\overline{T_j}\partial_{\overline{T_j}}\right).
\end{split}\end{equation}
\end{thm}
Let us note that according to the infinitesimal generator $\mathcal{L}^{\omega,s}_\xi$ defined by \eqref{Ls}, the diffusion phenomena happen on the support of $f_s$ describing the location of the random perturbations of the surface of the waveguide \eqref{randompert}. 

Here, we have considered the classical complex derivative with the following notation: If $v=v_1 + i v_2$, then $\partial_v =\frac{1}{2}\left(\partial_{v_1} -i \partial_{v_2} \right)$ and $\partial_{\overline{v}}=\frac{1}{2}\left(\partial_{v_1} +i \partial_{v_2} \right)$. The coupling coefficients are given by the following notations. Let us denote
\begin{equation}\label{defQ}\begin{split}
Q(r_1,s_1&,r_2,s_2,A)=\frac{4}{\pi}(n_1^2-1)^2\phi_{r_1}(d)\phi_{r_2}(d)\phi_{s_1}(d)\phi_{s_2}(d)e^{-2d\vert A\vert}\\
&+\frac{8}{\pi}(n_1^2-1)\vert A\vert\phi_{r_1}(d)\phi_{s_1}(d)\int_0^{+\infty} \phi_{r_2}(v)\phi_{s_2}(v)n^2(v)e^{-(d+v)\vert A\vert}dv\\
&+\frac{4}{\pi}A^2 \int_0^{+\infty}\int_0^{+\infty} \phi_{r_1}(v)\phi_{s_1}(v) \phi_{r_2}(v')\phi_{s_2}(v')n^2(v)n^2(v')e^{-(v+v')\vert A\vert}dvdv'.
\end{split}\end{equation}
We have for all $(j,l)\in\big\{1,\dots,N\big\}^{2}$ and $j\not=l$
\begin{equation}\label{coefcouplagelimit}\begin{split}
\Gamma^{c,s}_{jl}&= \frac{k^4}{2\beta_j\beta_l}Q(j,l,j,l,\beta_j-\beta_l)\int_{0}^{+\infty}R_{s}(z)\cos\big((\beta_l-\beta_j)z\big)dz,\\
\Gamma^{s,s}_{jl}&= \frac{k^4}{2\beta_j \beta_l}Q(j,l,j,l,\beta_j-\beta_l)\int_{0}^{+\infty}R_{s}(z)\sin\big((\beta_l-\beta_j)z\big)dz, \\
\Gamma^{c,s}_{jj}&=-\sum_{\substack{l=1\\l\not=j}}^{N}\Gamma^{c}_{jl},\qquad \Gamma^{s,s}_{jj}=-\sum_{\substack{l=1\\ l\not= j}}^{N}\Gamma^{s}_{jl},
\end{split}\end{equation}
and for all
$(j,l)\in\big\{1,\dots,N \big\}^{2}$,
\begin{equation}\label{coefatt}\begin{split}
\Lambda^{c,s,\xi}_{j}&= \int_{\xi}^{k^2}\frac{k^4 }{2\sgap\beta_j}Q(\ga',j,\ga',j,\sgap-\beta_j)\int_{0}^{+\infty}R_{s}(z)\cos\big((\sgap-\beta_j)z\big)dzd\ga', \\
\Lambda^{s,s,\xi}_{j}&= \int_{\xi}^{k^2} \frac{k^4 }{2\sgap \beta_j}Q(\ga',j,\ga',j,\sgap-\beta_j)\int_{0}^{+\infty}R_{s}(z)\sin\big((\sgap-\beta_j)z\big)dzd\ga',\\
\kappa^{s,\xi}_j &=\int_{-\infty}^{-\xi}\frac{i k^{4}}{2\beta_j \sqrt{\lvert\ga\rvert}} \int_{0}^{+\infty}G^{(1)}_{j\ga}(\omega,z)\cos\big(\beta_j z\big)e^{-\sqrt{\lvert\ga\rvert}z}dzd\ga
+\frac{ik^2}{2\beta_j}G^{(2)}_j .
\end{split}\end{equation}
Moreover, we have
 \[
G^{(1)}_{j\ga}(z)=\int_0^{+\infty}I_s(u) Q(\ga,j,\ga,j,u)\cos(uz)du\]
 and	
 \[\begin{split}
 G^{(2)}_j&=(n_1^2-1)\phi_j^2(d)\Big[6\int_0^{+\infty} I_s(u) ue^{-2du}du+\mathcal{U}(\tilde{R}_s)(0)\Big]\\
 &+\frac{4\pi^2}{d^2}(n_1^2-1)\phi_j(d)\phi'_j(d)\int_0^{+\infty} I_s(u) e^{-2du}du+\int_0^{+\infty}  \phi_j^2(v)n^2(v)\int_0^{+\infty} I_s(u) e^{-2uv}(1+u^2)du.
 \end{split}\]
 In the previous expressions $n(v)$ is defined by \eqref{indexrefraction2}, $R_s$ stands for the autocorrelation function of the random process $V_s$, and $\tilde{R}_s$ is the correlation function of the processes $V_s$ and $\frac{d}{dz}V_s$, which is $\frac{d}{dz}R_s$. 
Let us remark that the first (resp. the second and the third) term in the right hand side of \eqref{defQ}, and $G^{(1)}$ are produced by the coupling mechanism between the propagating modes at the bottom of the waveguide (resp. between the bottom and the transverse section, and  through the transverse section of the waveguide).

In the proof of Theorem \ref{thasymptP21}, the random operator $\textbf{G}^{aa}_{\e,\xi}$ in \eqref{eqfora} can be treated following the ideas of \cite{khasminskii} as it has been done in \cite{gomez2,gomez3}. In fact, the technique developed in \cite{khasminskii} does not require any mixing properties. However, the random operator $\textbf{H}^{aa}_{\e,\xi}$ requires to introduce an approximation of the random coefficients \eqref{randomcoef} to turn this formula in a more convenient form to use the $\phi$-mixing properties of the filtration \eqref{filtration}. The approximations of the random operators $\textbf{H}^{aa}_{\e,\xi}$ and  $\textbf{G}^{aa}_{\e,\xi}$ are given in Lemma \ref{approxcoef}, Section \ref{approxrandomcoef}. The remaining of the proofs are based on a martingale approach using the perturbed-test-function method. In a first step we show the tightness of the processes, and in a second step we characterize all the subsequence limits by mean of a well-posed martingale problem in a Hilbert space. Finally, let us remark that the convergence in Theorem \ref{thasymptP21} holds also on $\mathcal{C}([0,L], (\esp_{\xi},\|.\|_{\esp_{\xi}}))$ only for the $\N{}$-discrete propagating mode amplitudes.

\subsubsection{The Random Surface and Random Bottom}

In this section we present the asymptotic analysis of the forward mode amplitudes $\widehat{a}^{\e,\xi}(\omega, u)$, where we consider a randomly perturbed bottom topography in addition to the randomly perturbed free surface. This section presents no additional difficulties since the random perturbations of the surface and the bottom are assumed to be independent, and can be treated directly using the method of \cite{gomez2,gomez3}.

\begin{thm}\label{thasymptP12}
Under Assumptions \ref{hyp1}--\ref{hyp5} and the mixing conditions \eqref{mixcond}, for all $ L>0$ the family $\big(\widehat{a}^{\e,\xi}(\omega, \cdot)\big)_{\e\in(0,1)}$, unique solution of \eqref{eqfora}, converges in distribution on 
$\mathcal{C}([0,L], \esp_{\xi,w})$ as $\e \to 0$ to a limit denoted by $\widehat{a}^\xi(\omega, \cdot)$. Here $\esp_{\xi,w}$ stands for the Hilbert space $\esp_{\xi}$ equipped with the weak topology. This limit is the unique diffusion process on $\esp_{\xi}$, starting from $\widehat{a}^\xi_0(\omega)$ defined by \eqref{condinith}, and associated to the infinitesimal generator
\begin{equation}\label{geneapprox}\mathcal{L}^{\omega}_\xi(u)=\mathcal{L}^{\omega,s}_\xi(u)+\mathcal{L}^{\omega,b}_\xi(u),\end{equation}
where $\mathcal{L}^{\omega,s}_\xi(u)$ is defined by \eqref{Ls}, and
\begin{equation}\label{Lb}\mathcal{L}^{\omega,b}_\xi(u)=\mathcal{L}^{\omega,b}_1(u)+\mathcal{L}^{\omega,b}_{2,\xi}(u)+\mathcal{L}^{\omega,b}_{3,\xi}(u),\end{equation}
with
\begin{equation}\label{geneb}\begin{split}
\mathcal{L}^{\omega,b}_{1}(u)&=\frac{f_b^2(u)}{2}\sum_{\substack{j,l =1\\ j\not=l}}^{\N{}}\Gamma^{c,b}_{jl}(\omega)\left(T_j\overline{T_j}\partial_{T_l}\partial_{\overline{T_l}}+T_l\overline{T_l}\partial_{T_j}\partial_{\overline{T_j}}-T_jT_l\partial_{T_j}\partial_{T_l}-\overline{T_j}\overline{T_l}\partial_{\overline{T_j}}\partial_{\overline{T_l}}\right)\\
&+\frac{f_b^2(u)}{2}\sum_{j,l =1}^{\N{}}\Gamma^{1,b}_{jl}(\omega)\left(T_j\overline{T_l}\partial_{T_j}\partial_{\overline{T_l}}+\overline{T_j}T_l\partial_{\overline{T_j}}\partial_{T_l}-T_jT_l\partial_{T_j}\partial_{T_l}-\overline{T_j}\overline{T_l}\partial_{\overline{T_j}}\partial_{\overline{T_l}}\right)\\
&+\frac{f_b^2(u)}{2}\sum_{j=1}^{\N{}} \big(\Gamma^{c,b}_{jj}(\omega)-\Gamma^{1,b}_{jj}(\omega)\big)\left(T_j \partial_{T_j}+\overline{T_j}\partial_{\overline{T_j}}\right)+i \frac{f_b^2(u)}{2}\sum_{j=1}^{\N{}} \Gamma^{s,b}_{jl}(\omega)\left(T_j \partial_{T_j}-\overline{T_j}\partial_{\overline{T_j}}\right),
\end{split}\end{equation}
and
\begin{equation}\label{driftlb}\begin{split}
\mathcal{L}^{\omega,b}_{2,\xi}&=-\frac{f_b^2(u)}{2}\sum_{j=1}^{\N{}} (\Lambda^{c,b,\xi}_j(\omega) +i\Lambda^{s,b,\xi}_j(\omega))    T_j \partial_{T_j}+ (\Lambda^{c,b,\xi}_j(\omega) -i\Lambda^{s,b,\xi}_j(\omega))\overline{T_j}\partial_{\overline{T_j}},\\
\mathcal{L}^{\omega,b}_{3,\xi}&=i f_b^2(u) \sum_{j=1}^{\N{}} \kappa^{b,\xi}_j(\omega) \left(T_j \partial_{T_j}-\overline{T_j}\partial_{\overline{T_j}}\right).
\end{split}\end{equation}
\end{thm}

The coupling coefficients in \eqref{geneb} and \eqref{driftlb} are defined as follows: for all $(j,l)\in\big\{1,\dots,\N{} \big\}^{2}$ and $j\not=l$
\begin{equation}\label{coefcouplagelimit2}\begin{split}
\Gamma^{c,b}_{jl}&= \frac{2k^4(n_1^2-1)^2\phi_j^2(d)\phi_l^2(d)}{\beta_j\beta_l}\int_{0}^{+\infty}R_b(z) \cos\big((\beta_l-\beta_j)z\big)dz,\\
\Gamma^{s,b}_{jl}&= \frac{2k^4(n_1^2-1)^2\phi_j^2(d)\phi_l^2(d)}{\beta_j\beta_l}\int_{0}^{+\infty}R_b(z) \sin\big((\beta_l-\beta_j)z\big)dz, \\
\Gamma^{c,b}_{jj}&=-\sum_{\substack{l=1\\l\not=j}}^{N}\Gamma^{c,b}_{jl},\qquad\Gamma^{s,b}_{jj}=-\sum_{\substack{l=1\\ l\not= j}}^{N}\Gamma^{s,b}_{jl},
\end{split}\end{equation}
and for all $(j,l)\in\big\{1,\dots,N \big\}^{2}$,
\begin{equation}\label{coefatt2}\begin{split}
\Gamma^{1,b}_{jl}&=  \frac{2k^4(n_1^2-1)^2\phi_j^2(d)\phi_l^2(d)}{\beta_j\beta_l}\int_{0}^{+\infty}R_b(z) dz,\\
\Lambda^{c,b,\xi}_{j}&= \int_{\xi}^{k^2}\frac{2k^4(n_1^2-1)^2\phi_j^2(d)\phi_{\ga'}^2(d)}{\sgap \eta_j}\int_{0}^{+\infty}R_b(z)\cos\big((\sgap-\beta_j)z\big)dzd\ga', \\
\Lambda^{s,b,\xi}_{j}&= \int_{\xi}^{k^2} \frac{2k^4(n_1^2-1)^2\phi_j^2(d)\phi_{\ga'}^2(d)}{\sgap \beta_j}\int_{0}^{+\infty}R_b(z)\sin\big((\sgap-\beta_j)z\big)dzd\ga',\\
\kappa^{b,\xi}_j &=i\int_{-1/\xi}^{-\xi}\frac{2k^4(n_1^2-1)^2\phi_j^2(d)\phi_{\ga'}^2(d)}{\beta_j \sqrt{\lvert\ga\rvert}} \int_{0}^{+\infty}R_b(z)\cos\big(\beta_j z\big)e^{-\sqrt{\lvert\ga\rvert}z}dzd\ga  \\
& \hspace{0.5cm}+iR_b(0)\frac{8k^2(n_1^2-1)^2\phi^2_j(d){\phi'}^2_j(d)}{\beta_j}.
\end{split}\end{equation}

The drift ($\mathcal{L}^{\omega,s}_{2,\xi}$, $\mathcal{L}^{\omega,s}_{3,\xi}$, $\mathcal{L}^{\omega,b}_{2,\xi}$, and $\mathcal{L}^{\omega,b}_{3,\xi}$ defined by \eqref{driftls} and \eqref{driftlb}) and the initial value $\widehat{a}^\xi_0(\omega)$ defined by  \eqref{condinith} of the diffusion process $\widehat{a}^\xi(\omega, \cdot)$ still depend on the parameter $\xi$ introduced in Assumption \ref{hyp2}. Before discussing the meaning of each terms in \eqref{Ls} and \eqref{Lb}, we give the following result regarding the asymptotic $\xi\to 0$. 

\begin{thm}\label{thasympP2}
For all $L>0$, the family $\big(\widehat{a}^\xi(\omega, \cdot)\big)_{\xi\in(0,1)}$ converges in distribution on 
$\mathcal{C}([0,L], (\esp_{0},\|.\|_{\esp_{0}}))$ as $\xi \to 0$ to a limit denoted by $\widehat{a}(\omega, \cdot)$. This limit is the unique diffusion process on $\esp_{0}$, starting from 
\begin{equation}\label{sourcexi}\widehat{a}_0(\omega)=\lim_{\xi\to 0}\widehat{a}^\xi_0(\omega),\end{equation}
where $\widehat{a}^\xi_0(\omega)$ is defined by \eqref{condinith}, and associated to the infinitesimal generator
\begin{equation}\label{geneapprox2}\mathcal{L}^\omega(u)=\mathcal{L}^\omega_1(u)+\mathcal{L}^\omega_2(u)+\mathcal{L}^\omega_3(u),\end{equation}
where
\[\begin{split}
\mathcal{L}^{\omega}_{2}&=-\frac{1}{2}\sum_{j=1}^{\N{}}f_s^2(u) \Big((\Lambda^{c,s}_j(\omega) +i\Lambda^{s,s}_j(\omega))  T_j \partial_{T_j}+(\Lambda^{c,s}_j(\omega) -i\Lambda^{s,s}_j(\omega))\overline{T_j}\partial_{\overline{T_j}}\Big)\\
&\hspace{1.5cm}+ f_b^2(u) \Big((\Lambda^{c,b}_j(\omega) +i\Lambda^{s,b}_j(\omega))  T_j \partial_{T_j}+(\Lambda^{c,b}_j(\omega) -i\Lambda^{s,b}_j(\omega))\overline{T_j}\partial_{\overline{T_j}}\Big),\\
\mathcal{L}^{\omega}_{3}&=i  \sum_{j=1}^{\N{}}\left (f_s^2(u) \kappa^{s}_j(\omega)+f_b^2(u) \kappa^{b}_j(\omega)\right) \left(T_j \partial_{T_j}-\overline{T_j}\partial_{\overline{T_j}}\right),
\end{split}\]
with for all $j\in\big\{1,\dots,\N{} \big\}$
\begin{equation}\label{lim1}
\Lambda^{c,s}_{j}(\omega)=\lim_{\xi\to 0}\Lambda^{c,s,\xi}_{j}(\omega),\quad
\Lambda^{s,s}_{j}(\omega)=\lim_{\xi\to 0}\Lambda^{s,s,\xi}_{j}(\omega),\quad
\kappa^s _j (\omega)=\lim_{\xi \to 0}\kappa ^{s,\xi}_j (\omega),
\end{equation}
and
\begin{equation}\label{lim2}
\Lambda^{c,b}_{j}(\omega)=\lim_{\xi\to 0}\Lambda^{c,b,\xi}_{j}(\omega),\quad
\Lambda^{s,b}_{j}(\omega)=\lim_{\xi\to 0}\Lambda^{s,b,\xi}_{j}(\omega),\quad
\kappa^b _j (\omega)=\lim_{\xi \to 0}\kappa ^{b,\xi}_j (\omega).
\end{equation}
\end{thm}
Let us note that in Theorem \ref{thasympP2} the limits in \eqref{lim1} and \eqref{lim2} with respect to $\xi$ are well defined thanks to \eqref{phiga}-- \eqref{aga}.

Consequently, Theorem \ref{thasymptP21}, Theorem \ref{thasymptP12}, and Theorem \ref{thasympP2} describe the asymptotic behavior, as $\e\to 0$ first and then $\xi \to 0$, of the statistical properties of the forward mode amplitudes $\widehat{a}^{\e,\xi}(\omega, \cdot)$ in terms of a diffusion process with infinitesimal generator \eqref{geneapprox2} and starting from \eqref{sourcexi}. The infinitesimal generator $\mathcal{L}^\omega$ is composed of three parts which represent different behaviors on the diffusion process. However, we can remark that this infinitesimal generator depends only on the $\N{}$-discrete coordinates, so that the radiating part of the forward mode amplitudes remains constant. This result has already been obtained in \cite{gomez2,gomez3} in a different setup. The first operator $\mathcal{L}^\omega_1$ describes the mode coupling between the $\N{}$-propagating modes. This part is of the form of the infinitesimal generator obtained in \cite{book,papa}, and for which the total energy is conserved.  The second operator $\mathcal{L}^\omega_2$ describes the coupling between the propagating modes with the radiating modes. This part implies a mode-dependent and frequency-dependent attenuation on the $\N{}$-propagating modes, and a mode-dependent and frequency-dependent phase modulation. The third operator $\mathcal{L}^\omega_3$ describes the coupling between the propagating and the evanescent modes, and implies a mode-dependent and frequency-dependent phase modulation. The purely imaginary part of the operator $\mathcal{L}^\omega$ does not remove energy from the propagating modes but gives an effective dispersion. Finally, let us remark that their is no $\Gamma^{1,s}$ term for the surface fluctuations. The reason is that $I_s(0)=0$ for \eqref{formuleI}  (Pierson-Neuman and Pierson-Moscovitz spectra), that is a mode cannot be coupled with himself.

\subsection{Mean Mode Amplitudes}\label{meanmode}

In this section we describe the effects of the random perturbations on the forward mean mode amplitudes. From Theorem \ref{thasympP2}, we get the following result.
\begin{prop}\label{coherence}
For all $z\in[0,L]$ and $ j\in\big\{1,\dots,\N{} \big\}$, we have 
\[\begin{split}
\lim_{\xi\to0}\lim_{\e\to0}&\mathbb{E}\Big[\widehat{a}^{\xi,\e}_j( \omega,z)\Big]=
\mathbb{E}\Big[\widehat{a}_j( \omega,z)\Big]\\
&=\exp\left[\left(\frac{\Gamma^{c,s}_{jj}(\omega)-\Gamma^{1,s} _{jj}(\omega) - \Lambda^{c,s}_j (\omega)+ i(\Gamma^{s,s}_{jj}(\omega)-\Lambda^{s,s}_{jj}(\omega)+2k^s_j(\omega))}{2}\right)\int_0^z f^2_s(u)du \right]\\
&\times\exp\left[\left(\frac{\Gamma^{c,b}_{jj}(\omega)-\Gamma^{1,b} _{jj}(\omega) - \Lambda^{c,b}_j (\omega)+ i(\Gamma^{s,b}_{jj}(\omega)-\Lambda^{s,b}_{jj}(\omega)+2k^b_j(\omega))}{2}\right)\int_0^z f^2_b(u)du \right]\widehat{a}_{j,0} (\omega),
\end{split}\]
where $\widehat{a}_{0} (\omega)$ is defined by \eqref{sourcexi}.
\end{prop}
First, let us note that the mean amplitude of the radiating part remains constant on $L^2(0,\ko{})$, since the diffusion process only holds for the propagating modes (see Theorem \ref{thasympP2}). Second, for all $ j\in \big\{1,\dots,\N{}\big\}$, the coefficients $(\Gamma^1 _{jj}(\omega)+ \Lambda^{c}_j (\omega)-\Gamma^{c}_{jj}(\omega))/2$ are nonnegative according to the Bochner's theorem \cite{reed} or \cite[Section 6.3.6]{book}.

The decay rate for the mean $j$th-propagating mode is given by
\begin{equation}\label{decayrate}
\Big\lvert\mathbb{E}\Big[ \widehat{a}_j( \omega,z)\Big]\Big\rvert= \big\lvert \widehat{a}_{j,0} (\omega) \big\rvert e^{-( \Lambda^{c,s}_j (\omega)-\Gamma^{c,s} _{jj}(\omega)) \int_0^L f^2_s(u)du/2-(\Gamma^{1,b}_{jj}(\omega)+ \Lambda^{c,b}_j (\omega)-\Gamma^{c,b} _{jj}(\omega)) \int_0^L f^2_b(u)du/2}, 
\end{equation}
which depends on the effective coupling between the propagating modes, and the coupling between the propagating and the radiating modes. This decay describes the effective attenuation of the mode amplitudes caused by the cumulative effects of the random perturbations given by $f_s$ and $f_b$. The two following propositions describe the behavior of the decay rates \eqref{decayrate} in the limit of a large number of propagating modes $\N{}\gg 1$ corresponding to $  k(\omega) \gg 1$ (see \eqref{number}).
\begin{prop} \label{coefdec} 
We have the following asymptotic behaviors.
\begin{enumerate}
\item For $j=[\nu\N{}^{\eta_1}]$ with $\eta_1\in[0,1/2)$ and $\nu> 0$, we have
\[  \Gamma^{c,b}_{jj}(\omega)\underset{k(\omega)\gg1}{\sim} -   k^{3/2}(\omega) \frac{j^2/\N{}^2}{(1-\theta^2j^2/\N{}^2)^{1/4}}  \frac{2^{9/2}n_1^{5/2} \theta^2}{\pi d }\int_{0}^{+\infty}\sqrt{v}I_b(v)dv. \]
\item For $j=[\nu\N{}^{1/2}]$ and $\nu> 0$, we have
\[  \Gamma^{c,b}_{jj}(\omega)\underset{k(\omega)\gg1}{\sim} -   k^{3/2}(\omega) \frac{j^2/\N{}^2}{\sqrt{1-\theta^2j^2/\N{}^2}}  \frac{2^4 n_1^{5/2} \theta^2}{\pi d }\int_{-\nu^2 \pi \theta/(2d)}^{+\infty}\sqrt{\nu^2 \pi \theta/d+ 2v}I_b(v)dv. \]
\item  For $j=[\nu\N{}^{\eta_2}]$ with $\eta_2\in(1/2,1]$ and $\nu> 0$  ($\nu\leq 1$ if $\eta_2=1$), we have
\[  \Gamma^{c,b}_{jj}(\omega)\underset{k(\omega)\gg1}{\sim} - k^{2}(\omega) \frac{j^3/\N{}^3}{\sqrt{1-\theta^2j^2/\N{}^2}}\frac{2^5 n_1^3\theta^3}{\pi d}\int_{0}^{+\infty}I_b(v)dv.\]
\item For $j=[\nu\N{}^{\eta_3}]$ with $\eta_3\in[0,1]$ and $\nu> 0$ ($\nu\leq 1$ if $\eta_3=1$), we have
\[ \Gamma^{1,b}_{jj}(\omega)\underset{k(\omega)\gg1}{\sim}  k^{2}(\omega) \frac{j^4/\N{}^4}{1-\theta^2j^2/\N{}^2} \frac{4n_1^2\theta ^4}{d^2}\int_{0}^{+\infty}R_b(v)dv.\]
\item For $j=[\nu\N{}^{\eta_4}]$ with $\eta_4\in[0,1]$ and $\nu> 0$ ($\nu < 1$ if $\eta_4=1$), we have for $k(\omega)\gg 1$
\[\frac{C_1}{k(\omega)^{\alpha_{I_b}-3}} \frac{j^2/\N{}^2}{\sqrt{1-\theta^2j^2/\N{}^2}}\leq \Lambda^{c,b}_j(\omega)\leq \frac{C_2}{k(\omega)^{\alpha_{I_b}-3}}\frac{j^2/\N{}^2}{\sqrt{1-\theta^2j^2/\N{}^2}},\]
where $C_1$ and $C_2$ are two positive constants.
 \item For $j=\N{}-[\nu]$ with $\nu> 0$, we have for $k(\omega)\gg 1$
\[C_3 k(\omega)^{5/2} \leq \Lambda^{c,b}_j(\omega)\leq C_4 k(\omega)^{5/2},\]
where $C_3$ and $C_4$ are two positive constants.
\end{enumerate}
\end{prop}
The proof of this proposition is given in Section \ref{proofcoefdec}. Regarding the influence of the surface fluctuations we have the following result.
\begin{prop} \label{coefdec2} 
We have the following asymptotic behaviors.
\begin{enumerate}
\item For $j=[\nu\N{}^{\eta_1}]$ with $\eta_1\in[0,1/2)$ and $\nu> 0$, we have
\[\begin{split}  \Gamma^{c,s}_{jj}(\omega)\underset{k(\omega)\gg1}{\sim} - C_1 k^{3/2}(\omega) &\frac{ j^2}{\N{}^2\pi d}\int_{0}^{+\infty}\sqrt{v}I_s(v)dv\\
&\times  \Big[\big(1-\frac{ \theta^2j^2}{\N{}^2}\big)^{-1/2}+\big(1-\frac{ \theta^2j^2}{\N{}^2}\big)^{-5/4}+\big(1-\frac{ \theta^2j^2}{\N{}^2}\big)^{-9/4} \Big].
\end{split}\]
\item For $j=[\nu\N{}^{1/2}]$ and $\nu> 0$, we have
\[\begin{split}  \Gamma^{c,s}_{jj}(\omega)\underset{k(\omega)\gg1}{\sim} - C_2 k^{3/2}(\omega) &\frac{ j^2}{\N{}^2\pi d} \int_{-\nu^2 \pi \theta/(2d)}^{+\infty}\sqrt{\nu^2 \pi \theta/d+ 2v}I_b(v)dv\\
&\times \Big[(1-\frac{ \theta^2j^2}{\N{}^2})^{-1/2}+(1-\frac{ \theta^2j^2}{\N{}^2})^{-3/2}+(1-\frac{ \theta^2j^2}{\N{}^2})^{-5/2} \Big].
\end{split}\]
\item  For $j=[\nu\N{}^{\eta_2}]$ with $\eta_2\in(1/2,1]$ and $\nu> 0$  ($\nu\leq 1$ if $\eta_2=1$), we have
\[\begin{split}  \Gamma^{c,s}_{jj}(\omega)\underset{k(\omega)\gg1}{\sim} - C_3 k^{2}(\omega)& \frac{ j^3}{\N{}^3\pi d} \int_{0}^{+\infty}I_s(v)dv\\
&\times \Big[\big(1-\frac{ \theta^2j^2}{\N{}^2}\big)^{-1/2}+\big(1-\frac{ \theta^2j^2}{\N{}^2}\big)^{-3/2}+\big(1-\frac{ \theta^2j^2}{\N{}^2}\big)^{-5/2} \Big].
\end{split}\]
 \item For $j=[\nu\N{}^{\eta_3}]$ with $\eta_3\in[0,1]$ and $\nu> 0$ ($\nu < 1$ if $\eta_4=1$), we have for $k(\omega)\gg 1$
\[\frac{C_1}{k(\omega)^{\alpha_{I_s}-3}}\frac{j^2/\N{}^2}{\sqrt{1-\theta^2j^2/\N{}^2}}\leq \Lambda^{c,s}_j(\omega)\leq \frac{C_2}{k(\omega)^{\alpha_{I_s}-3}}\frac{j^2/\N{}^2}{\sqrt{1-\theta^2j^2/\N{}^2}},\]
where $C_1$ and $C_2$ are two positive constants.
 \item For $j=\N{}-[\nu]$ with $\nu> 0$, we have for $k(\omega)\gg 1$
\[C_3 k(\omega)^{5/2} \leq \Lambda^{c,s}_j(\omega)\leq C_4 k(\omega)^{5/2},\]
where $C_3$ and $C_4$ are two positive constants.
\end{enumerate}
\end{prop}
The proof of this proposition follows the lines of the one of Proposition \ref{coefdec} but with lengthier computations. For the points $1-3$ the three terms into brackets correspond respectively to the three terms in \eqref{defQ}. From these two propositions we can first remark that the surface and bottom fluctuations affect the amplitude of the propagating mode in the same way. Second, the decay rates become larger as the order of the propagating mode $j$ increase. The reason is that the higher the propagating mode $j$ is the more it bounces on the randomly perturbed boundaries, and therefore the more it is scattered. However, while $j$ is not of order $\N{}$ the amplitude attenuation is mainly produced by the mode coupling between the propagating modes themselves through $\Gamma^{c,s}_{jj}$ and $\Gamma^{c,b}_{jj}$. Now, when $j\sim \N{}$, the modes can coupled significantly with the radiating modes (see \eqref{coefatt}, \eqref{coefatt2}, and \eqref{spectrumRP2}) through $\Lambda^{c,s}_{j}$ and $\Lambda^{c,b}_{j}$, which produce more important losses in the bottom of the waveguide. In fact, for these modes the decay rate produced by the mode coupling between the propagating modes is of order $k^2(\omega)$ (points 3-4 in Proposition  \ref{coefdec} and point 3 in Proposition \ref{coefdec2}), while the one produced by the mode coupling with the radiating modes is of order $k^{5/2}(\omega)$ which is significantly larger.

\section*{Conclusion}

In this paper we have analyzed wave propagation in an acoustic waveguide with a randomly perturbed free surface and uneven topography, and the resulting coupling mechanism between the three kinds of modes (propagating, radiating, and evanescent). 
We have shown that the evolution of the forward propagating mode amplitudes can be described in term of a diffusion process (Theorem \ref{thasymptP21}, Theorem \ref{thasymptP12}, and Theorem \ref{thasympP2}) taking into account the main coupling mechanisms: the coupling with the evanescent modes induces a mode-dependent and frequency-dependent phase modulation on the propagating modes, the coupling with the radiating modes, in addition to a mode-dependent and frequency-dependent phase modulation, induces a mode-dependent and frequency-dependent attenuation on the propagating modes. Moreover, we have observed (Proposition \ref{coefdec} and Proposition \ref{coefdec2}) that the surface and bottom fluctuations affect the amplitude of the forward propagating modes in the same way. However, the amplitudes of highest propagating mode is more affected because of an efficient coupling with the radiating modes.

\section*{Acknowledgment}
 
This work was supported by AFOSR FA9550-10-1-0194 Grant. 
 
\section{Appendix}\label{appendix}

\subsection{Spectral Decomposition in Unperturbed Waveguides}\label{spectralP2}

This section is devoted to the presentation of the spectral decomposition of the Pekeris operator $\partial ^2 _x +\ko{}n^2 (x)$, where the index of refraction $n(x)$ is defined by \eqref{indexrefraction2}. Here, $H=L^2(0,+\infty)$ is equipped with the inner product defined by
\[\forall (h_1,h_2)\in H\times H,\quad \big<h_1,h_2\big>_H=\int_0^{+\infty}h_1(x)\overline{h_2(x)}dx.\]
\begin{defi}
Let us note by $R(\omega)$ the Pekeris operator defined by 
\begin{equation}\label{pekerisop}R(\omega)(y)=\frac{d^2}{ dx^2} y +\ko{}n^2 (x)y\quad \forall y\in \mathcal{D}(R(\omega)),\end{equation}
which is an unbounded operator on $H$ with domain
\[ \mathcal{D}(R(\omega))= H^1 _0 (0,+\infty)\cap H^2(0,+\infty).\]
\end{defi}
The following result regarding the spectral decomposition of the Pekeris operator  $R(\omega)$ is proved in \cite{wilcox}.
\begin{thm}
$R(\omega)$ is a self-adjoint operator on $H$, and its spectrum is given by
\begin{equation}\label{spectrumRP2}
Sp\big(R(\omega)\big)=\left(-\infty,\ko{}\right]\cup\big\{ \beta^2 _{\N{}}(\omega),\dots,\beta_1^2 (\omega)\big\},\end{equation}
where for all $j\in\big\{1,\dots,\N{}\big\}$, the modal wavenumbers $\beta_j (\omega)$ are positive and ordered in an decreasing way : 
\[\ko{}<\beta^2 _{\N{}}(\omega)<\cdots<\beta_1^2 (\omega)<n_1^2 \ko{}.\]
Moreover, let $\Pi_\omega$ be the resolution of the identity associated to $R(\omega)$, we have for all $ y\in H$ and for all $ r\in\mathbb{R}$,

\begin{equation}\label{dec1}\begin{split}
 \Pi_\omega(r,+\infty)(y) (x) =&  \sum_{j=1}^{\N{}} \big<y,\phi_j(\omega,.)\big>_H \phi_j(\omega,x)\textbf{\emph{1}}_{(r,+\infty)}\left(\Bh{j}{}^2\right) \\
                                                                      & + \int_{r}^{\ko{}}  \big<y,\phi_\ga(\omega,.) \big>_H \phi_\ga (\omega,x)d\ga \textbf{\emph{1}}_{\left(-\infty,\ko{} \right)}(r),
  \end{split}\end{equation}
and for all $y\in \mathcal{D}(R(\omega))$,
  \begin{equation}\label{dec2}\begin{split}                                                                    
 \Pi_\omega(r,+\infty)(R(\omega)(y)) (x) =&  \sum_{j=1}^{\N{}}\Bh{j}{}^2 \big<y,\phi_j(\omega,.)\big>_H \phi_j(\omega,x)\textbf{\emph{1}}_{(r,+\infty)}\left(\Bh{j}{}^2\right) \\
                                                                      & + \int_{r}^{\ko{}}  \ga \big<y,\phi_\ga(\omega,x) \big>_H \phi_\ga (\omega,x)d\ga \textbf{\emph{1}}_{\left(-\infty,\ko{} \right)}(r).
\end{split}\end{equation}   
\end{thm}
Let us describe more precisely the decompositions \eqref{dec1} and \eqref{dec2}.

\paragraph{Discrete part of the decomposition}

For all $j\in \big\{1,\dots,\N{}\big\}$, the $j$th eigenvector is given \cite{wilcox} by 
\begin{equation}\label{phij}\phi_j(\omega, x)=\left\{ \begin{array}{ccl}
A_j(\omega)\sin(\sigma_j(\omega) x/d) & \mbox{ if } & 0\leq x \leq d \\
A_j(\omega)\sin(\sigma_j(\omega)  )e^{-\zeta_j (\omega) \frac{x-d}{d}}& \mbox{ if } & d\leq x,  \end{array} \right.\end{equation}
where
\begin{equation}\label{sigj}\sigma_j (\omega)=d \sqrt{n_1 ^2 \ko{}-\beta^2 _j( \omega)}, \quad \zeta_j(\omega) =d \sqrt{\Bh{j}{}^2-\ko{}},\end{equation}
and 
\begin{equation}\label{Aj}
A_j(\omega)=\sqrt{\frac{2/d}{1+\frac{\sin^2 (\sigma_j (\omega))}{\zeta_j (\omega) }-\frac{\sin(2 \sigma_j (\omega)  )}{2\sigma_j (\omega) } }}.
\end{equation}
According to \cite{wilcox},  $\sigma_1 (\omega),\dots,\sigma_{\N{}}(\omega)$ are the solutions on $(0,n_1k(\omega)d\theta)$ of the equation 
\begin{equation}\label{eqsigma}
\tan(y )=-\frac{y }{\sqrt{( n_1 k d  \theta) ^2-y^2}},
\end{equation}
such that $0<\sigma_1 (\omega)<\cdots<\sigma_{\N{}}(\omega)< n_1 k(\omega)d \theta$,  and with  
\begin{equation}\label{deftheta}\theta =\sqrt{1-1/n_1 ^2}=\sqrt{1-c_1^2/c_0 ^2}.\end{equation}
 This last equation admits exactly one solution over each interval of the form $\big(\pi/2+(j-1)\pi, \pi/2 +j\pi \big)$ for $ j \in \{1,\dots, N(\omega) \}$, so that the number of eigenvectors is 
\begin{equation}\label{number}\N{}=\left[ \frac{n_1  k(\omega)d}{\pi} \theta \right],\end{equation}
where $[\cdot]$ stands for the integer part.

\paragraph{Continuous part of the decomposition}

For $\ga \in (-\infty, \ko{})$, we have \cite{wilcox}
\begin{equation}\label{phiga}
\begin{split}
\phi_\ga& (\omega, x)=\\ 
&\left\{ 
\begin{array}{ccl}
A_\gamma(\omega) \sin(\eta(\omega)  x/d ) & \mbox{ if } & 0\leq x \leq d \\
A_\gamma(\omega) \left(\sin(\eta(\omega)  )\cos\big(\xi(\omega) \frac{x-d}{d}\big)+\frac{\eta(\omega) }{\xi(\omega) }\cos(\eta(\omega) )\sin\big(\xi(\omega) \frac{x-d}{d}\big)\right)& \mbox{ if } & d\leq x, \end{array}
 \right.
\end{split}
\end{equation}
where
\begin{equation}\label{etaga}\eta (\omega) =d\sqrt{n_1 ^2 \ko{}-\gamma }, \quad \xi(\omega)  =d\sqrt{\ko{}-\gamma },\end{equation}
and
\[A_\gamma(\omega) =\sqrt{\frac{d \xi(\omega) }{\pi\big(\xi ^2(\omega) \sin^{2}(\eta(\omega))+\eta ^2(\omega) \cos^2 (\eta(\omega))\big)}}.\]
Let us note that we have
\begin{equation}\label{aga} A_\ga(\omega)\underset{\ga\to-\infty}{\sim}\frac{1}{\sqrt{\pi}\lvert \ga \rvert^{1/4}}.\end{equation}
We can also remark that $\phi_\ga(\omega, .)$ does not belong to $H$, so that $\big<y,\phi_\ga (\omega,.)\big>_H$ in \eqref{dec1} and \eqref{dec2} is not defined in the classical sense, but in the following way
\[ \big<y,\phi_\ga (\omega,.) \big>_H=\lim_{M\to+\infty}\int_{0}^M y(x) \phi_\ga (\omega, x)dx\quad  \text{ in }L^2\big(-\infty,\ko{}\big). \]
As a result, we have
\[ \|y \|^2 _H = \sum_{j=1}^{\N{}}\big \lvert \big<y,\phi_j(\omega,.) \big>_H \big \rvert^2+\int_{-\infty}^{\ko{}}\big\lvert \big<y,\phi_\ga (\omega,.) \big>_H\big\rvert^2 d\ga,\]
and then,
\[
\begin{array}{rcc}
\Theta_\omega:H&\longrightarrow &\mathcal{H}^\omega\\
y& \longrightarrow & \Big(\big(\big<y,\phi_j(\omega,.)\big>_H\big)_{ j=1,...,\N{}},\big(\big<y,\phi_\ga (\omega,.)\big>_H\big)_{\ga\in(-\infty,\ko{})}\Big)
\end{array}
\]
is an isometry, from $H$ onto $\esp=\mathbb{C}^{\N{}}\times L^2\big(-\infty,\ko{}\big)$.

In this paper the pressure field \eqref{modedec} can be decomposed according to the resolution of the identity $\Pi_\omega$ introduced in this section.
\[
\widehat{p}_0(\omega,u,v)=\underbrace{\sum_{j=1}^{\N{}} \widehat{p}_j (\omega,u )\phi_j(\omega,v)}_{\text{propagating modes}}+\underbrace{\int_{0}^{\ko{}}\widehat{p}_\ga (\omega,u)\phi_\ga (\omega,v)d \ga}_{\text{radiating modes}}+\underbrace{\int_{-\infty}^{0}\widehat{p}_\ga (\omega,u)\phi_\ga (\omega,v)d \ga}_{\text{evanescent modes}},
\]
where the amplitudes of the propagating modes are defined by
\begin{equation}\label{defmode1} 
\forall j\in \big\{1,\dots,\N{}\big\}, \quad \widehat{p}_j (\omega,u )=\Pi_\omega( \{ \Bh{j}{} \} )(\widehat{p}_0(\omega,u,\cdot))=\int_{0}^{+\infty} \widehat{p}_0(\omega,u,v) \phi_\ga (\omega, v)dv,
\end{equation}
and the amplitudes of the radiating and evanescent modes are defined by
\begin{equation}\label{defmode2} 
\forall \ga\in(-\infty,\ko{}), \quad \widehat{p}_\ga (\omega,u)=\lim_{M\to+\infty}\int_{0}^M \widehat{p}_0(\omega,u,v) \phi_\ga (\omega, v)dv \quad  \text{ in }L^2\big(-\infty,\ko{}\big).
\end{equation}
Consequently, the mode amplitudes are the projections of the pressure field $\widehat{p}_0(\omega,u,v)$ over the eigenelements of the Pekeris operator.

\subsection{Proof of Proposition \ref{propconfmap2}}\label{proofpropconfmap}

As we have remarked right after the statement of Proposition \ref{propconfmap2} the conformal transformation $\Phi$ can be extended as a homeomorphism from the closure of $D_0$ onto the closure of $D_{V_s}$. Consequently, letting $v\to0$ in the definition of $z(u,v)$ in \eqref{conformmap}, we obtain the following equation 
\[z(u,0)=u-\mathcal{U}(V^{\e}_s(z(\cdot,0)))(u),\]
which is well defined since $V^{\e}_s$ has a compact support and $z(u,0)$ goes to $+\infty$ (reps. $-\infty$) as $u$ goes to $+\infty$ (reps. $-\infty$). This equation also can be recast as follows
\begin{equation}\label{invert}q(u)-Q(q)(u) =\tilde{q}(u) \end{equation}
where $q(u)=z(u,0)-u$, $\tilde{q}(u)=-\mathcal{U}(V^{\e}_s)(u)$, and 
\[
Q(q)(u)=\frac{1}{\pi}\text{p.v.} \int \frac{V^{\e}_s(\tilde{u}+q(\tilde{u}))-V^{\e}_s(\tilde{u})}{u-\tilde{u}}d\tilde{u}.
\]
The goal is to invert \eqref{invert} in $L^p(\mathbb{R})$ ($p>2$) by showing that $Id-Q$ is homeomorphism, where $Id$ stands for the identity map. However, to show that we need to know first if $q\in L^p(\mathbb{R})$. In fact, for all $ v>0$ one can show that 
\[ \int \lvert (\Phi-Id)(u+iv)\rvert^p du<+\infty, \]
where $\Phi$ is the conformal transform, so that according to \cite[Theorem 93]{titch} we have $q\in L^p(\mathbb{R})$. Moreover,
\[\|Q(q)\|_{L^p(\mathbb{R})} \leq \se K \| q \|_{L^p(\mathbb{R})} \quad \text{and}\quad  \| \tilde{q} \|_{L^p(\mathbb{R})} \leq \e^{1/2-1/p} K \]
thanks to the boundedness of the Hilbert transform operator in $L^p(\mathbb{R})$, so that $\|q\|_{L^p(\mathbb{R})}\leq K$ thanks to \eqref{invert} for $\e$ small enough. Here, $K$ denotes a deterministic constant which is independent of $\e$. Now, let us consider $\Theta$ the application from $L^p(\mathbb{R})$ to itself defined by
\[ \Theta (l)=\sum_{n\geq 0} \underbrace{Q\circ\dots\circ Q(l)}_{n\text{ times}}.\]
This application is well defined since $L^p(\mathbb{R})$ is a Banach space and $\|Q(l)\|_{L^p(\mathbb{R})} \leq \se K \| l \|_{L^p(\mathbb{R})}$, and we have
\[(Id-Q)\circ \Theta =Id\quad\text{and}\quad\Theta\circ(Id-\Theta)=Id.\] 
Moreover, one can easily check that both $(Id-Q)$ and $\Theta$ are Lipschitz function, so that  $(Id-Q)$ is an homeomorphism. As a result, we obtain from \eqref{invert} that 
\[ q(u) = \Theta(\tilde{q})(u) = \tilde{q}(u) + Q(\Theta(\tilde{q}))(u)\quad\text{with}\quad \|Q(\Theta(\tilde{q}))\|_{L^p(\mathbb{R})}\leq \e^{1-1/p} K,\]
which concludes the proof of this proposition.

\subsection{Proof of Proposition \ref{asymptconf}} \label{proofasymptconf}

First, under Assumption \ref{hyp1} and thanks to the Holder's inequality, we have
\[\begin{split}
\lvert x(u,v)-v\rvert&\leq \Big(\frac{v}{\pi}\int \frac{\lvert x(\tilde{u},0)-x(u,0)\rvert}{(u-\tilde{u})^2+v^2}d\tilde{u}+ \lvert x(u,0)\rvert \Big)\1_{(v<\mu)}\\
&+v\Big[\int \lvert x(\tilde{u},0) \rvert^pd\tilde{u}\Big]^{1/p}\Big[\int \frac{\tilde{u}}{((u-\tilde{u})^2+v^2)^q}d\tilde{u} \Big]^{1/q}\1_{(v>\mu)}\\
&\leq \e^{1/2-1/p} C_1 \Big(\frac{v}{\pi}\int_{\lvert \bar{u} \rvert>1} \frac{d\bar{u}}{\bar{u}^2}+v \int_{\lvert \bar{u} \rvert<1/v} \frac{\lvert \bar{u}\rvert }{\bar{u}^2+1}\bar{u}+ 1 \Big)\1_{(v<\mu)}+\e^{1/2-1/p}\frac{C_2}{v^{1-1/q}} \1_{(v>\mu)}\\
&\leq \e^{1/2-1/p} C_\mu,
\end{split}\]
since $p>2$, $\sup_{u\geq 0}\,\lvert \mathcal{U}(V^\e_s)(u)\rvert+\lvert \mathcal{U}((d/du)V^\e_s)(u)\rvert\leq C \e^{1/2-1/p}$,
and
\[\begin{split}
\lvert z(u,v)-u\rvert&\leq \Big(\frac{1}{\pi}\int_{\lvert u-\tilde{u}\rvert>1} \frac{\lvert u-\tilde{u}\rvert \lvert x(\tilde{u},0)\rvert}{(u-\tilde{u})^2+v^2}d\tilde{u}+\frac{1}{\pi}\int_{\lvert\bar{u}\rvert<1} \frac{\lvert \bar{u}\rvert \lvert x(u+\bar{u},0)-x(u-\bar{u},0)\rvert}{\bar{u}^2+v^2} d\bar{u}\Big)\1_{(v<\mu)}\\
&+\Big[\int \lvert x(\tilde{u},0) \rvert^pd\tilde{u}\Big]^{1/p}\Big[\int\lvert \frac{u-\tilde{u}}{(u-\tilde{u})^2+v^2} \rvert ^qd\tilde{u}\Big]^{1/q}\1_{(v>\mu)}\\
&\leq \e^{1/2-1/p} C_1 \Big(\Big[\int_{\lvert \bar{u} \rvert>1} \frac{\bar{u}}{\bar{u}^q}d\bar{u}\Big]^{1/q}+v \int_{\lvert \bar{u} \rvert<1/v} \frac{\lvert \bar{u}\rvert^2 }{\bar{u}^2+1}d\bar{u} \Big)\1_{(v<\mu)}+\e^{1/2-1/p}\frac{C_2}{v^{1-1/q}} \1_{(v>\mu)}\\
&\leq \e^{1/2-1/p} C_\mu.
\end{split}\]

For the second part of the proof, for all $\mu>0$, we have
 \[\begin{split}
 \int_\mu^{+\infty} \E\big[ \lvert z(u,v)-u\rvert^2 \big]&\leq 2\e\int_\mu^{+\infty} dv \int d\bar{u} \int d\tilde{u} R_s(\bar{u}-\tilde{u}) \frac{f_s(\e \bar{u})f_s(\e \tilde{u})(u-\bar{u})(u-\tilde{u})}{((u-\bar{u})^2+v^2)((u-\tilde{u})^2+v^2)}\\
 &+ \e^2 C  \int_\mu^{+\infty} dv\E\Big[\Big\lvert \int_0^{L/\e} d\tilde{u} \frac{\tilde{u}\,\, \mathcal{U}(V^\e_s)(\tilde{u})}{\tilde{u}^2+v^2} \Big\rvert^2\Big]\\
 &\leq 2 \e\int_0^{L/\e}dw R_s(w) \int_\mu^{+\infty} dv\frac{1}{v}  \int_w^{L/\e} d\tilde{u}\frac{\tilde{u}}{\tilde{u}^2+v^2}\\
 &+ \e^{2(1-1/q)} C \Big[\int_0^L \lvert f_s(\tilde{u})\rvert^q d\tilde{u} \Big]^{2/q} \Big[\int_0^{+\infty} \Big\lvert \frac{\tilde{u}}{\tilde{u}^2+1}\Big\rvert^p d\tilde{u} \Big]^{2/p} \int_\mu^{+\infty}dv \frac{1}{v^{1-2/p}}\\
 &\leq \e \int_0^{+\infty}dw R_s(w) \int_\mu^{+\infty} dv \frac{1}{v}  \ln\Big(\frac{L}{\e^2 v^2}+1\Big) + \e^{2(1-1/q)} C, 
 \end{split} \]
with $p>2$ and $1<q<2$, and $R_s$ is the autocorrelation function of the surface fluctuations $V_s$ introduced in \eqref{randompert}. However, we have
\[ \e  \int_\mu^{+\infty} dv \frac{1}{v}  \ln\Big(\frac{L}{\e^2 v^2}+1\Big)\leq \e  \int_1^{+\infty} dv \frac{1}{v}  \ln\Big(\frac{L}{ v^2}+1\Big) + \e  \int_{\e\mu}^{1} dv \frac{1}{v}  \ln\Big(\frac{L}{ v^2}+1\Big)\leq C(\e+\e^{1-\nu}), \]
for some $\nu\in(0,1)$. Now, we have
\[\begin{split}
 \int_\mu^{+\infty} \E\big[ \lvert x(u,v)-v\rvert^2 \big]&\leq2 \e\int_\mu^{+\infty} dv \int d\bar{u} \int d\tilde{u} R_s(\bar{u}-\tilde{u}) \frac{f_s(\e \bar{u})f_s(\e \tilde{u})v^2}{((u-\bar{u})^2+v^2)((u-\tilde{u})^2+v^2)}\\
  &+ \e^2 C  \int_\mu^{+\infty} dv\, v^2\E\Big[\Big\lvert \int_0^{L/\e} d\tilde{u} \frac{\mathcal{U}(V^\e_s)(\tilde{u})}{(u-\tilde{u})^2+v^2} \Big\rvert^2\Big]\\
 &\leq 2 \e \int_0^{L/\e}dw R_s(w) \int_\eta^{+\infty} dv \int_w^{L/\e} d\tilde{u}\frac{1}{\tilde{u}^2+v^2}\\
&+ \e^{2(1-1/q)} C \Big[\int_0^L \lvert f_s(\tilde{u})\rvert^q d\tilde{u} \Big]^{2/q} \Big[\int_0^{+\infty} \frac{1}{(\tilde{u}^2+1)^p} d\tilde{u} \Big]^{2/p} \int_\mu^{+\infty}dv \frac{1}{v^{1-2/p}},\\
 &\leq  2 \e \int_0^{+\infty}dw R_s(w) \int_\eta^{+\infty} dv \frac{1}{v} \arctan\Big(\frac{L}{\e v}\Big)+\e^{2(1-1/q)} C,
 \end{split}\]
with $p>2$ and $1<q<2$, and 
\[ \e  \int_\mu^{+\infty} dv \frac{1}{v}  \arctan\Big(\frac{L}{\e v}\Big)\leq \e  \int_1^{+\infty} dv \frac{1}{v}  \arctan\Big(\frac{L}{ v}\Big) + \e  \int_{\e\mu}^{1} dv \frac{1}{v}  \arctan\Big(\frac{L}{ v}\Big)\leq C\Big(\e+\e\ln\Big(\frac{1}{\e}\Big)\Big), \]
which concludes the proof of Proposition  \ref{asymptconf} thanks to the first part of the proof.

\subsection{Proof of Proposition \ref{propresidual}}\label{proofpropresidual}

The proof given in this section follows closely the one proposed in \cite[Section 4.3]{gomez2}, but we still present the main ideas for the sake of completeness. First of all, we recall that $\Theta_\omega\circ\Pi_\omega(-1/\xi,-\xi)\big(\widehat{p}_0(\omega,.,u)\big)$ represents the evanescent part of the pressure field $\widehat{p}_0(\omega,.,u)$, where $\Theta_\omega$ and $\Pi_\omega$ are defined in Section \ref{spectralP2}. In this section we consider the Banach space $F=L^1 (-1/\xi,-\xi)$
equipped with the norm $\|y\|_F=\int_{-1/\xi}^{-\xi}\lvert y_\ga\rvert d\ga$.
According to \eqref{evaeqdiff} and the radiation condition given by Assumption \ref{hyp6}, we have
\begin{equation}\label{exp1} \widehat{p}_\ga(\omega,u)=\frac{1}{2\sgaa}\int g_\ga(\omega,u+v)e^{-\sgaa\lvert v \rvert}dv+\tilde{p}_{1,\ga} (\omega,u),\end{equation}
where $g_\ga(\omega,u)$ is defined by \eqref{g}, and
\[\tilde{p}_{1,\ga} (\omega,u)=\frac{\widehat{a}^{\e,\xi}_{\ga,0}(\omega) }{\gaa^{1/4}}e^{-\sgaa u}\1_{(-1/\xi,-\xi)}(\ga).\]
Therefore, substituting \eqref{g} into \eqref{exp1}, we obtain the following relation in $\big(\mathcal{C}\big([0,+\infty),F\big),\|.\|_{\infty,F}\big)$
\begin{equation}\label{eqevmP2}
(Id- \Phi^\omega)\Big(\Theta_\omega\circ\Pi_\omega(-1/\xi,-\xi)\big(\widehat{p}_0(\omega,.,.)\big)\Big)=\tilde{p}(\omega,.)+\tilde{p}_1 (\omega,.),
\end{equation}
where $\|y\|_{\infty,F}=\sup_{z\geq 0} \|y(z)\|_F$, for all $y\in \mathcal{C}\big([0,+\infty),F\big)$.
In \eqref{eqevmP2}, $\Phi^\omega$ is a linear bounded operator, from $\big(\mathcal{C}\big([0,+\infty),F\big),\|.\|_{\infty,F}\big)$ to itself,
defined by
\[
\Phi^\omega_\ga (y)(u)=\frac{\ko{}}{2\sgaa}\int \int_{-1/\xi}^{-\xi} C^\e_{\ga \ga'}(\omega,u+v)y_{\ga'} (u+v)d\ga' e^{-\sgaa \lvert v\rvert }dv
\]
and 
\[
\tilde{p}_\ga(\omega,u)=\frac{\ko{}}{2\sgaa}\int\sum_{l=1}^{\N{}}\frac{C^\e _{\ga l}(\omega,u+v)}{\sqrt{\beta_l}}\big(\widehat{a}_l (\omega,u+v)e^{i\beta_l (u+v)}+\widehat{b}_l (\omega,u+v)e^{-i\beta_l (u+v)}\big)e^{-\sgaa \lvert v\rvert}du\
\]
for all $u \in [0,+\infty)$ and for almost every $\ga \in (-1/\xi,-\xi)$.
We can check that the norm of the operator $\Phi^\omega$ is bounded by
\[ \| \Phi^\omega\|\leq C \frac{\e^{3/2-1/p}}{\xi}\big(  \sup_{u}\lvert V_s(u) \rvert+ \sup_{u}\lvert \frac{d}{du}V_s(u) \rvert+ \sup_{u}\lvert V_b(u) \rvert  \big), \]
so that $\lim_{\e\to 0}\mathbb{P}(Id- \Phi^\omega \text{ is invertible})=1$, and the condition ($Id-\Phi^\omega$ is invertible) is satisfied in the asymptotic $\e \to 0$. As a result, on the event ($Id- \se \,\Phi^\omega$ is invertible), we have
\[\begin{split}
\Theta_\omega\circ\Pi_\omega(-1/\xi,-\xi)&\big(\widehat{p}_0(\omega,.,.)\big)=\big(Id-\Phi^\omega\big)^{-1}(\tilde{p}(\omega, .)+ \tilde{p}_{1} (\omega,.))\\
&=\tilde{p}(\omega, .)+ \tilde{p}_{1} (\omega,.)+\Phi^\omega(\tilde{p}_{1} (\omega,.))\\
& \quad+\sum_{j=1}^{+\infty}(\Phi^\omega)^j \big(\tilde{p}(\omega,.)+\Phi^\omega(\tilde{p}_{1} (\omega,.)) \big),
\end{split}\]
so that
\begin{equation}\label{expensioneva}\begin{split}\Theta_\omega&\circ\Pi_\omega(-1/\xi,-\xi)\big(\widehat{p}_0(\omega,.,.)\big)=  \tilde{p}(\omega, .) + \tilde{p}_{1} (\omega,.)+\Phi^\omega(\tilde{p}_{1} (\omega,.))\\
&+\mathcal{O}\Big( \frac{\e^{3/2-1/p}}{\xi}\big(  \sup_{u}\lvert V_s(u) \rvert+ \sup_{u}\lvert \frac{d}{du}V_s(u) \rvert+ \sup_{u}\lvert V_b(u) \rvert  \big) \sup_{u\in[0,L /\e]}\|\widehat{a}(\omega,u) \|_{\esp_\xi}+\|\widehat{b}(\omega,u) \|_{\esp_\xi}\Big)
\end{split}\end{equation}
in  $\mathcal{C}\big([0,+\infty),F\big)$. Moreover, introducing
\[
\tilde{p}_{\ga,2}(\omega,u)=\frac{\ko{}}{2\sgaa}\int  \sum_{l=1}^{\N{}}\frac{C^\e_{\ga l}(\omega, u+v)}{\sqrt{\beta_l}}\big(\widehat{a}_l (\omega,u)e^{i\beta_l (u+v)}+\widehat{b}_l (\omega,u)e^{-i\beta_l (u+v)}\big) e^{-\sgaa \vert v \rvert }dv
\]
for all $u \in [0,+\infty)$, we have according to \eqref{mcP21}-\eqref{mcP24}, 
\[\begin{split}
 \| \tilde{p}(\omega,.)- &\tilde{p}_2(\omega, .) \|_{\infty,F}\leq C \frac{\e^{3/2-1/p}}{\xi}\big(  \sup_{u}\lvert V_s(u) \rvert+ \sup_{u}\lvert \frac{d}{du}V_s(u) \rvert+ \sup_{u}\lvert V_b(u) \rvert  \big)^2\\
 &\times \Big( \sup_{z\in[0,L /\e]}\|\widehat{a}(\omega,z) \|_{\esp_\xi}+\|\widehat{b}(\omega,z) \|_{\esp_\xi} +\|\Theta_\omega \circ\Pi_\omega (-1/\xi,-\xi)\big(\widehat{p}_0(\omega,.,.)\big) \|_F \Big).
 \end{split}\]
Therefore, since $V_s$, $V_b$ and their derivative are assumed to be bounded, and according to \eqref{expensioneva}, we obtain
\[\begin{split} \Theta_\omega& \circ\Pi_\omega (-1/\xi,-\xi)\big(\widehat{p}_0(\omega,.,.)\big)= \tilde{p}_2(\omega, .) + \tilde{p}_{1} (\omega,.)+\Phi^\omega(\tilde{p}_{1} (\omega,.))\\
&+\mathcal{O}\Big( \frac{\e^{3/2-1/p}}{\xi}\big(  \sup_{u}\lvert V_s(u) \rvert+ \sup_{u}\lvert \frac{d}{du}V_s(u) \rvert+ \sup_{u}\lvert V_b(u) \rvert  \big)^2\sup_{z\in[0,L /\e]}\|\widehat{a}(\omega,z) \|_{\esp_\xi}+\|\widehat{b}(\omega,z) \|_{\esp_\xi} \Big)\end{split}\]
in  $\mathcal{C}\big([0,+\infty),F\big)$, which concludes the proof of Proposition  \ref{propresidual}.

\subsection{Proof of Theorem \ref{thasymptP21}}\label{proofth}

The proof of this theorem is based on a martingale approach using the perturbed-test-function method. However, the process $\big(\widehat{a}^{\xi,\e}(u)\big)_{u \geq 0}$ is not adapted with respect to the filtration $\mathcal{F}^\e _u=\mathcal{F}_{u/\e}$, where $\mathcal{F}_{u}$ is defined by \eqref{filtration}. In fact, neither the random operator $\textbf{H}^{aa}_{\e,\xi}$ nor $\textbf{G}^{aa}_{\e,\xi}$ (\eqref{haajP2}-\eqref{haagaP2} and \eqref{gaajP2}-\eqref{gaagaP2}) are adapted  with respect to the filtration $\mathcal{F}^\e _u$. The proof of this theorem is in three parts. The first part of the proof consists in approximating the random coefficients $\textbf{H}^{aa}_{\e,\xi}$ and $\textbf{G}^{aa}_{\e,\xi}$ with new ones from which we can use the mixing property.  The second part follows the ideas of \cite{khasminskii} and consists in simplifying the coupled mode equation and introducing a new process for which the martingale approach can be used. Finally, the third part of the proof is based on a martingale approach using the perturbed-test-function method and  follows the ideas developed in \cite{carmona, gomez2}.

\subsubsection{Approximation of the Random Coefficients}\label{approxrandomcoef}

This section is devoted to the approximation of the random operators $\textbf{H}^{aa}_{\e,\xi}$ and  $\textbf{G}^{aa}_{\e,\xi}$. In fact, these operators do not satisfy the mixing properties coming from the random perturbations of the waveguide surface. the reason is that the conformal map \eqref{conformmap} involves the complete trajectory of the random surface of the waveguide.

Let $\alpha_0\in(1/2,1)$, 
\begin{equation}\label{coef1/2}\begin{split}
C^{1/2,\e}_{rs}\Big(\frac{u}{\e}\Big)=&\,2(n_1^2-1)\phi_r(d)\phi_s(d)\frac{d}{\pi}\int _{\lvert w \rvert < 1/\e^{\alpha_0} }\frac{V_s(u/\e+w)f_s(u+\e w)}{w^2+d^2}dw\\
&-\frac{2}{\pi}\int_0^{+\infty}\phi_r (v)\phi_s (v)n^2(v)\int_{\lvert w\rvert < 1/\e^{\alpha_0}}\frac{V_s(u/\e+w)f_s(u+\e w)}{w^2+v^2}dwdv\\
&+ \frac{4}{\pi}\int_0^{+\infty}\phi_r (v)\phi_s (v)v^2n^2(v)\int_{\lvert w\rvert < 1/\e^{\alpha_0}}\frac{V_s(u/\e+w)f_s(u+\e w)}{(w^2+v^2)^2}dwdv,
\end{split}
\end{equation}
and 
\begin{equation}\label{filtration2}
\tilde{\mathcal{F}}_{s,t}=\sigma( C^{1/2,\e} (u),\quad s\leq u\leq t),
\end{equation}
the $\sigma$-algebra generated by $C^{1/2,\e}$. Let us remark that $\tilde{\mathcal{F}}_{s,t}$ is still $\tilde{\phi}_\e$-mixing. In fact, for $s>2/\e^{\alpha_0}$, if $A\in \tilde{\mathcal{F}}_{t+s,+\infty}$ and $B\in \tilde{\mathcal{F}}_{0,t}$, then  $A\in \mathcal{F}_{t+s-\e^{1-\alpha_0},+\infty} $ and  $B\in \mathcal{F}_{0,t+\e^{1-\alpha_0}} $,
and therefore we have
\[
\sup_{\substack{t\geq 0\\A\in\tilde{ \mathcal{F}}_{t+s,+\infty}\\B\in\tilde{\mathcal{F}}_{0,t}}}\lvert \mathbb{P}(A\vert B)-\mathbb{P}(A)\rvert\leq \phi(s-2/\e^{\alpha_0}) =\tilde{\phi}_\e(s).\]
According to \cite{kushner} we have the following results regarding mixing processes :
\begin{equation}\label{phimixing}
\big\lvert \E[ C^{1/2,\e}(t+s_1) \vert \tilde{\mathcal{F}}_{0,t} ] \big \rvert\leq 2 \tilde{\phi}_\e(s_1),\end{equation}
and
\begin{equation}\label{phimixing2}\begin{split}
\big\vert \E[ C^{1/2,\e}(t+s_1+s_2) C^{1/2,\e}(t+s_1) &\vert \tilde{\mathcal{F}}_{0,t} ]- \E[ C^{1/2,\e}(t+s_1+s_2) C^{1/2,\e}(t+s_2) ]  \big \vert\\
&\leq 4 \tilde{\phi}^{1/2}_\e(s_1)\tilde{\phi}^{1/2}_\e(s_2),
\end{split}\end{equation}
for all $ s_1,s_2 >2/\e^{\alpha_0}$. Let us also consider
\begin{equation}\label{coef1/22}\begin{split}
C&^{1,\e}_{rs}(u)=2(n_1^2-1)\phi_r(d)\phi_s(d)\Big[\frac{6d^3}{\pi^2}\int \frac{V_s(\tilde{u})f_s(\e \tilde{u})}{(\tilde{u}-u)^2+d^2}d\tilde{u}\int \frac{V_s(\tilde{u})f_s(\e \tilde{u})}{((\tilde{u}-u)^2+d^2)^2}d\tilde{u} \\
&\hspace{3cm}-\frac{3d}{\pi^2}\Big(\int \frac{V_s(\tilde{u})f_s(\e \tilde{u})}{(\tilde{u}-u)^2+d^2}d\tilde{u}\Big)^2+\frac{d}{\pi}\int \frac{\mathcal{U}(V_s f_s(\e\cdot))(\tilde{u})\frac{d}{dz}V_s(\tilde{u})f_s(\e \tilde{u})}{(\tilde{u}-u)^2+d^2}d\tilde{u}\Big]\\
&+(\phi_r(d)\phi'_s(d)+\phi'_r(d)\phi_s(d))\Big(\frac{d}{\pi}\int \frac{V_s(\tilde{u})f_s(\e \tilde{u})}{(\tilde{u}-u)^2+d^2} d\tilde{u}\Big)^2+\int_0^{+\infty} \phi_r(v)\phi_s(v)n^2(v)\\
& \hspace{2cm}\times \Big[\frac{4v^2}{\pi^2}\Big(\int \frac{(\tilde{u}-u)V_s(\tilde{u})f_s(\tilde{u})}{((\tilde{u}-u)^2+v^2)^2}d\tilde{u}\Big)^2-\frac{2}{\pi}\int \frac{\mathcal{U}(V_s f_s(\e \cdot))(\tilde{u})\frac{d}{dz}V_s(\tilde{u})f_s(\e \tilde{u})}{(\tilde{u}-u)^2+v^2}d\tilde{u}\\
&+\frac{4v^2}{\pi}\int \frac{\mathcal{U}(V_s f_s(\e \cdot))(\tilde{u})\frac{d}{dz}V_s(\tilde{u})f_s(\e \tilde{u})}{((\tilde{u}-u)^2+v^2)^2}d\tilde{u}-\frac{2v^2}{\pi^2}\int \frac{V_s(\tilde{u})f_s(\e \tilde{u})}{(\tilde{u}-u)^2+v^2}d\tilde{u}\int \frac{V_s(\tilde{u})f_s(\e \tilde{u})}{((\tilde{u}-u)^2+v^2)^2}d\tilde{u}\\
&+\frac{1}{\pi^2}\Big(\int \frac{V_s(\tilde{u})f_s(\e \tilde{u})}{(\tilde{u}-u)^2+v^2}d\tilde{u}\Big)^2+\frac{4v^4}{\pi^2}\Big(\int \frac{V_s(\tilde{u})f_s(\e \tilde{u})}{((\tilde{u}-u)^2+v^2)^2}d\tilde{u}\Big)^2\Big].
\end{split}\end{equation} 
We have the following result.
\begin{lem}\label{approxcoef} For all $(r,s)\in\{1,\dots,N\}\cup (\xi,k^2)$, we have
\[C^\e_{rs}(u)=\se \,\,C^{1/2,\e}_{rs}(u)+\e \,\,C^{1,\e}_{rs}(u)+o(\e),\]
uniformly with respect to $u\in(0,L/\e)$, and where $C^\e(u)$ is defined by \eqref{randomcoef}.
\end{lem}
This Lemma gives an approximation of the random operators $\textbf{H}^{aa}_{\e,\xi}$ and  $\textbf{G}^{aa}_{\e,\xi}$ (\eqref{haajP2}-\eqref{haagaP2} and \eqref{gaajP2}-\eqref{gaagaP2}) by new operators $\tilde{\textbf{H}}^{aa}_{\e,\xi}$ and  $\tilde{\textbf{G}}^{aa}_{\e,\xi}$ defined by 
 \[\begin{split}
 \tilde{\textbf{H}}^{aa}_{j,\e,\xi}(u)(y)&=\frac{ik^2}{2}\Big[\sum_{l=1}^{N} \frac{C^{1/2,\e}_{jl}(u)}{\sqrt{\beta_j-\beta_l}}y_l e^{i(\beta_l -\beta_j)u}+ \int_{\xi}^{k^2}\frac{C^{1/2,\e}_{j\ga'}(u)}{\sqrt{\beta_j \sgap}}y_{\ga'}e^{i(\sgap -\beta_j)u}d\ga' \Big],\\
 \tilde{\textbf{H}}^{aa}_{\ga,\e,\xi}(u)(y)&=\frac{ik^2}{2}\sum_{l=1}^{N} \frac{C^{1/2,\e}_{\ga l}(u)}{\sqrt{\sga \beta_l}}y_l e^{i(\beta_l -\sga)u},\\
\tilde{\textbf{G}}^{aa}_{j,\e,\xi}(u)(y)=&\,\frac{ik^4}{4}\sum_{l=1}^{N}\int_{-1/\xi}^{-\xi} \frac{C^{1/2,\e}_{j\ga'}(u)C^{1/2,\e}_{\ga' l}(\tilde{u}+u)}{\sqrt{\beta_j\lvert \ga'\rvert \beta_l}}e^{i\beta_l (\tilde{u}+u) -\sqrt{\lvert \ga'\rvert}\lvert \tilde{u}\rvert -i\beta_j u}d\ga'd\tilde{u}\\
&+\frac{ik^2}{2}\Big[\sum_{l=1}^{N} \frac{C^{1,\e}_{jl}(u)}{\sqrt{\beta_j \beta_l}}y_l e^{i(\beta_l -\beta_j)u}+ \int_{\xi}^{\ko{}}\frac{C^{1,\e}_{j\ga'}(u)}{\sqrt{\beta_j \sgap}}y_{\ga'}e^{i(\sgap -\beta_j)u}d\ga' \Big],
\end{split}\]
where $C^{1/2,\e}(u)=0$ if $u\not\in[0,L/\e]$ and $\tilde{\textbf{G}}^{aa}_{\ga,\e,\xi}=0$ for $\ga\in(\xi,\ko{})$.
The following lemma is a direct consequence of Lemma \ref{approxcoef}.
\begin{lem}\label{approxop}
We have
\[
 \frac{1}{\e}\emph{\textbf{H}}^{aa}_{\e,\xi}\Big(\frac{u}{\e}\Big)(y)=\frac{1}{\se}\tilde{\emph{\textbf{H}}}^{aa}_{\e,\xi}\Big(\frac{u}{\e}\Big)(y)+o(1)\quad \text{and}\quad\frac{1}{\e}\emph{\textbf{G}}^{aa}_{\e,\xi}\Big(\frac{u}{\e}\Big)(y)=\tilde{\emph{\textbf{G}}}^{aa}_{\e,\xi}\Big(\frac{u}{\e}\Big)(y)+o(1),
\]
uniformly in $u\in[0,L/\e]$ and for each $y\in \esp_\xi$.
\end{lem}
This second lemma will be used in the next section to approximate the forward mode amplitudes $\widehat{a}^{\e,\xi}$ unique solution of \eqref{eqfora}.

\begin{preuve}[of Lemma \ref{approxcoef}]
The proof consists only in doing expansions of the perturbations produced by the waveguide transformation.

\[
C^\e_{rs}(u)=\int_0^{+\infty} (J_\e(u,v)n^2_\e(u,v)-n^2(v))  \phi_r(v)\phi_s(v)dv=D^{1,\e}_{rs}(u)+D^{2,\e}_{rs}(u)+D^{3,\e}_{rs}(u),
\]
where $J_\e$ is defined by \eqref{jacodef}, and
\[\begin{split}
D^{1,\e}_{rs}(u)&=\int_0^{+\infty}  \phi_r(v)\phi_s(v)(J_\e(u,v)-1)n^2(v)dv,\\
D^{2,\e}_{rs}(u)&=\int_0^{+\infty}  \phi_r(v)\phi_s(v)(n^2_\e(u,v)-n^2(v))dv,\\
D^{3,\e}_{rs}(u)&=\int_0^{+\infty}  \phi_r(v)\phi_s(v)(J_\e(u,v)-1)(n^2_\e(u,v)-n^2(v))dv.\\\\
\end{split}\]
First, we have
\begin{equation}\label{defD1}\begin{split}
D^{1,\e}_{rs}(u)=&\, \se \int_0^{+\infty}  \phi_r(v)\phi_s(v) n^2(v)\Big[\frac{4v^2}{\pi} \int \frac{V_s(\tilde{u})f_s(\tilde{u})}{((\tilde{u}-u)^2+v^2)^2}d\tilde{u} -\frac{2}{\pi}\int \frac{V_s(\tilde{u})f_s(\tilde{u})}{(\tilde{u}-u)^2+v^2}d\tilde{u}\Big]\\
&+\e \int_0^{+\infty}  \phi_r(v)\phi_s(v) n^2(v)\Big[\frac{4v^2}{\pi^2}\Big(\int \frac{(\tilde{u}-u)V_s(\tilde{u})f_s(\tilde{u})}{((\tilde{u}-u)^2+v^2)^2}d\tilde{u}\Big)^2\\
&-\frac{2}{\pi}\int \frac{\mathcal{U}(V_s(\cdot)f_s(\e \cdot))(\tilde{u})\frac{d}{dz}V_s(\tilde{u})f_s(\e \tilde{u})}{(\tilde{u}-u)^2+v^2}d\tilde{u}+\frac{4v^2}{\pi}\int \frac{\mathcal{U}(V_s(\cdot)f_s(\e \cdot))(\tilde{u})\frac{d}{dz}V_s(\tilde{u})f_s(\e \tilde{u})}{((\tilde{u}-u)^2+v^2)^2}d\tilde{u}\\
&-\frac{4v^2}{\pi^2}\int \frac{V_s(\tilde{u})f_s(\e \tilde{u})}{(\tilde{u}-u)^2+v^2}d\tilde{u}\int \frac{V_s(\tilde{u})f_s(\e \tilde{u})}{((\tilde{u}-u)^2+v^2)^2}d\tilde{u}\\
&+\frac{1}{\pi^2}\Big(\int \frac{V_s(\tilde{u})f_s(\e \tilde{u})}{(\tilde{u}-u)^2+v^2}d\tilde{u}\Big)^2+\frac{4v^4}{\pi^2}\Big(\int \frac{V_s(\tilde{u})f_s(\e \tilde{u})}{((\tilde{u}-u)^2+v^2)^2}d\tilde{u}\Big)^2\Big]\\
&\hspace{0.5cm}+o(\e).
\end{split}\end{equation}
Second, we have
\[
D^{2,\e}_{rs}(u)= \int_{A_{1,\e}}  (1-n_1^2)\phi_r(v)\phi_s(v)dv+\int_{A_{2,\e}}(n_1^2-1)  \phi_r(v)\phi_s(v)dv,
\]
where
\[
A_{1,\e}=\big\{ v<d,\quad x(u,v)>d\big\}\quad\text{and}\quad A_{2,\e}=\big\{ v>d,\quad x(u,v)<d\big\}.
\]
The real part $x(u,\cdot)$ of the conformal map \eqref{conformmap} is a nondecreasing bijection since
\[ \sup_{\substack{u\geq0\\ v\in[\eta,M]}}\lvert \partial_v x(u,v)-1\rvert\leq \se K. \]
Moreover, we have
\[
\lvert x^{-1}(u,d)-d -(d-x(u,d))\partial_v x^{-1}(u,x(u,d))\rvert \leq \frac{1}{2}\lvert d-x(u,d)\rvert^2 \sup_{\tilde{v}\in[d,x(u,d)]}\lvert \partial^2_v x^{-1}(u,\tilde{v})\rvert,
\]
with
\[
\sup_{u\geq 0} \sup_{\tilde{v}\in[d,x(u,d)]}\lvert \partial^2_v x^{-1}(u,\tilde{v})\rvert\leq K\se,
\quad\text{and}\quad
\sup_{u\geq 0} \lvert d-x(u,d)\rvert \leq K\se.
\]
Then, $D^{2,\e}_{rs}$ can be approximated as follows
\[
D^{2,\e}_{rs}(u)= \int_{\tilde{A}_{1,\e}}  (1-n_1^2)\phi_r(v)\phi_s(v)dv+\int_{\tilde{A}_{2,\e}}(n_1^2-1)  \phi_r(v)\phi_s(v)dv+o(\e),
\]
where
\[
\tilde{A}_{1,\e}=\big\{v<d,\quad v>d+Y^\e(u,d) \big\}\qquad\text{and}\qquad\tilde{A}_{2,\e}=\big\{ v>d,\quad v<d+Y^\e(u,d)\big\},
\]
with
\[
Y^\e(u,d)= \frac{d}{\pi}\int \frac{x(u,0)}{(\tilde{u}-u)^2+d^2}d\tilde{u}\Big[1-\frac{1}{\pi}\int \frac{x(u,0)}{(\tilde{u}-u)^2+d^2}d\tilde{u}+2\frac{d^2}{\pi}\int \frac{x(u,0)}{((\tilde{u}-u)^2+d^2)^2}d\tilde{u}\Big].
\]
Consequently, we have
\[\begin{split}
D^{2,\e}_{rs}(u)= &\,\se\, 2(n_1^2-1)\phi_r(d)\phi_s(d)\frac{d}{\pi}\int \frac{V_s(\tilde{u})f_s(\e \tilde{u})}{(\tilde{u}-u)^2+d^2}d\tilde{u}\\
&+\e\,2(n_1^2-1)\phi_r(d)\phi_s(d)\Big[2\frac{d^3}{\pi^2}\int \frac{V_s(\tilde{u})f_s(\e \tilde{u})}{(\tilde{u}-u)^2+d^2}d\tilde{u}\int \frac{V_s(\tilde{u})f_s(\e \tilde{u})}{((\tilde{u}-u)^2+d^2)^2}d\tilde{u} \\
&\hspace{1cm}-\frac{d}{\pi^2}\Big(\int \frac{V_s(\tilde{u})f_s(\e \tilde{u})}{(\tilde{u}-u)^2+d^2}d\tilde{u}\Big)^2+\frac{d}{\pi}\int \frac{\mathcal{U}(V_s(\cdot)f_s(\e\cdot))(\tilde{u})\frac{d}{dz}V_s(\tilde{u})f_s(\e \tilde{u})}{(\tilde{u}-u)^2+d^2}d\tilde{u}\Big]\\
&+\e4(n_1^2-1)(\phi_r(d)\phi'_s(d)+\phi'_r(d)\phi_s(d))\Big(\frac{d}{\pi}\int \frac{V_s(\tilde{u})f_s(\e \tilde{u})}{(\tilde{u}-u)^2+d^2} d\tilde{u}\Big)^2+o(\e).
\end{split}\]
Finally, according to the previous step, we have
\[\begin{split}
D^{3,\e}_{rs}(u)=\e\phi_r(d)\phi_s(d)\Big[\frac{4d^2}{\pi} \int \frac{V_s(\tilde{u})f_s(\e \tilde{u})}{((\tilde{u}-u)^2+d^2)^2}d\tilde{u} -&\frac{2}{\pi}\int \frac{V_s(\tilde{u})f_s(\e \tilde{u})}{(\tilde{u}-u)^2+d^2}d\tilde{u}\Big]\\
&\times \frac{d}{\pi}\int \frac{V_s(\tilde{u})f_s(\e \tilde{u})}{(\tilde{u}-u)^2+d^2}d\tilde{u}+o(\e).
\end{split}\]
To finish the proof of this lemma, we have to approximate in $D^{1,\e}(u)$ and $D^{2,\e}(u)$ the terms of order $\se$ to obtain a more convenient form to exploit the mixing properties of the random perturbations. Let $\alpha_0\in(1/2,1)$. First, we have
\[
\int \frac{V_s(\tilde{u})f_s(\e \tilde{u})}{(\tilde{u}-u)^2+d^2}d\tilde{u}=\int \frac{V_s(u+w)f_s(\e(u+w))}{w^2+d^2}dw,
\]
and then
\[
\Big\lvert \int_{\substack{\lvert w\rvert >1/\e^{\alpha_0} \\ u+w\geq 0}} \frac{V_s(u+w)f_s(\e(u+w))}{w^2+d^2}dw\Big\rvert  \leq
K\int_{1/(d\e^{\alpha_0})}^{+\infty}\frac{1}{w^2+1}dw\leq K\e^{\alpha_0}.
\]
In the same way, we have
\[
\int_0^{+\infty}\phi_r(v)\phi_s(v)\int \frac{V_s(\tilde{u})f_s(\e \tilde{u})}{(\tilde{u}-u)^2+v^2}d\tilde{u}dv=\int_0^{+\infty}\phi_r(v)\phi_s(v)\int \frac{V_s(u+w)f_s(\e(\tilde{u}+w))}{w^2+v^2}dw dv,
\]
so that
\[\begin{split}
\Big\lvert \int_0^{+\infty}\phi_r(v)\phi_s(v) &\int_{\substack{\lvert w\rvert >1/\e^{\alpha_0} \\ u+w\geq 0}} \frac{V_s(u+w)f_s(\e(u+w))}{w^2+v^2}dwdv\Big\rvert \\ 
&\leq K \int_0^{+\infty}\frac{\lvert \phi_r(v)\phi_s(v)\rvert }{v}\int_{1/(v\e^{\alpha_0})}^{+\infty}\frac{1}{w^2+1}dwdv\\
&\leq \int_0^{+\infty}\lvert \phi_r(v)\phi_s(v)\rvert dv K\e^{\alpha_0}.
\end{split}\]
 Let us note that $ \int_0^{+\infty}\lvert \phi_r(v)\phi_s(v)\rvert dv$ is well defined thanks to Assumption \ref{hyp3} and the form of $\phi_j$ given in Section \ref{spectralP2}. The last term can be treated exactly as the previous one. $\square$
\end{preuve}

\subsubsection{Approximation of the Transfer Operator}

The forward mode amplitudes $\widehat{a}^{\e,\xi}$ unique solution of \eqref{eqfora} cannot be studied directly using the perturbed test function method introduced below because the random operator $\textbf{H}^{aa}_{\e,\xi}$ and $\textbf{G}^{aa}_{\e,\xi}$ are not adapted to the filtration of the random perturbations. Consequently, we approximate $\widehat{a}^{\e,\xi}$ by a new process from which we can exploit the mixing properties of the random perturbations. Let us introduce $\tilde{a}^{\e,\xi}$ the unique solution of the differential equation
\begin{equation}\label{approxtransf}
\frac{d}{du} \tilde{a}^{\e,\xi}(u)=\frac{1}{\sqrt{\e}} \tilde{\textbf{H}}^{aa} _{\e,\xi} \left(\frac{u}{\e}\right)
 \tilde{a}^{\e,\xi} (u)+ \big<\tilde{\textbf{G}}^{aa}\big>(u) \tilde{a}^{\e,\xi} (u),
\end{equation}
with $ \tilde{a}^{\e,\xi} (0)=  \widehat{a}^\xi_0$ defined by \eqref{condinith}, and where $\big<\textbf{G}^{aa}\big>$ is defined for all $ y \in \mathcal{H}_\xi$ and $j\in \big\{1,\dots,N\big\}$ by
\[
 \big<\textbf{G}^{aa}\big> _{j}(u,y)=\int_{-\infty}^{-\xi}\frac{i k^{4}f^2_s(u)}{2\beta_j \sqrt{\lvert\ga\rvert}} \int_{0}^{+\infty}G^{(1)}_{j\ga}(z)\cos\big(\beta_j z\big)e^{-\sqrt{\lvert\ga\rvert}z}dzd\ga y_j 
+\frac{ik^2f^2_s(u)}{2\beta_j}G^{(2)}_j y_j,  
\]
 and $\big<\textbf{G}^{aa}\big> _{\ga}(y)=0$ for almost every $\ga\in(\xi,k^2)$. Here, we have
 \[\begin{split}
G^{(1)}_{j\ga}(z)=& 4(n_1-1)^2\phi^2_j(d)\phi^2_\ga(d)\frac{d^2}{\pi^2}\iint \frac{R_s(z+w-w')}{(w^2+d^2)({w'}^2+d^2)}dwdw'\\
&+\frac{4}{\pi^2}\iint\phi_j (v)\phi_{\ga} (v)\phi_j (v')\phi_\ga (v')\iint\frac{R_s(z+w-w')}{(w^2+v^2)({w'}^2+{v'}^2)}dwdw'dvdv'\\
&+ \frac{16}{\pi^2}\iint\phi_j (v)\phi_{\ga} (v)\phi_j (v')\phi_\ga (v')v^2{v'}^2\iint\frac{R_s(z+w-w')}{(w^2+v^2)^2({w'}^2+{v'}^2)^2}dwdw'dvdv'\\
&-\frac{8(n_1-1)}{\pi}\phi_j(d)\phi_\ga(d)\int\phi_j (v)\phi_{\ga} (v)\int\frac{R_s(z+w-w')}{(w^2+v^2)({w'}^2+{d}^2)}dwdw'dv\\
&+\frac{16(n_1-1)}{\pi}\phi_j(d)\phi_\ga(d)\int\phi_j (v)\phi_{\ga} (v)v^2\int\frac{R_s(z+w-w')}{(w^2+v^2)^2({w'}^2+{d}^2)}dwdw'dv\\
&- \frac{8}{\pi^2}\iint\phi_j (v)\phi_{\ga} (v)\phi_j (v')\phi_\ga (v')v^2\iint\frac{R_s(z+w-w')}{(w^2+v^2)^2({w'}^2+{v'}^2)}dwdw'dvdv'\\
&- \frac{8}{\pi^2}\iint\phi_j (v)\phi_{\ga} (v)\phi_j (v')\phi_\ga (v'){v'}^2\iint\frac{R_s(z+w-w')}{(w^2+v^2)({w'}^2+{v'}^2)^2}dwdw'dvdv',
 \end{split}\]
 and
 \[\begin{split}
 G&^{(2)}_j=2(n_1^2-1)\phi_j^2(d)\\
 &\times\Big[\frac{6d^3}{\pi^2}\iint \frac{R_s(w-w')}{(w^2+d^2)({w'}^2+d^2)^2}dwdw' -\frac{3d}{\pi^2}\iint \frac{R_s(w-w')}{(w^2+d^2)({w'}^2+d^2)}dwdw'+\mathcal{U}(\tilde{R}_s)(0)\Big]\\
 &+4(n_1^2-1)\phi_j(d)\phi'_j(d)\iint \frac{R_s(w-w')}{(w^2+d^2)({w'}^2+d^2)}dwdw'\\
&+\int_0^{+\infty} \phi_j^2(v)n^2(v)\Big[\frac{4v^2}{\pi^2}\iint \frac{ww' R_s(w-w')}{({w'}^2+v^2)^2(w^2+v^2)^2}dwdw'\\
&-\frac{4v^2}{\pi^2}\iint \frac{R_s(w-w')}{(w^2+v^2)({w'}^2+v^2)^2}dwdw'\\
&+\frac{1}{\pi^2}\iint \frac{R_s(w-w')}{(w^2+v^2)({w'}^2+v^2)}dwdw'+\frac{4v^4}{\pi^2}\iint \frac{R_s(w-w')}{(w^2+v^2)^2({w'}^2+v^2)^2}dwdw'\Big],
 \end{split}\]
 where $R_s$ is the autocorrelation function of the random process $V_s$, and $\tilde{R}_s$ is the correlation function of the processes $V_s$ and $\frac{d}{dz}V_s$, that is $\frac{d}{dz}R_s$. Let us note $\mathcal{U}(\tilde{R}_s)(0)$ is well defined. In fact, $R_s$ is an even function so that $\tilde{R}_s(0)=0$.
 
 Moreover, let us remark that the terms involving the Hilbert transform in \eqref{defD1} have disappeared. In fact, taking the expectation, using the stationarity of $V_s$ by making a change of variable, and then passing to the limit in $\e$, we have
 \[\begin{split}
 \lim_{\e\to0} \int d\tilde{u}& f_s(\e\tilde{u})\frac{1}{\pi}\int \frac{\tilde{R}_s(w)}{w} f_s(\e(\tilde{u}-w)) \Big[\frac{4v^2/\pi}{\big((\tilde{u}-u/\e)^2+v^2\big)^2}-\frac{2/\pi}{(\tilde{u}-u/\e)^2+v^2}\Big]\\
 &=f^2_s(u)\mathcal{U}(\tilde{R}_s)(0)\Big[\frac{4v^2}{\pi}\int \frac{d\tilde{u}}{\big(\tilde{u}^2+v^2\big)^2}-\frac{2}{\pi}\int \frac{d\tilde{u}}{\tilde{u}^2+v^2}\Big]=f^2_s(u)\mathcal{U}(\tilde{R}_s)(0)\Big[\frac{2}{v}-\frac{2}{v}\Big]=0.
\end{split} \]

 We have the following proposition that describes the relation between the two transfer processes $\widehat{a}^{\e,\xi}$ and $\tilde{a}^{\e,\xi}$.
\begin{prop}\label{propoP2}
\[\forall \eta>0,\quad \lim_{\e \to 0}\mathbb{P}\left( \sup_{u\in[0,L]}\|\widehat{a}^{\e,\xi}(u)-\tilde{a}^{\e,\xi}(u)\|^2_{\esp_{\xi}} > \eta \right)=0.\]
\end{prop}  
Thanks to Proposition \ref{propoP2} and \cite[Theorem 3.1]{biling}, one can study the new process $(\tilde{a}^{\e,\xi})_{\e}$ instead of $(\widehat{a}^{\e,\xi})_{\e}$.
Let us remark that $\tilde{a}^{\e,\xi}$ is adapted to the filtration $\tilde{\mathcal{F}}^\e _{0,u}$ defined by \eqref{filtration2} and $\|\tilde{a}^{\e,\xi} (u)\|^2_{\mathcal{H}_{\xi}}=\|\widehat{a}^\xi_0\|^2_{\mathcal{H}_{\xi}}$, for all $ u\geq 0$. Consequently, $(\tilde{a}^{\e,\xi})_{\e}$ can be studied using the perturbed test function method. 
\begin{preuve}[of Proposition \ref{propoP2}]
This proof is in two steps. The first step consists in using Lemma \ref{approxop} to approximate the random operator $\textbf{H}^{aa}_{\e,\xi}$ and $\textbf{G}^{aa}_{\e,\xi}$ by $\tilde{\textbf{H}}^{aa}_{\e,\xi}$ and $\tilde{\textbf{G}}^{aa}_{\e,\xi}$ from which we can exhibit the mixing properties of the random perturbations. The second step consists in using the ideas developed in \cite{khasminskii} to exploit the fast phase of the random operator $\textbf{G}^{aa}_{\e,\xi}$. The proof of the second step follows closely the proof of Proposition 6.3 in \cite{gomez2bis}.
$\blacksquare$
\end{preuve}

\subsubsection{The Perturbed Test Function Method}

In this section we study the limit in distribution of the process $(\tilde{a}^{\e,\xi})_{\e}$ unique solution of the differential equation \eqref{approxtransf}, using the perturbed-test-function method and a martingale technique. The proof is in two steps. First, we prove the tightness of the process $(\tilde{a}^{\e,\xi})_{\e}$, afterward we identify all the subsequence limits thanks to a well posed martingale problem. According to \cite[Theorem 4]{kushner}, the proof of the tightness is decomposed in Lemma \ref{f1P2}, Lemma \ref{f2P2} , and Lemma \ref{A1P2}. The subsequence limits are identify, first by a complex martingale problem (Lemma \ref{propmgP2}), and then by a real martingale problem (Lemma \ref{uniquemgP2}) to use the classical uniqueness results in \cite{inf} for instance.

Now, let us defined the good space to apply the asymptotic analysis. Let $r=\|\widehat{a}^\xi_0\|_{\mathcal{H}_{\xi}}$,
\[\mathcal{B}_{r, \mathcal{H}_{\xi}}=\left\{\lambda \in \mathcal{H}_\xi, \|\lambda\|_{\mathcal{H}_\xi}=\sqrt{\left<\lambda,\lambda\right>_{\mathcal{H}_\xi}} \leq r_y\right\}\] 
the closed ball with radius $r$, and $\{g_n, n\geq 1\}$ a dense subset of $\mathcal{B}_{r,\mathcal{H}_\xi}$. We equip  $\mathcal{B}_{r, \mathcal{H}_\xi}$ with the distance $d_{\mathcal{B}_{r, \mathcal{H}_\xi}}$ defined by
\[d_{\mathcal{B}_{r,\mathcal{H}_\xi}}(\lambda, \mu)=\sum_{j=1}^{+\infty}\frac{1}{2^j}\left\lvert\big<\lambda-\mu,g_n\big>_{\mathcal{H}_\xi}\right\rvert\]
for all $ (\lambda,\mu)\in\mathcal{B}_{r,\mathcal{H}_\xi}$. As a result, $(\mathcal{B}_{\mathcal{H}_\xi} ,d_{\mathcal{B}_{r, \mathcal{H}_\xi}})$ is a compact metric space.
From the definition of the metric $d_{\mathcal{B}_{r,\mathcal{H}_\xi}}$, we have the following criterion.
\begin{thm}\label{crit}
A family of processes $(X^\e)_{\e\in(0,1)}$ is tight in $\mathcal{C}([0,+\infty),(\mathcal{B}_{r,\mathcal{H}_\xi},d_{\mathcal{B}_{r, \mathcal{H}_\xi}} ))$ if and only if $\big(\big<X^\e,\lambda\big>_{\mathcal{H}_\xi}\big)_{\e\in(0,1)}$ is tight on $\mathcal{C}([0,+\infty),\mathbb{C})$ for all $\lambda \in \mathcal{H}_\xi$.
\end{thm}
This last theorem looks like the tightness criterion of Mitoma and Fouque \cite{mitoma,fouque}.
For any $\lambda \in \mathcal{H}_\xi$, we set $\tilde{a}^{\e,\xi}_\lambda (u)=\big<\tilde{a}^{\e,\xi}(u),\lambda\big>_{\mathcal{H}_\xi}$. According to Theorem \ref{crit}, the family $(\tilde{a}^{\e,\xi}(.))_{\e}$ is tight on $\mathcal{C}([0,+\infty), (\mathcal{B}_{r,\mathcal{H}_\xi},d_{\mathcal{B}_{r,\mathcal{H}_\xi}} ))$ if and only if the family $(\tilde{a}^{\e,\xi}_\lambda(.))_{\e }$ is tight on $\mathcal{C}([0,+\infty),\mathbb{C} )$  for all $\lambda \in \mathcal{H}_\xi$. Furthermore, $(\tilde{a}^{\e,\xi}(.))_{\e}$ is a family of continuous processes and then it is sufficient to prove that for all $\lambda \in \mathcal{H}_\xi$, $(\tilde{a}^{\e,\xi}_\lambda(.))_{\e}$ is tight in the space of cad-lag functions $\mathcal{D}([0,+\infty),\mathbb{C} )$ equipped with the Skorokhod topology \cite[Theorem 13.4]{biling}. 

Using the notion of a pseudogenerator, we prove tightness and characterize all subsequence limits. Let us recall the techniques developed by Kurtz and Kushner. Let $\mathcal{M}^\e $ be the set of all $\mathcal{F}^\e$-measurable functions $f(u)$ for which $\sup_{u\leq L} \mathbb{E}\left[\lvert f(u) \rvert \right] <+\infty $ and where $L>0$ is fixed. The $p-\lim$ and the pseudogenerator are defined as follows. Let $f$ and $f^\delta$ in $\mathcal{M}^\e $ for all $\delta>0$. We say that $f=p-\lim_\delta f^\delta$ if
\[ \sup_{u, \delta }\mathbb{E}[\lvert f^\delta(u)\rvert]<+\infty\quad \text{and}\quad \lim_{\delta\rightarrow 0}\mathbb{E}[\lvert f^\delta (u) -f(u)\rvert]=0 \quad \forall u.\]
The domain of $\mathcal{A}^\e$ is denoted by $\mathcal{D}\left(\mathcal{A}^\e\right)$. We say that $f\in \mathcal{D}\left(\mathcal{A}^\e\right)$ and $\mathcal{A}^\e f=g$ if $f$ and $g$ are in $\mathcal{D}\left(\mathcal{A}^\e\right)$ and 
\[
p-\lim_{\delta \to 0} \left[ \frac{\mathbb{E}^\e _u [f(u+\delta)]-f(u)}{\delta}-g(u) \right]=0,
\]
where $\mathbb{E}^\e _u$ is the conditional expectation given $\mathcal{F}^\e _u$ and $\mathcal{F}^\e _u =\mathcal{F}_{u/\e}$. A useful result about $\mathcal{A}^\e$ is given by the following theorem.
\begin{thm}\label{martingaleP2}
Let $f\in \mathcal{D}\left(\mathcal{A}^\e\right)$. Then,
\begin{equation*}
M_f ^\e (u)=f(u)-\int _0 ^u  \mathcal{A}^\e f(v)dv
\end{equation*}
is an $\left( \mathcal{F}^\e _u \right)$-martingale.
\end{thm} 
In what follows, we consider the classical complex derivative with the following notation: If $v=\alpha +i\beta$, then $\partial_v =\frac{1}{2}\left(\partial_\alpha -i \partial_ \beta \right)$ and $\partial_{\overline{v}}=\frac{1}{2}\left(\partial_\alpha +i \partial_ \beta \right)$.
\begin{prop}\label{tightth1P2}
For all $ \lambda \in \mathcal{H}_\xi$, the family $\big(\tilde{a}^{\e,\xi} _{\lambda} (.) \big)_{\e\in(0,1)}$ is tight in $\mathcal{D}\left([0,+\infty),\mathbb{C} \right)$.
\end{prop}
\begin{preuve}
First, we easily obtain for all $T>0$,  \[\lim_{M\to+\infty }\varlimsup_{\e \to 0} \mathbb{P}\left( \sup_{0 \leq u \leq T} \lvert \tilde{a}^{\e,\xi}_\lambda (u) \rvert \geq M \right)=0,\]
since $(\tilde{a}^{\e,\xi}(.))_\e$ is a bounded process. To show the tightness of the process $\big(\tilde{a}^{\e,\xi}(u)\big)_{u \geq 0}$, according to \cite[Theorem 4]{kushner}, we need to show the three following lemmas. Let $\lambda \in \mathcal{H}_\xi$, $f$ be a smooth function, and $f_0 ^\e (u)=f\left(\tilde{a}^{\e,\xi} _\lambda (u) \right)$. We have,
\[
\begin{split}
\mathcal{A}^\e f_0 ^\e (u)&= \partial_v f\left(\tilde{a}^{\e,\xi} _{\lambda} (t) \right)\left[\frac{1}{\sqrt{\e}}H _{\lambda}\left(\tilde{a}^{\e,\xi}  (u),C^{1/2,\e}\left(\frac{u}{\e}\right),\frac{u}{\e} \right)+G_\lambda\left(\tilde{a}^{\e,\xi} (u)\right)\right] \\
&+ \partial_{\overline{v}} f\left(\tilde{a}^{\e,\xi} _{\lambda} (u)\right)\overline{ \left[\frac{1}{\sqrt{\e}}H _{\lambda}\left(\tilde{a}^{\e,\xi}  (u),C^{1/2,\e}\left(\frac{u}{\e}\right),\frac{u}{\e} \right)+G_\lambda\left(\tilde{a}^{\e,\xi} (u)\right)\right]},
\end{split}
\]
where
\[
H _{\lambda}\Big(\tilde{a}^{\e,\xi}  (u),C^{1/2,\e}\Big(\frac{u}{\e}\Big),\frac{u}{\e} \Big)=\big<\tilde{\textbf{H}}^{aa}_{\e,\xi}\Big(\frac{u}{\e}\Big)\tilde{a}^{\e,\xi}  (u)(y),\lambda\big>_{\esp_\xi},
\]
and $G_\lambda(\tilde{a}^{\e,\xi} (u))=\big<\big<\tilde{\textbf{G}}^{aa}\big>(u)\tilde{a}^{\e,\xi}  (u),\lambda\big>_{\esp_\xi}$. 
Let us consider
\[
\begin{split}
f^\e _1 (u)&=\frac{1}{\sqrt{\e}}\partial_v f\left(\tilde{a}^{\e,\xi} _{\lambda} (u) \right)\int_{u}^{+\infty} \mathbb{E}^\e _u\left[ H _{\lambda}\left(\tilde{a}^{\e,\xi}  (u),C^{1/2,\e}\left(\frac{w}{\e}\right),\frac{w}{\e} \right)\right] dw  \\
&+\frac{1}{\sqrt{\e}}\partial_{\overline{v}} f\left(\tilde{a}^{\e,\xi} _{\lambda} (u) \right)\int_{t}^{+\infty} \mathbb{E}^\e _u\left[ \overline{H_{\lambda}\left(\tilde{a}^{\e,\xi} (u),C^{1/2,\e}\left(\frac{w}{\e}\right),\frac{w}{\e} \right)}\right]dw. \\
\end{split}
\]
Here, $\E^\e_u$ stands for the conditional expectation with respect to the filtration $\tilde{\mathcal{F}}^\e_{0,u}$ defined by \eqref{filtration2}.
\begin{lem}\label{f1P2}  For all $T>0$, we have
$\lim_{\e} \sup_{0\leq u \leq T} \lvert f^\e _1 (u)\rvert =0$ almost surely, and $\sup_{u\geq 0}\mathbb{E}\left[\lvert f^\e _1 (u) \rvert \right]=\mathcal{O}(\e^{1/2-\eta})$ for all $\eta>0 $.
\end{lem}
\begin{preuve}[of Lemma \ref{f1P2}]
Using a change of variable and thanks to \eqref{phimixing},  we obtain
\[\begin{split}
\lvert f^{\e}_1(u)\rvert& \leq K \se \int_0^{+\infty} \phi(w)dw \int_0^{+\infty}\psi(w)(1+1/(2w))dw \\
&+K\se \ln(1+1/(d^2 \e^{2\alpha_0}))+K\se \int_0^{+\infty}\psi(w)(\ln(1+1/(w^2 \e^{2\alpha_0}))+1/(2w))(1+1/w)dw,
\end{split}\]
where, thanks to Assumption \ref{hyp3}, $\psi$ is a function coming from the propagating modes $\phi_j$ such that the integral in $w$ is well defined. More precisely, $\psi(w)\sim w^2$ as $w$ goes to $0$ and $\psi$ decay exponentially at infinity. The second and third part of the previous expression come from the fact that in \eqref{phimixing} $s_1$ has to be greater than $2/\e^{\alpha_0}$. Indexes $s_1$ smaller than $2/\e^{\alpha_0}$ are not covered by the mixing property and then lead us to explicit computation using a primitive of $\arctan(w)$ given by $w\arctan(w) -1/2\ln(1+w^2) $. That concludes the proof of Lemma \ref{f1P2}.$\square$
\end{preuve}
A computation gives us
\[\begin{split}
\mathcal{A}^\e& \big(f^\e _0 +f^\e _1\big)(u)=\partial_v f\big(\tilde{a}^{\e,\xi} _{\lambda} (u) \big)\Big[\int_{0}^{+\infty} H^{(1)} _{\lambda}\Big(\tilde{a}^{\e,\xi}  (u), \mathbb{E}^\e _u \Big[C^{1/2,\e}\Big(\frac{u}{\e}\Big) \otimes C^{1/2,\e}\Big(\frac{u}{\e}+w\Big)\Big] ,\frac{u}{\e},\frac{u}{\e}+w \Big)dw \\
&\hspace{3cm}+G_\lambda\Big(\tilde{a}^{\e,\xi} (u)\Big)\Big]\\
&+\partial_{\overline{v}} f\big(\tilde{a}^{\e,\xi} _{\lambda} (u) \big)\Big[\int_{0}^{+\infty} \overline{ H^{(1)} _{\lambda}\Big(\tilde{a}^{\e,\xi}  (u),\mathbb{E}^\e _u\Big[C^{1/2,\e}\Big(\frac{u}{\e}\Big) \otimes C^{1/2,\e}\Big(\frac{u}{\e}+w\Big)\Big],\frac{u}{\e},\frac{u}{\e}+w \Big)} dw \\
&\hspace{3cm}+\overline{G_\lambda\Big(\tilde{a}^{\e,\xi} (u)\Big)}\Big]
\end{split}\]
\[\begin{split}
\hspace{0.5cm}&+\partial^2_v f\big(\tilde{a}^{\e,\xi} _{\lambda} (u) \big)\int_{0}^{+\infty} H^{(2)} _{\lambda}\Big(\tilde{a}^{\e,\xi}  (u),\mathbb{E}^\e _u\Big[ C^{1/2,\e}\Big(\frac{u}{\e}\Big)\otimes C^{1/2,\e}\Big(\frac{u}{\e}+w\Big)\Big],\frac{u}{\e},\frac{u}{\e}+w \Big) dw\\
&+\partial^2_{\overline{v}} f\big(\tilde{a}^{\e,\xi} _{\lambda} (u) \big)\int_{0}^{+\infty} \overline{ H^{(2)} _{\lambda}\Big(\tilde{a}^{\e,\xi}  (u),\mathbb{E}^\e _u\Big[ C^{1/2,\e}\Big(\frac{u}{\e}\Big)\otimes C^{1/2,\e}\Big(\frac{u}{\e}+w\Big)\Big],\frac{u}{\e},\frac{u}{\e}+w \Big)}dw\\
&+\partial_{\overline{v}}\partial_v f\big(\tilde{a}^{\e,\xi} _{\lambda} (u) \big)\int_{0}^{+\infty} H^{(3)} _{\lambda}\Big(\tilde{a}^{\e,\xi}  (u),\mathbb{E}^\e _u\Big[ C^{1/2,\e}\Big(\frac{u}{\e}\Big)\otimes C^{1/2,\e}\Big(\frac{u}{\e}+w\Big)\Big],\frac{u}{\e},\frac{u}{\e}+w \Big) dw\\
&+\partial_v \partial_{\overline{v}} f\big(\tilde{a}^{\e,\xi} _{\lambda} (u) \big)\int_{0}^{+\infty} \overline{ H^{(3)} _{\lambda}\Big(\tilde{a}^{\e,\xi}  (u),\mathbb{E}^\e _u\Big[ C^{1/2,\e}\Big(\frac{u}{\e}\Big)\otimes C^{1/2,\e}\Big(\frac{u}{\e}+w\Big)\Big],\frac{u}{\e},\frac{u}{\e}+w \Big)}dw\\
&+o(1),
\end{split}
\]
where
\[\begin{split}
H^{(1)} _\lambda& (\textbf{T},C,s,\tilde{s})=-\frac{k^4}{4}\sum_{j=1}^N \overline{\lambda_{j}}\Big[\sum_{l ,l'=1}^{N}\frac{C_{jlll'}}{\sqrt{\beta _j \beta _l ^2 \beta _{l'}}}e^{i (\beta _{l}-\beta _{j})\tilde{s} + i (\beta _{l'}-\beta _{l})s}\textbf{T}_{l'} \\
&+\sum_{l=1}^N\int_{\xi}^{k^2} \frac{C_{jll\ga'}}{ \sqrt{\beta _j \beta _l^2  \sqrt{\ga'}}}e^{i (\beta_l-\beta _j)\tilde{s}+ i (\sqrt{\ga'} -\beta _l)s} \textbf{T}_{\ga' }+\frac{C_{j\ga'\ga' l}}{ \sqrt{\beta _j \ga' \beta _{l }}}e^{i (\sgap -\beta_{j})\tilde{s}+i (\beta _{l} -\sgap )s} \textbf{T}_{l} d\ga' \Big]\\
&-\frac{k^4}{4}\int_\xi^{k^2} \Big[\sum_{l ,l'=1}^{N}\frac{C_{\ga lll'}}{\sqrt{\sga \beta _l ^2 \beta _{l'}}}e^{i (\beta _{l}-\sga)\tilde{s}+i (\beta _{l'}-\beta_l)s} \textbf{T}_{l'}\\
&+\sum_{l=1}^N\int_{\xi}^{k^2} \frac{C_{\ga ll\ga''}}{ \sqrt{\sga \beta _l^2  \sqrt{\ga''}}}e^{i (\beta_l -\sga)\tilde{s}+i (\sqrt{\ga''} -\beta_l)s} \textbf{T}_{\ga'' }d\ga'' \Big]\overline{\lambda_{\ga}}d\ga,
\end{split}\]

\[\begin{split}
H^{(j)} _\lambda &(\textbf{T},C,s,\tilde{s})= (-1)^{j-1}\frac{k^4}{4}\sum_{j,j'=1}^N \Big[\sum_{l ,l'=1}^{N}\frac{C_{jlj'l'}}{\sqrt{\beta _j \beta _l \beta_{j'} \beta _{l'}}}e^{i (\beta _{l}-\beta _{j})\tilde{s}-i(-1)^{j-1}(\beta _{l'}-\beta _{j'})s} \textbf{T}_{l}\textbf{T}_{l'} \\
&+\sum_{l=1}^N\int_{\xi}^{k^2} \frac{C_{jlj'\ga'_2}}{ \sqrt{\beta _j \beta _l \beta_{j'}  \sqrt{\ga'_2}}}e^{i(\beta _{l}-\beta _{j})\tilde{s}-i(-1)^{j-1}(\sqrt{\ga'_2}-\beta _{j'})s} \textbf{T}_{l}\textbf{T}_{\ga'_2} d\ga'_2 \\
&+\int_\xi^{k^2}\sum_{l'=1}^N\frac{C_{j\ga'_1j' l'}}{ \sqrt{\beta _j \sqrt{\ga'_1}\beta_{j'} \beta _{l' }}}e^{i (\sqrt{\ga'_1}-\beta_j)\tilde{s}-i(-1)^{j-1}(\beta _{l'} -\beta_{j'})s} \textbf{T}_{\ga'_1}\textbf{T}_{l'}d\ga'_1 \\
&+\int_\xi^{k^2}\int_\xi^{k^2}\frac{C_{j\ga'_1j'\ga'_2}}{ \sqrt{\beta _j \sqrt{\ga'_1} \beta_{j'}\sqrt{\ga'_2}}}e^{i (\sqrt{\ga'_1} -\beta _{j})\tilde{s}-i(-1)^{j-1}(\sqrt{\ga'_2}-\beta_{j'})s} \textbf{T}_{\ga'_1}\textbf{T}_{\ga'_2}  d\ga'_1 d\ga'_2\Big]\overline{\lambda_{j}\lambda_{j'}}
\end{split}\]
\[\begin{split}
\hspace{1cm}& -\frac{k^4}{4}\sum_{j=1}^N \int_\xi ^{k^2}\Big[\sum_{l ,l'=1}^{N}\frac{C_{jl \ga_2 l'}}{\sqrt{\beta _j \beta _l \sqrt{\ga_2} \beta _{l'}}}e^{i (\beta _{l}-\beta _{j})\tilde{s}-i(-1)^{j-1}(\beta _{l'}-\sqrt{\ga_2})s} \textbf{T}_{l}\textbf{T}_{l'}\\
&+\int_\xi^{k^2}\sum_{l'=1}^N\frac{C_{j\ga'_1\ga_2 l'}}{ \sqrt{\beta _j \sqrt{\ga'_1\ga_2} \beta _{l' }}}e^{i (\sqrt{\ga'_1}-\beta_j)\tilde{s}-i(-1)^{j-1}(\beta _{l'} -\sqrt{\ga_2})s} \textbf{T}_{\ga'_1}\textbf{T}_{l'} d\ga'_1 \Big]\overline{\lambda_{j}\lambda_{\ga_2}}d\ga_2\\
&-\frac{k^4}{4}\int_\xi ^{k^2}\sum_{j'=1}^N \Big[\sum_{l ,l'=1}^{N}\frac{C_{\ga_1 lj'l'}}{\sqrt{\sqrt{\ga_1} \beta _l \beta_{j'} \beta _{l'}}}e^{i (\beta _{l}-\sqrt{\ga_1})\tilde{s}-i(-1)^{j-1}(\beta _{l'}-\beta _{j'})s} \textbf{T}_{l}\textbf{T}_{l'} \\
&+\sum_{l=1}^N\int_{\xi}^{k^2} \frac{C_{\ga_1 lj'\ga'_2}}{ \sqrt{\sqrt{\ga_1} \beta _l \beta_{j'}  \sqrt{\ga'_2}}}e^{i(\beta _{l}-\sqrt{\ga_1})\tilde{s}-i(-1)^{j-1}(\sqrt{\ga'_2}-\beta _{j'})s} \textbf{T}_{l}\textbf{T}_{\ga'_2} d\ga'_2 \Big]\overline{\lambda_{\ga_1}\lambda_{j'}}\\
& -\frac{k^4}{4}\int_\xi^{k^2} \int_\xi ^{k^2}\Big[\sum_{l ,l'=1}^{N}\frac{C_{\ga_1l \ga_2 l'}}{\sqrt{\sqrt{\ga_1} \beta _l \sqrt{\ga_2} \beta _{l'}}}e^{i (\beta _{l}-\sqrt{\ga_1})\tilde{s}-i(-1)^{j-1}(\beta _{l'}-\sqrt{\ga_2})s} \textbf{T}_{l}\textbf{T}_{l'} \Big]\overline{\lambda_{\ga_1}\lambda_{\ga_2}}d\ga_1d\ga_2,
\end{split}\]
where $j\in\{2,3\}$. Let us note that a quick computation shows that $\mathcal{A}^\e(f^\e _0 +f^\e _1)(u)$ is not necessarily uniformly integrable. In fact,  we obtain for $\mathbb{E}\big[ \lvert \mathcal{A}^\e(f^\e _0 +f^\e _1)(u) \rvert ^2 \big]$  the same bound as the one obtained for  $f^{\e}_1(u)$ in Lemma \ref{f1P2} but without $\se$, so that this bound blows up as $\e$ goes to $0$. The problem comes from the indexes not covered by the mixing property \eqref{phimixing}. To show the tightness of the process  $\tilde{a}^{\e,\xi}$ we have to introduce a second perturbed test function to correct the problem. This second perturbed test function $f_2^\e$ is given by
\[\begin{split}
f_2^\e(u)=&\int_u^{+\infty} \tilde{F}_\lambda\Big(\tilde{a}^{\e,\xi}  (u),\E_u^\e\Big[\textbf{C}^\e\Big(\frac{w}{\e}\Big)\Big]- \E\Big[\textbf{C}^\e\Big(\frac{w}{\e}\Big)\Big],\frac{w}{\e}\Big)dw\end{split},\]
where
\[\begin{split}
\tilde{F}_\lambda(T,C,w)=&\,\partial_v f(T_{\lambda}  \big)\int_{0}^{+\infty}H^{(1)} _{\lambda}(T,C(w) \otimes C(w+\tilde{w}),w,w+\tilde{w} )d\tilde{w} \\
&+\partial_{\overline{v}} f(T _{\lambda} )\int_{0}^{+\infty} \overline{ H^{(1)} _{\lambda}(T,C(w) \otimes C(w+\tilde{w}),w,w+\tilde{w} )} d\tilde{w} \\
&+\partial^2_v f(T _{\lambda} )\int_{0}^{+\infty} H^{(2)} _{\lambda}(T,C(w) \otimes C(w+\tilde{w}),w,w+\tilde{w})d\tilde{w}\\
&+\partial^2_{\overline{v}} f(T _{\lambda} )\int_{0}^{+\infty} \overline{ H^{(2)} _{\lambda}(T,C(w) \otimes C(w+\tilde{w}),w,w+\tilde{w} )}d\tilde{w}\\
&+\partial_{\overline{v}}\partial_v f(T _{\lambda} )\int_{0}^{+\infty} H^{(3)} _{\lambda}(T,C(w) \otimes C(w+\tilde{w}),w,w+\tilde{w} )d\tilde{w}\\
&+\partial_v \partial_{\overline{v}} f(T_{\lambda})\int_{0}^{+\infty} \overline{ H^{(3)} _{\lambda}(T,C(w) \otimes C(w+\tilde{w}),w,w+\tilde{w} )}d\tilde{w}.
\end{split}\]

\begin{lem}\label{f2P2}  For all $ T>0$, we have
$\lim_{\e} \sup_{0\leq u \leq T} \lvert f^\e _2 (u)\rvert =0$ almost surely, and $\sup_{u\geq 0}\mathbb{E}[\lvert f^\e _2 (u) \rvert ]=\mathcal{O}(\e^{1-\eta})$, for all $\eta >0$.
\end{lem}
\begin{preuve}[of Lemma \ref{f2P2}]
Using a change of variable and thanks to \eqref{phimixing2},  we obtain
\[\begin{split}
\lvert f^{\e}_2(u)\rvert& \leq K \e\Big( \int_0^{+\infty} \phi^{1/2}(w)dw \int_0^{+\infty}\psi(w)(1+1/(2w))dw\Big)^2 \\
&+K\e\Big( \ln(1+1/(d^2 \e^{2\alpha_0}))\Big)^2+K\e \Big(\int_0^{+\infty}\psi(w)(\ln(1+1/(w^2 \e^{2\alpha_0}))+1/(2w))(1+1/w)dw\Big)^2,
\end{split}\]
where $\psi$ has been described in the proof of Lemma \ref{f1P2}. We recall that indexes $s_1$ smaller than $2/\e^{\alpha_0}$ in \eqref{phimixing} are not covered by the mixing property. These indexes lead us to explicit computations using a primitive of $\arctan(w)$ given by $w\arctan(w) -1/2\ln(1+w^2) $. That concludes the proof of Lemma \ref{f2P2}.$\square$
\end{preuve}

Finally, the tightness of $\tilde{a}^{\e,\xi}$ is given by the following lemma.

\begin{lem}\label{A1P2}
$\left\{\mathcal{A}^\e \big(f^\e _0 +f^\e _1+ f^\e _2\big)(u), \e \in(0,1), 0\leq u\leq T\right\}$ is uniformly integrable for all $T>0$, where
\[\begin{split}
\mathcal{A}^\e& \big(f^\e _0 +f^\e _1+f^\e_2 \big)(u)=\partial_v f\big(\tilde{a}^{\e,\xi} _{\lambda} (u) \big)\Big[\int_{0}^{+\infty} H^{(1)} _{\lambda}\Big(\tilde{a}^{\e,\xi}  (u), \mathbb{E}[C^{1/2,\e}(0) \otimes C^{1/2,\e}(w)] ,\frac{u}{\e},\frac{u}{\e}+w \Big)dw \\
&\hspace{3cm}+G_\lambda\Big(\tilde{a}^{\e,\xi} (u)\Big)\Big]\\
&+\partial_{\overline{v}} f\big(\tilde{a}^{\e,\xi} _{\lambda} (u) \big)\Big[\int_{0}^{+\infty} \overline{ H^{(1)} _{\lambda}\Big(\tilde{a}^{\e,\xi}  (u),\mathbb{E}[C^{1/2,\e}(0) \otimes C^{1/2,\e}(w)],\frac{u}{\e},\frac{u}{\e}+w \Big) dw +G_\lambda\Big(\tilde{a}^{\e,\xi} (u)\Big)}\Big]
\end{split}\]
\[\begin{split}
\hspace{0.7cm}&+\partial^2_v f\big(\tilde{a}^{\e,\xi} _{\lambda} (u) \big)\int_{0}^{+\infty} H^{(2)} _{\lambda}\Big(\tilde{a}^{\e,\xi}  (u),\mathbb{E}[C^{1/2,\e}(0) \otimes C^{1/2,\e}(w)],\frac{u}{\e},\frac{u}{\e}+w \Big) dw\\
&+\partial^2_{\overline{v}} f\big(\tilde{a}^{\e,\xi} _{\lambda} (u) \big)\int_{0}^{+\infty} \overline{ H^{(2)} _{\lambda}\Big(\tilde{a}^{\e,\xi}  (u),\mathbb{E}[C^{1/2,\e}(0) \otimes C^{1/2,\e}(w)],\frac{u}{\e},\frac{u}{\e}+w \Big)}dw\\
&+\partial_{\overline{v}}\partial_v f\big(\tilde{a}^{\e,\xi} _{\lambda} (u) \big)\int_{0}^{+\infty} H^{(3)} _{\lambda}\Big(\tilde{a}^{\e,\xi}  (u),\mathbb{E}[C^{1/2,\e}(0) \otimes C^{1/2,\e}(w)],\frac{u}{\e},\frac{u}{\e}+w \Big) dw\\
&+\partial_v \partial_{\overline{v}} f\big(\tilde{a}^{\e,\xi} _{\lambda} (u) \big)\int_{0}^{+\infty} \overline{ H^{(3)} _{\lambda}\Big(\tilde{a}^{\e,\xi}  (u),\mathbb{E}[C^{1/2,\e}(0) \otimes C^{1/2,\e}(w)],\frac{u}{\e},\frac{u}{\e}+w \Big)}dw\\
&+o(1).
\end{split}
\]
\end{lem}
\begin{preuve}[of Lemma \ref{A1P2}]
This lemma is a consequence of a long but straightforward computation. The terms of order one in $\mathcal{A}^\e \big(f^\e _0 +f^\e _1+f^\e _2\big)$ are uniformly integrable. Moreover, the terms bring by $f^\e _2$ in $o(1)$ can be bounded using the bound obtain in the proof of Lemma \ref{f2P2} with $\se$ instead of $\e$. Therefore, the unbounded logarithm which was a problem in $\mathcal{A}^\e \big(f^\e _0 +f^\e _1)$ can be killed. $\square$ That completes the proof of Proposition \ref{tightth1P2}. 
$\blacksquare$
\end{preuve}
\end{preuve}
Now, we characterize all the subsequence limits (Proposition \ref{subseq}), using first by a complex martingale problem (Lemma \ref{propmgP2}), and then a real martingale problem (Lemma \ref{uniquemgP2}) to use the classical uniqueness results in \cite{inf} for instance. 
\begin{prop}\label{subseq}
All the subsequence limits of $(\tilde{a}^{\e,\xi}  (.))_{\epsilon \in (0,1)}$ are solution of a well posed martingale problem associated to the infinitesimal generator defined by \eqref{Ls}.
\end{prop}
\begin{preuve}[of Proposition \ref{subseq}]
To do that, we consider a converging subsequence of $(\tilde{a}^{\e,\xi}  (.))_{\epsilon \in (0,1)}$ which converges to a limit $\widehat{a}^{\xi} (.)$. For the sake of simplicity we denote by $(\tilde{a}^{\e,\xi}  (.))_{\epsilon \in (0,1)}$ such a subsequence. To exhibit the well-posed martingale problem, we need first to extract all the fast oscillating phase, as done in \cite[Proposition 6.5]{gomez2bis}, in $\mathcal{A}^\e \big(f^\e _0 +f^\e _1+f^\e _2\big)$ using Assumption \ref{hyp5}.
\begin{lem}\label{propmgP2}
For all $\lambda \in \mathcal{H}_\xi$ and for all $f$ smooth test function,
\[
\begin{split}
f\big(\widehat{a}^{\xi}_{\lambda}&(u)\big)-\int_0 ^u  \partial_{v}f\big(\widehat{a}^{\xi}_{\lambda}(w)\big)\left<J^\xi(\widehat{a}^{\xi}(w)),\lambda\right>_{\mathcal{H}_\xi} +\partial_{\overline{v}}f\big(\widehat{a}^{\xi}_{\lambda}(w)\big) \overline{\left<J^\xi(\widehat{a}^{\xi}(w)),\lambda\right>_{\mathcal{H}_\xi}}\\
+&\partial ^2 _{v} f\big(\widehat{a}^{\xi}_{\lambda}(w)\big)\left<K \big(\widehat{a}^{\xi}(w)\big) (\lambda),\lambda \right>_{\mathcal{H}_\xi}+\partial ^2_{\overline{v}}f\big(\widehat{a}^{\xi}_{\lambda}(w)\big) \overline{\left<K \big(\widehat{a}^{\xi}(w)\big) (\lambda),\lambda \right>_{\mathcal{H}_\xi}}\\
+&\partial  _{\overline{v}} \partial  _{v}f\big(\widehat{a}^{\xi}_{\lambda}(w)\big)\left<L\big(\widehat{a}^{\xi}(w)\big) (\lambda),\lambda \right>_{\mathcal{H}_\xi}+ \partial  _{v} \partial  _{\overline{v}}f\big(\widehat{a}^{\xi}_{\lambda}(w)\big)\overline{\left<L \big(\widehat{a}^{\xi}(w)\big) (\lambda),\lambda \right>_{\mathcal{H}_\xi}}dw
\end{split}
\]
is a martingale, where
\begin{equation}\label{mgfinal}
\begin{split}
J^\xi(\emph{\textbf{T}})_{j}&= \left[\frac{\Gamma^{c,s} _{jj}}{2}-\Lambda^{c,s,\xi}_j+i\left(\frac{\Gamma^{s,s}_{jj}}{2}-\Lambda^{s,s,\xi}_j+\kappa^{s,\xi}_j\right)\right]\emph{\textbf{T}}_j,\\
K (\emph{\textbf{T}}) (\lambda)_{j}&= - \frac{1}{2}\sum_{\substack{l=1\\l\not=j}}^N\left( \Gamma^{c,s}_{jl}+i\Gamma^{s,s}_{jl} \right)\emph{\textbf{T}}_j \emph{\textbf{T}}_l\overline{\lambda_l},\qquad L (\emph{\textbf{T}}) (\lambda)_{j}= \frac{1}{2}\sum_{\substack{l=1\\l\not=j}}^N \Gamma^{c,s}_{jl} \emph{\textbf{T}}_l \overline{\emph{\textbf{T}}_l}\lambda_j,
\end{split}
\end{equation}
and $J^\xi(\emph{\textbf{T}})_{\ga}=K (\emph{\textbf{T}}) (\lambda)_{\ga}=L (\emph{\textbf{T}}) (\lambda)_{\ga}=0$ for almost every $\ga \in (\xi,k^2)$, and for $(\emph{\textbf{T}},\lambda)\in \mathcal{H}_\xi^2$. Here,  $\Gamma^{1,s}$, $ \Gamma^{c,s}$,  $\Gamma^{s,s}$,  $\Lambda^{c,s,\xi}$, $\Lambda^{s,s,\xi}$, and $\kappa^{s,\xi}$ are defined in Section \ref{randomsurfsec}.
\end{lem} 
In order to prove uniqueness, we decompose  $\widehat{a}^{\xi} (.)$ into real and imaginary parts. Then, let us consider the new process
\[\textbf{Y} ^\xi(u)=\begin{bmatrix} \textbf{Y}^{1,\xi} (u) \\ \textbf{Y}^{2,\xi}(u) \end{bmatrix},\qquad \text{where}\qquad   \textbf{Y}^{1,\xi} (u)=Re\big( \widehat{a}^{\xi}(u)\big) \qquad \text{and}\quad \textbf{Y}^{2,\xi} (u)=Im\big(\widehat{a}^{\xi}(u) \big).\]
This new process takes its values in $\mathcal{G}_\xi\times \mathcal{G}_\xi$, where $\mathcal{G}_\xi=\mathbb{R}^N\times L^2((\xi,k^2),\mathbb{R})$, and we introduce the operator 
\begin{equation*}
\begin{split}
\Upsilon:\,&\mathcal{G}_\xi\times \mathcal{G}_\xi\longrightarrow \mathcal{G}_\xi\times \mathcal{G}_\xi,\\
&\begin{bmatrix} \textbf{T}^1 \\ \textbf{T}^2 \end{bmatrix} \longmapsto \begin{bmatrix} \textbf{T}^2 \\ -\textbf{T}^1 \end{bmatrix}.
\end{split}
\end{equation*}

\begin{lem}\label{uniquemgP2}
For all $f \in \mathcal{C}^2 _b( \mathcal{G}_\xi\times \mathcal{G}_\xi)$,
\[
M^\xi_f(u)=f(\emph{\textbf{Y}}^\xi(u))-\int_0 ^u L^\xi f(\emph{\textbf{Y}}^\xi(w))dw
\]
is a continuous martingale, where for all $(\emph{\textbf{Y}},\lambda) \in(\mathcal{G}_\xi\times \mathcal{G}_\xi)^2$
\begin{equation}\label{genereal} L^\xi f(\emph{\textbf{Y}})=\frac{1}{2}trace\left( A(\emph{\textbf{Y}})D^2 f(\emph{\textbf{Y}}) \right) +\left< B^\xi(\emph{\textbf{Y}}), Df(\emph{\textbf{Y}})\right>_{\mathcal{G}_\xi\times \mathcal{G}_\xi},\end{equation}
with $A(\emph{\textbf{Y}})(\lambda)=A_1(\emph{\textbf{Y}})(\lambda)+A_2(\emph{\textbf{Y}})(\lambda)$. Moreover, for $j\in\{1,\dots,N\}$,
\[\begin{split}
B^\xi(\emph{\textbf{Y}})_j&=\left[\frac{\Gamma^{c,s} _{jj}}{2}-\Lambda^{c,s,\xi}_j\right]\emph{\textbf{Y}}_j-\left[\frac{\Gamma^{s,s}_{jj}}{2}+\kappa^{s,\xi}_j-\Lambda^{s,s,\xi}_j\right]\Upsilon_j(\emph{\textbf{Y}})\\
A_1(\emph{\textbf{Y}})(\lambda)_j&=-\emph{\textbf{Y}}_j\sum_{\substack{l=1\\ l\not=j}}^N \Gamma^{c,s}_{jl}\big[\emph{\textbf{Y}}^1_l\lambda^1_l+\emph{\textbf{Y}}^2_l\lambda^2_l\big]+\Upsilon_j(\emph{\textbf{Y}})\sum_{\substack{l=1\\ l\not=j}}^N \Gamma^{c,s}_{jl}\big[\Upsilon^1_l(\emph{\textbf{Y}})\lambda^1_l+\Upsilon^2_l(\emph{\textbf{Y}})\lambda^2_l\big]\\
A_2(\emph{\textbf{Y}})(\lambda)_j&=\lambda_j\sum_{\substack{l=1\\l\not=j}}^N\Gamma^{c,s}_{jl}\big[(\emph{\textbf{Y}}^1_l)^2+(\emph{\textbf{Y}}^2_l)^2\big],
\end{split}\]
and $B^\xi_\ga(\emph{\textbf{Y}})=A_\ga(\emph{\textbf{Y}})(\lambda)=A_\ga(\emph{\textbf{Y}})(\lambda)=A_\ga(\emph{\textbf{Y}})(\lambda)=0$
for almost every $\ga\in(\xi,k^2)$.
Moreover, the martingale problem associated to the generator $L^\xi$ is well-posed.
\end{lem}
\begin{preuve}[of Lemma \ref{uniquemgP2}]
Following the proof of Theorem 4.1.4 in \cite{inf}, to prove that $M^\xi_f$ is a martingale it suffices to show that
\[
\big< M^\xi (u), \lambda\big>_{\mathcal{G}_\xi\times\mathcal{G}_\xi}=M^\xi_\lambda (u)= \big< \textbf{Y}^\xi (u)-\int_0 ^u B^\xi(\textbf{Y}^\xi(w))dw ,\lambda\big>_{\mathcal{G}_\xi\times \mathcal{G}_\xi}
\]
is a continuous martingale with quadratic variation 
\[<M^\xi_\lambda >(u)=\int _0 ^u \big< A(\textbf{Y}^\xi(w))(\lambda),\lambda\big>_{\mathcal{G}_\xi\times \mathcal{G}_\xi}dw.\]
Moreover, for all $(\textbf{Y},\lambda) \in (\mathcal{G}_\xi\times \mathcal{G}_\xi)^2$, we have
$\left< A(\textbf{Y})(\lambda),\lambda\right>_{\mathcal{G}_\xi\times \mathcal{G}_\xi}\geq 0$ and $trace(A(\textbf{Y}))<+\infty$. According to Theorem 3.2.2 and 4.4.1 in \cite{inf}, the martingale problem associated to $L^\xi$ is therefore well-posed.
$\square$
\end{preuve}
Finally, to recover rigorously the form of the generator \eqref{geneapprox} for the complex diffusion process $\widehat{a}^{\xi}(\omega,\cdot)$, we use formula \eqref{genereal} where
 \[  g(\textbf{Y})=f(\textbf{T})\quad \text{and}\quad \textbf{T}=\textbf{Y}^1+i\textbf{Y}^2,\]
with 
\[\partial_{\textbf{Y}^1_j}g(\textbf{Y})=(\partial_{\textbf{T}_j}+\partial_{\overline{\textbf{T}_j}})f(\textbf{T})\quad\text{and}\quad\partial_{\textbf{Y}^2_j}g(\textbf{Y})=i(\partial_{\textbf{T}_j}-\partial_{\overline{\textbf{T}_j}})f(\textbf{T}).\]
In fact, for the drift term, we have
\[\begin{split}
\left< B^\xi(\emph{\textbf{Y}}), Dg(\textbf{Y})\right>_{\mathcal{G}_\xi\times \mathcal{G}_\xi}=&\sum_{j=1}^N \left[\frac{\Gamma^{c,s} _{jj}}{2}-\Lambda^{c,s,\xi}_j\right]\textbf{Y}^1_j\partial_{\textbf{Y}^1_j}g(\textbf{Y})-\left[\frac{\Gamma^{s,s}_{jj}}{2}+\kappa^{s,\xi}_j-\Lambda^{s,s,\xi}_j\right]\textbf{Y}^2_j\partial_{\textbf{Y}^1_j}g(\textbf{Y})\\
&+\left[\frac{\Gamma^{c,s} _{jj}}{2}-\Lambda^{c,s,\xi}_j\right]\textbf{Y}^2_j\partial_{\textbf{Y}^2_j}g(\textbf{Y})+\left[\frac{\Gamma^{s,s}_{jj}}{2}+\kappa^{s,\xi}_j-\Lambda^{s,s,\xi}_j\right]\textbf{Y}^1_j\partial_{\textbf{Y}^2_j}g(\textbf{Y})\\ 
=& \sum_{j=1}^N   \left[\frac{\Gamma^{c,s} _{jj}}{2}-\Lambda^{c,s,\xi}_j  +i\left(  \frac{\Gamma^{s,s}_{jj}}{2}+\kappa^{s,\xi}_j-\Lambda^{s,s,\xi}_j \right)    \right] \textbf{T}_j \partial_{\textbf{T}_j}g+c.c,
\end{split}\]
where $c.c$ means complex conjugate of the previous term. To obtain the diffusion term, a formal but fast way to do consists in replacing the component of the test function $\lambda_j$ by $\partial_{\overline{\textbf{T}_j}}$ in $\left<K \big(\widehat{a}^{\xi}(w)\big) (\lambda),\lambda \right>_{\mathcal{H}_\xi}$ and $\left<L\big(\widehat{a}^{\xi}(w)\big) (\lambda),\lambda \right>_{\mathcal{H}_\xi}$ defined by \eqref{mgfinal}. In fact, in the proof of Proposition \ref{uniquemgP2}, one can see that the real part $\lambda^1_j$ corresponds to $\partial_{\textbf{Y}^1_j}$ and the imaginary $\lambda^2_j$ corresponds to $\partial_{\textbf{Y}^2_j}$, so that $\lambda_j=\lambda^1_j+i\lambda^2_j$ correspond to $\partial_{\textbf{Y}^1_j}+i\partial_{\textbf{Y}^2_j}=2\partial_{\overline{\textbf{T}_j}}$. 

Finally, using the fact that
\[R_s(z+w-w')=\frac{1}{2\pi}\int du\,  \widehat{R}_s(u) e^{iuz}e^{iu(w-w')},\]
and $\int \frac{e^{-iwu}}{A^2+w^2}dw=\frac{\pi}{A}e^{-A\vert u\vert}$, $\int \frac{we^{-iwu}}{(A^2+w^2)^2}dw=-iu\frac{\pi}{2A}e^{-A\vert u\vert}$, $\int \frac{e^{-iwu}}{(A^2+w^2)^2}dw=\frac{\pi(1+A\vert u\vert)}{2A^3}e^{-A\vert u\vert}$ for $A>0$, we obtain after long but straightforward computations the expression of the diffusion coefficients in Theorem \ref{thasymptP21}.
That completes the proof of Proposition \ref{subseq}. $\blacksquare$
\end{preuve}

\subsection{Proof of Proposition \ref{coefdec}}\label{proofcoefdec}

The proof of this proposition follows the idea of the one of \cite[Proposition 5.1]{alonso}. We restrict the proof to the case of a randomly perturbed bottom. For the case of a randomly perturbed surface the proof remains the same but with more terms and lengthier computations. Proposition \ref{coefdec} is a consequence of the three following lemmas allowing us to study the behavior of the coupling coefficients $\Gamma^{1,b}_{jj}(\omega)$, $\Gamma^{c,b} _{jj}(\omega)$, and dissipation coefficients $\Lambda^{c,b}_j (\omega)$ in the asymptotic $k(\omega)\to+\infty$. However, we restrict the proof to the ones of $\Gamma^{c,b} _{jj}(\omega)$, and $\Lambda^{c,b}_j (\omega)$ since for $\Gamma^{1,b}_{jj}(\omega)$ it is a simple application of Lemma \ref{asympsigma}.

\begin{lem}\label{asympgamma}
Let us denote 
\[\tilde{ \Gamma}_j(\omega)=\sum_{\substack{l=1 \\ l\not = j} }^{\N{}} \frac{\phi^2_l(\omega,d)}{\beta_l(\omega)}I_b( \beta_j(\omega)-\beta_l(\omega) ) .\]
We have the three following asymptotic behaviors.
 \begin{enumerate}
 \item For  $\eta_1\in[0,1/2)$ and $\nu_1> 0$, we have
\[  \sup_{j\in\{1,\dots,[\nu_1\N{}^{\eta_1}]\}}\Big \lvert\tilde{ \Gamma}_j(\omega)- \frac{\sqrt{2}}{\pi\theta^2 (n_1 k(\omega))^{3/2}}\int_{0}^{+\infty}\sqrt{v}I_b(v)dv\Big \rvert =o\Big(\frac{1}{\N{}^{3/2}}\Big).\]
\item For $\eta_2\in(1/2,1]^2$ and $\mu_2>\nu_2> 0$ (and $\nu_2<\mu_2\leq 1$ if $\eta_2=1$), we have
\[  \sup_{j\in\{[\nu_2\N{}^{\eta_2}],\dots,[\mu_2\N{}^{\eta_2}]\}}\Big \lvert \tilde{ \Gamma}_j(\omega)-\frac{2j}{\pi^2  \N{}^2}\int_{-\infty}^{+\infty}I_b(v)dv  \Big \rvert =o\Big(\frac{1}{\N{}^{2-\eta_2}}\Big).\]
 \item For $j=[\nu \N{}^{1/2}]$, we have
\[ \Big \lvert\tilde{ \Gamma}_j(\omega)- \frac{2}{\pi\theta^2 (n_1 k(\omega))^{3/2}}\int_{-\nu^2 \pi \theta/(2d)}^{+\infty}\sqrt{\frac{\nu^2\theta \pi}{d} +2v}I_b(v)dv\Big \rvert =o\Big(\frac{1}{\N{}^{3/2}}\Big).\]
\end{enumerate}
\end{lem}
Regarding the dissipation coefficients, we have the following result.
\begin{lem}\label{asymplambda}
Let $\eta\in[0,1]$, $\nu_1>0,\nu_2\geq0$ with $\nu_1<1$ if $\eta=1$, and let us denote
\[ \tilde{\Lambda}_j(\omega)=\int_0^{k^2(\omega)} \frac{\phi^2_{\ga'}(\omega,d)}{\sqrt{\ga'}}I_b(\beta_j(\omega)-\sqrt{\ga'})d\ga'.\]
We have the two following asymptotic behaviors. 
\begin{enumerate}
\item For $j=[\nu\N{}^{\eta}]$ with $\eta\in[0,1]$ and $\nu> 0$ ($\nu < 1$ if $\eta=1$), we have
\[\frac{C_1}{(n_1k(\omega))^{\alpha_I}} \leq \tilde{\Lambda}_j\leq \frac{C_2}{(n_1k(\omega))^{\alpha_I}},\]
where $C_1$ and $C_2$ are two positive constants.
\item For $j=\N{}-[\nu]$ with $\nu> 0$, we have
\[\frac{C_3}{\N{}^{1/2}} \leq \tilde{\Lambda}_j\leq \frac{C_4}{\N{}^{1/2}}.\]
\end{enumerate}
\end{lem}
\begin{preuve}[of Lemma \ref{asympgamma}]
To prove this results, we need the following lemma.
\begin{lem}\label{asympsigma}
We have for all $l\in\{1,\dots,\N{}\}$
\[\begin{split}
\frac{\theta}{\pi(\N{}+1)}\frac{2}{1+\frac{1}{\pi \N{} \sqrt{1-\frac{(l-1)^2}{(\N{}+1)^2}}}} & \frac{\frac{(l-1)^2}{(\N{}+1)^2}}{\sqrt{1-\theta^2\frac{(l-1)^2}{(\N{}+1)^2}}} \leq \frac{\phi^2_l(\omega,d)}{\beta_l(\omega)}\leq \frac{2\theta}{\pi\N{}}  \frac{\frac{l^2}{\N{}^2}}{\sqrt{1-\theta^2\frac{l^2}{\N{}^2}}}.
\end{split}\]
\end{lem}
\begin{preuve}[of Lemma \ref{asympsigma}]
In fact, using relation \eqref{eqsigma} and the relation $\sin(\arctan(x))=x/\sqrt{1+x^2}$ we have
\[
\sin(\sigma_j(\omega))=-\sin(\arctan(\sigma_j(\omega)/\sqrt{(n_1k(\omega)d\theta)^2-\sigma_j^2(\omega)}))=-\sigma_j(\omega)/(n_1k(\omega)d\theta).
\]
In the same way, using the relation $\cos(\arctan(x))=1/\sqrt{1+x^2}$, we have
\[A^2_j(\omega)=\frac{2/d}{1+\frac{1}{n_1k(\omega)d\theta} \sqrt{1-\sigma^2_j(\omega)/(n_1k(\omega)d\theta)^2}}\]
where $A_j(\omega)$ is defined by \eqref{Aj}. Therefore, using the fact that $\sigma_j(\omega)\in[(j-1)\pi,j\pi]$ and the definitions \eqref{sigj} and \eqref{number}, we obtain the desired result.
$\square$
\end{preuve}
According to Lemma \ref{asympsigma}  we have for all $j\in\{1,\dots,\N{}\}$ the two following inequalities 
 \begin{equation}\label{doubleineq}\begin{split}
\frac{\theta}{\pi(\N{}+1)}& \sum_{\substack{l=1 \\ l\not = j} }^{\N{}}\frac{2}{1+\frac{1}{\pi \N{} \sqrt{1-(l-1)^2/(\N{}+1)^2}}} \frac{(l-1)^2/(\N{}+1)^2}{\sqrt{1-\theta^2(l-1)^2/(\N{}+1)^2}}\\
&\hspace{3cm}\times I_b( n_1k(\omega)(\sqrt{1-\theta^2 (j-1)^2/(\N{}+1)^2} -\sqrt{1-\theta^2 l^2/\N{}^2} ))   \\
&\leq  \sum_{\substack{l=1 \\ l\not = j} }^{\N{}} \frac{\phi^2_l(\omega,d)}{\beta_l(\omega)}I_b( \beta_j(\omega)-\beta_l(\omega) ) \leq\frac{2\theta}{\pi\N{}} \sum_{\substack{l=1 \\ l\not = j} }^{\N{}} \frac{l^2/\N{}^2}{\sqrt{1-\theta^2l^2/\N{}^2}}\\
&\hspace{2cm}\times I_b( n_1k(\omega)(\sqrt{1-\theta^2 j^2/\N{}^2} -\sqrt{1-\theta^2 (l-1)^2/(\N{}+1)^2} )) 
\end{split}\end{equation}
so that for all $\eta\in (0,1)$ and $j\in \{1,\dots,[\nu \N{}^\eta]\}$
\begin{equation}\label{estimatej}\begin{split}
\Big\vert \sum_{\substack{l=1 \\ l\not = j} }^{\N{}} &\frac{\phi^2_l(\omega,d)}{\beta_l(\omega)}I( \beta_j(\omega)-\beta_l(\omega)) -\frac{2\theta}{\pi}\int_0^1 \frac{y^2}{\sqrt{1-\theta^2y^2}}I_b( n_1k(\omega)(\sqrt{1-\theta^2 j^2/\N{}^2} -\sqrt{1-\theta^2 y^2} ))dy  \Big\vert\\
&\leq C\Big[\frac{1}{\N{}}\Big[\frac{j+\N{}^\alpha}{\N{}}\Big]^2+\frac{1}{\N{}}\frac{1}{[\N{}^{\alpha-1}(\N{}^{\alpha}+2j)]^{\alpha_I }}\Big].
 \end{split}\end{equation}
Here $\alpha_I$ is the decaying power of the function $I_b$. To obtain this estimate we have split the set of indexes into two parts $\vert j-l\vert\leq \N{}^\alpha$ and $\vert j-l\vert > \N{}^\alpha$. The reason is that for the first set of indexes the term $I_b( n_1k(\omega)(\sqrt{1-\theta^2 j^2/\N{}^2}-\sqrt{1-\theta^2 l^2/\N{}^2} ))$ is considered to be bounded, but decays fast for the second set  of indexes. Moreover, for all $\eta_1\in (0,1/2)$ and $j\in \{1,\dots,[\nu_1\N{}^{\eta_1}]\}$, using a change of variable, we have
\[\begin{split}
\frac{2\theta}{\pi}\int_0^1& \frac{y^2}{\sqrt{1-\theta^2y^2}}I_b( n_1k(\omega)(\sqrt{1-\theta^2 j^2/\N{}^2}-\sqrt{1-\theta^2 y^2})) dy\\
&=\frac{2}{\pi\theta^2 n_1k(\omega)}\int_{n_1k(\omega)(\sqrt{1-\theta^2j^2/\N{}^2}-1)}^{n_1k(\omega)(\sqrt{1-\theta^2j^2/\N{}^2}-\sqrt{1-\theta^2})}\sqrt{1-\Big(\sqrt{1-\theta^2j^2/\N{}^2}-v/(n_1k)\Big)^2}I_b(v)dv
\end{split}\]
so that
\[\begin{split}
\Big\vert \frac{2\theta}{\pi}\int_0^1& \frac{y^2}{\sqrt{1-\theta^2y^2}}I_b( n_1k(\omega)(\sqrt{1-\theta^2 j^2/\N{}^2}-\sqrt{1-\theta^2 y^2})) dy\\
&\hspace{2cm}-\frac{\sqrt{2}  (1-\theta^2 j^2/\N{}^2)^{1/4}  }{\pi\theta^2 (n_1 k(\omega))^{3/2}}\int_{n_1k(\omega)(\sqrt{1-\theta^2j^2/\N{}^2}-1)}^{n_1k(\omega)(\sqrt{1-\theta^2j^2/\N{}^2}-\sqrt{1-\theta^2})}\sqrt{\vert v\vert}I_b(v)dv\Big\vert\\
\leq&\frac{C}{(n_1 k(\omega))^{3/2}} \int_{n_1k(\omega)(\sqrt{1-\theta^2j^2/\N{}^2}-1)}^{n_1k(\omega)(\sqrt{1-\theta^2j^2/\N{}^2}-\sqrt{1-\theta^2})}\Big(\frac{j^2}{\N{} v}+\frac{v^2}{n_1k(\omega)}\Big)\sqrt{\vert v\vert }I_b(v)dv,
\end{split}\]
and choosing properly $\alpha\in(1/2+1/(4\alpha_I),3/4)$ in \eqref{estimatej}, all the previous error estimates are negligible with respect to $\N{}^{-3/2}$,
\[
\sup_{j\in \{1,\dots,[\nu_1\N{}^{\eta_1}]\}}\Big\vert \tilde{ \Gamma}_j(\omega) -\frac{\sqrt{2}}{\pi\theta^2 (n_1 k(\omega))^{3/2}}\int_{0}^{+\infty}\sqrt{v}I_b(v)dv\Big\vert=o\Big(\N{}^{-3/2}\Big).
\]
Now, for all $\eta_2\in (1/2,1)$ and $j\in \{[\nu_2\N{}^{\eta_2}],\dots,[\mu_3\N{}^{\eta_2}]\}$, using a change of variable, we have
\[\begin{split}
\frac{2\theta}{\pi}\int_0^1& \frac{y^2}{\sqrt{1-\theta^2y^2}}I_b( n_1k(\omega)(\sqrt{1-\theta^2 j^2/\N{}^2}-\sqrt{1-\theta^2 y^2})) dy\\
&=\frac{2}{\pi\theta^2 n_1k(\omega)}\int_{n_1k(\omega)(\sqrt{1-\theta^2j^2/\N{}^2}-1)}^{n_1k(\omega)(\sqrt{1-\theta^2j^2/\N{}^2}-\sqrt{1-\theta^2})}\sqrt{1-\Big(\sqrt{1-\theta^2j^2/\N{}^2}-v/(n_1k)\Big)^2}I_b(v)dv
\end{split}\]
so that
\[\begin{split}
\Big\vert \frac{2\theta}{\pi}\int_0^1& \frac{y^2}{\sqrt{1-\theta^2y^2}}I_b( n_1k(\omega)(\sqrt{1-\theta^2 j^2/\N{}^2}-\sqrt{1-\theta^2 y^2})) dy\\
&\hspace{4cm}-\frac{2jd}{\pi^2  \N{}^2}\int_{n_1k(\omega)(\sqrt{1-\theta^2j^2/\N{}^2}-1)}^{n_1k(\omega)(\sqrt{1-\theta^2j^2/\N{}^2}\sqrt{1-\theta^2})}I_b(v)dv\Big\vert\\
\leq&\frac{Cj}{\N{}^2} \int_{n_1k(\omega)(\sqrt{1-\theta^2j^2/\N{}^2}-1)}^{n_1k(\omega)(\sqrt{1-\theta^2j^2/\N{}^2}\sqrt{1-\theta^2})}\Big(\frac{v^2}{j^2}+\frac{vN}{j^2}\Big)I_b(v)dv,
\end{split}\]
and
\[
\Big\vert \int_{n_1k(\omega)(\sqrt{1-\theta^2j^2/\N{}^2}-1)}^{n_1k(\omega)(\sqrt{1-\theta^2j^2/\N{}^2}-\sqrt{1-\theta^2})}I_b(v)dv -\int_{-\infty}^{+\infty}I_b(v)dv \Big\vert \leq C \Big(\frac{\N{}}{j^2}\Big)^{\alpha_I-1}.
\]
Then, choosing properly $\alpha\in(\eta_2,(1+\eta_2)/2)$ and $\alpha>1/2+(1-\eta_2)/(2\alpha_I)$ in \eqref{estimatej}, we obtain
\[
\sup_{j\in \{[\nu_2\N{}^{\eta_2]},\dots,[\mu_2\N{}^{\eta_2}]\}}\Big\vert\tilde{ \Gamma}_j(\omega)-\frac{2jd}{\pi^2  \N{}^2}\int_{-\infty}^{+\infty}I_b(v)dv\Big\vert=o\Big(\N{}^{-2+\eta_2}\Big).
\]
Finally, for all $\nu\in (0,1)$ and $j\in \{[\nu\N{}],\dots,\N{}\}$, we have using \eqref{doubleineq}
\[\begin{split}
\Big\vert \tilde{ \Gamma}_j(\omega) -\frac{2\theta}{\pi}\int_0^1 \frac{y^2}{\sqrt{1-\theta^2y^2}}&I_b( n_1k(\omega)(\sqrt{1-\theta^2 j^2/\N{}^2} -\sqrt{1-\theta^2 y^2} ))dy  \Big\vert\\
&\leq \frac{C}{\N{}^{2-\alpha}}+ \frac{C}{\N{}^{1+\alpha_I(1-\alpha)}}.
\end{split}\]
The exponent $\alpha\in(0,1)$ comes from the fact that we have split the sum between $\{1,\dots,[\nu\N{}]-[\N{}^\alpha]\}$ and $\{[\nu\N{}]-[\N{}^\alpha],\dots,\N\}$ in the same sprit of \eqref{estimatej}. As in the previous case, we obtain
\[
\sup_{j\in\{\N{}-[\N{}^\eta],\dots,\N{} \}}\Big\vert \tilde{ \Gamma}_j(\omega)-\frac{2j}{\pi^2  \N{}^2}\int_{-\infty}^{+\infty}I_b(v)dv\Big\vert=o\Big(\N{}^{-1}\Big),
\]
which concludes the proof of Lemma \ref{asympgamma}. The case $\eta_3=1/2$ needs just a more accurate truncated expansion but brings no significant difficulty. 
$\square$
\end{preuve}
\begin{preuve}[of Lemma \ref{asymplambda}]
We have
\[\begin{split}
\tilde{\Lambda}_j&(\omega)=\int_0^{k^2(\omega)} \frac{\phi^2_{\ga'}(\omega,d)}{\sqrt{\ga'}}I_b(\beta_j(\omega)-\sqrt{\ga'})d\ga'\\
&=\frac{2}{\pi}\int_1^{1/\theta} \frac{v\sqrt{v^2-1}\sin^2(n_1k(\omega)d\theta v)}{(v^2-\sin^2(n_1k(\omega)d\theta v))\sqrt{1-\theta^2 v^2}}I_b\big( n_1k(\omega)\big( \sqrt{1-\frac{\sigma^2_j(\omega)}{(n_1 k(\omega)d)^2}}-\sqrt{1-\theta^2v^2} )  )dv,
\end{split}\]
so that 
\[J_1\leq \tilde{\Lambda}_j(\omega)\leq J_2,\]
where
\[\begin{split}
J_1&=\frac{2}{\pi}\int_1^{1/\theta}  \frac{\sqrt{v^2-1} \sin^2(n_1k(\omega)d\theta v)}{v\sqrt{1-\theta^2 v^2}}I_b( n_1k(\omega)( \sqrt{1-\sigma^2_j(\omega)/(n_1 k(\omega)d)^2}-\sqrt{1-\theta^2v^2} )  )dv\\
J_2&=\frac{2}{\pi}\int_1^{1/\theta} \frac{v}{\sqrt{v^2-1}\sqrt{1-\theta^2 v^2}}I_b( n_1k(\omega)( \sqrt{1-\sigma^2_j(\omega)/(n_1 k(\omega)d)^2}-\sqrt{1-\theta^2v^2} ) )dv.
\end{split}\]
First, for $J_2$, using the change of variable $u=\sqrt{1-\theta^2}-\sqrt{1-\theta^2v^2}$ we have
\begin{equation}\label{J2}\begin{split}
J_2 &=\frac{C_b}{2\pi^2\theta(n_1k(\omega))^{\alpha_I}}\int_0^{\sqrt{1-\theta^2}} \frac{1}{\Big[ (\sqrt{1-\sigma^2_j(\omega)/(n_1 k(\omega)d)^2}-\sqrt{1-\theta^2}+u)^2 +\frac{1}{(n_1k(\omega))^2}\Big]^{\alpha_I/2}}  \\
&\hspace{4cm}\times\frac{du}{\sqrt{2u\sqrt{1-\theta^2}-u^2}}du,
\end{split}\end{equation}
so that for all $\eta\in[0,1)$
\[\begin{split}
\sup_{j\in\{1,\dots,[\nu_1\N{}^\eta]\}}\Big\vert J_2-  \frac{C_b}{2\pi^2\theta(n_1k(\omega))^{\alpha_I}}\int_{0}^{\sqrt{1-\theta^2}} &\frac{1}{(1-\sqrt{1-\theta^2}+u)^{\alpha_I}\sqrt{2u\sqrt{1-\theta^2}-u^2}} du\Big\vert\\
=o(\N{}^{-\alpha_I}),
\end{split}\]
and if $j=[\nu_1\N{}]$
\[\begin{split}
\Big\vert J_2-  &\frac{C_b}{2\pi^2\theta(n_1k(\omega))^{\alpha_I}}\int_{0}^{\sqrt{1-\theta^2}} \frac{1}{(\sqrt{1-\theta^2\nu_1^2}-\sqrt{1-\theta^2}+u)^{\alpha_I}\sqrt{2u\sqrt{1-\theta^2}-u^2}} du\Big\vert\\
&=o(\N{}^{-\alpha_I}).
\end{split}\]
Moreover, for the case $j$ of the form $\N{}-[\nu]$, we have
\[\begin{split}
J_2 &\leq  \tilde{J}_2(\eta)= \frac{C_b}{2\pi^2\theta(n_1k(\omega))^{\alpha_I}}\int_0^{\sqrt{1-\theta^2}} \frac{1}{\Big[ (\sqrt{1-\theta^2(1-\nu/\N{})^2}-\sqrt{1-\theta^2}+u)^2 +\frac{1}{(n_1k(\omega))^2}\Big]^{\alpha_I/2}}  \\
&\hspace{4cm}\times\frac{du}{\sqrt{2u\sqrt{1-\theta^2}-u^2}}du,
\end{split}\]
so that using the change of variable $v=u\N{}$ we obtain
\[\begin{split}
\Big\vert \tilde{J}_2(\eta)- & \frac{C_b d^{\alpha_I}\theta^{\alpha_I-1}}{2^{3/2}\pi^{\alpha_I+2}(1-\theta^2)^{1/4}\N{}^{1/2}} \int_{0}^{+\infty} \frac{1}{[(\frac{\theta^2}{\sqrt{1-\theta^2}}+v)^2+(\frac{d\theta}{\pi})^2]^{\alpha_I/2}\sqrt{v}} dv\Big\vert=o(\N{}^{-1/2}).
\end{split}\]
Now, for $J_1$, let us note that we have to take care of the sinus. To do so, we split the integral in many parts over which the sinus will be close to $1$,  
\[
J_1\geq \frac{2\theta}{\pi}I_b( n_1k(\omega) \sqrt{1-\pi^2/(n_1 k(\omega)d)^2} )\sum_{m=m_1}^{m_2}  \int_{\frac{\pi}{n_1k(\omega)d\theta}(1/2-\nu+m)}^{\frac{\pi}{n_1k(\omega)d\theta}(1/2+\nu+m)}\sqrt{v^2-1} \sin^2(n_1k(\omega)d\theta v)dv,\]
where $\nu_0$ is such that $\sin^2(\pi/2\pm\nu_0 \pi)\geq1/2$, and with
\[m_1=\big[  \frac{n_1k(\omega)d\theta}{\pi}+\nu_0-1/2\big]\qquad\text{and}\qquad m_2=\big[  \frac{n_1k(\omega)d}{\pi}-\nu_0-1/2\big].\]
Then, for the case $j=[\nu \N{}^\eta]$, with $\eta\in[0,1]$ and $\nu\in(0,1)$ if $\eta=1$, we have
\[
J_1 \geq \frac{\theta}{4\pi^2(n_1k(\omega))^{\alpha_I} } \frac{C_b}{\Big[ 1 +\frac{1}{(n_1k(\omega))^2}\Big]^{\alpha_I/2}} \frac{\nu_0}{\N{}+1}\sum_{m=m_1}^{m_2}  \sqrt{m^2/(\N{}+1)^2-1}
\]
with
\[ \frac{1}{\N{}+1}\sum_{m=m_1}^{m_2}  \sqrt{m^2/(\N{}+1)^2-1} \underset{ \N{} \gg1}{\sim} \int_{0}^{1/\theta-1}\sqrt{u(u+2)}du>0. \]
Now, in the case $j=\N{}-[\nu]$ with $\nu>0$, we have for $\mu\in(0,1)$
\[ \begin{split}
\tilde{\Lambda}_j&(\omega)\geq  \frac{2}{\pi} \sum_{m=1}^{[\mu\N{}]} \int_{1+\frac{2m+1}{2\N{}}-\frac{1}{\N{}^{3/2}}}^{1+\frac{2m+1}{2\N{}}+\frac{1}{\N{}^{3/2}}} \frac{v\sqrt{v^2-1}\sin^2(\N{} \pi v)}{(v^2-\sin^2(\N{}\pi v))\sqrt{1-\theta^2 v^2}}\\
&\hspace{2cm}\times I_b\big( n_1k(\omega)\big( \sqrt{1-\theta^2 j^2/\N{}^2}-\sqrt{1-\theta^2}\big) \big)dv\\
&\geq  \frac{1}{\pi} \sum_{m=1}^{[\mu\N{}]} \int_{1+\frac{2m+1}{2\N{}}-\frac{1}{\N{}^{3/2}}}^{1+\frac{2m+1}{2\N{}}+\frac{1}{\N{}^{3/2}}} \frac{v\sqrt{v^2-1}}{(v^2-\sin^2(\N{}\pi v))\sqrt{1-\theta^2 v^2}}\\
&\hspace{2cm}\times I_b\big( n_1k(\omega)\frac{\theta^2\nu}{\sqrt{1-\theta^2}\N{}}  \big)dv.
\end{split}\]
Moreover, using the fact that $\sin^2(\N{}\pi v)\geq 1-\sin^2(\pi /\N{}^{1/2})\geq 1- \pi^2 /\N{}$, we obtain  
\[ \begin{split}
\tilde{\Lambda}_j(\omega)&\geq \frac{C_b(1-\theta^2)^{(\alpha_I-1)/2}d^{\alpha_I}}{4 \pi^{\alpha_I+2}(\nu \theta)^{\alpha_I}}  \frac{1}{\big[ 1+\frac{\N{}^{2}\sqrt{1-\theta^2}}{(n_1 k(\omega)\theta^2)^2}\big]^{\alpha_I/2}} \\
&\hspace{3cm}\times \frac{1}{\N{}^{3/2}}  \sum_{m=1}^{[\mu\N{}]} \frac{\sqrt{\frac{2m+1}{2\N{}}-\frac{1}{\N{}^{3/2}}}}{\frac{2m+1}{2\N{}}-\frac{1}{\N{}^{3/2}}+\pi^2 /\N{}}
\end{split}\]
where
\[
\frac{1}{[\mu\N{}]}  \sum_{m=1}^{[\mu\N{}]}  \frac{\sqrt{\frac{2m+1}{2\N{}}-\frac{1}{\N{}^{3/2}}}}{\frac{2m+1}{2\N{}}-\frac{1}{\N{}^{3/2}}+\pi^2 /\N{}} \underset{ \N{} \gg1}{\sim} \frac{1}{\sqrt{\mu}}\int_0^1 \frac{du}{\sqrt{u}},
\]
which concludes the proof of Lemma \ref{asymplambda}.
$\square$
\end{preuve}


\begin{thebibliography}{20} 
\bibitem{axler}
{\sc S. J. Axler, P. Bourdon, and W. Ramey}, {\em Harmonic function theory}, Springer, New York, 2001.

\bibitem{alonso}
{\sc R. Alonso, L. Borcea, and J. Garnier}, \emph{Wave propagation in waveguides with random boundaries}, Commun. Math. Sci., 11 (2012), pp. 233--267. 

\bibitem{basson}
 {\sc   A. Basson and D. Grard-Varet}
 {\em Wall laws for fluid flows at a boundary with random roughness}, Comm. Pure Appl. Math., 61 (2008) ,pp.~941--987.

\bibitem{biling}
{\sc P. Billingsley}, \emph{Convergence of probability measure}, 2nd ed., Wiley InterScience, 1999.


\bibitem{carmona}
{\sc R.~Carmona and J.-P.~Fouque}, {\em Diffusion-approximation for the advection-diffusion of a passive scalar by a space-time gaussian velocity field}, Seminar on Stochastic Analysis, Random Fields and Applications. (Edited by E. Bolthausen, M. Dozzi and F. Russo). Birkhauser, Basel (1994) 37-50. 

\bibitem{donvar} 
{\sc M.D.~Donsker and S.R.S.~Varadhan}, {\em Asymptotic evaluation of certain Markov process expectations for large time, I-IV}, Comm. Pure Appl. Math., 28 (1975), pp.~1--47, 279--301, 29 (1979), pp.~368--461, 36 (1983) pp.~183--212.

\bibitem{dozier}
{\sc L. B. Dozier}, {\em A numerical treatment of rough surface scattering for the parabolic wave equation}, J. Acoust. Soc. Am., 75 (1984), pp. 1415--1432.

\bibitem{dya}
{\sc A. I. Dyachenko, E. A. Kuznetsov, M. D. Spector, and V. E. Zakharov}, \emph{Analytical description of the free surface dynamics of an ideal fluid (canonical formalism and conformal mapping)}, Phys. Lett. A, 221 (1996), pp. 73--79.

 
\bibitem{fouque}
{\sc J.-P.~Fouque}, {\em La convergence en loi pour les processus \`a valeur dans un espace nucl{\'e}aire}, Ann. Inst. Henri Poincar\'e, 20 (1984), pp.~225--245.

\bibitem{book}
{\sc J.-P.~Fouque, J.~Garnier, G.~Papanicolaou, and K.~S{\o}lna}, {\em Wave propagation and time reversal in randomly layered  media}, Springer, New York, 2007.

\bibitem{garniereva}
{\sc J.~Garnier}, {\em The role of evanescent modes in randomly perturbed single-mode waveguides }, Discrete and Continuous Dynamical Systems-Series B, 8 (2007), pp.~455--472.

\bibitem{garnier2}
{\sc J.~Garnier and K.~S{\o}lna}, {\em Effective transport equations and enhanced backscattering in random waveguides }, SIAM J. Appl. Math, 68 (2008), pp.~1574--1599.


\bibitem{garniernach}
{\sc J. Garnier, J. C. Munoz Grajales, and A. Nachbin},
{\em Effective behavior of solitary waves over random topography}, SIAM Multiscale Model. Simul., 6 (2007), pp. 995--1025. 

\bibitem{papa}
{\sc J.~Garnier and G.~Papanicolaou}, {\em Pulse propagation and time reversal in random waveguides}, SIAM J. Appl. Math, 67 (2007), pp.~1718--1739.

\bibitem{gomez2}
{\sc C.~Gomez}, {\em Wave propagation in shallow-water acoustic random waveguides}, Commun. Math. Sci.,  9 (2011), pp. 81--125.

\bibitem{gomez2bis}
{\sc C.~Gomez}, {\em Wave propagation in shallow-water acoustic random waveguides}, arXiv:0911.5646

\bibitem{gomez3}
{\sc C.~Gomez}, {\em Wave propagation and time reversal in random waveguides}, Ph.D dissertation, Universit\'e Paris-Diderot (Paris7), 2009. http://tel.archives-ouvertes.fr/tel-00439576/fr/.

\bibitem{papanicolaou}
{\sc W.~Kohler and G.~Papanicolaou}, {\em Wave propagation in randomly inhomogeneous ocean}, Lecture Notes in Physics, Vol. 70, J. B. Keller and J. S. Papadakis, eds., Wave Propagation and Underwater acoustics, Springer-Verlag, Berlin, 1977.

\bibitem{kat}
{\sc B. Katsnelson, V. Petnikov, and J. Lynch}, \emph{Fundamentals of shallow water acoustics}, Springer, New York, 2012.


\bibitem{khasminskii}
{\sc R. Z.~Khasminskii}, {\em A limit theorem for solutions of differential equations with random right hand side}, Theory Probab. Appl. 11 (1966), pp.~390--406.

\bibitem{kuperman}
{\sc W. A. Kuperman}, {\em Coherent component of specular reflection and transmission at a randomly rough two-fluid interface}, J. Acoust. Soc. Am., 58 (1975), pp. 365--370.

\bibitem{kuperman2}
{\sc W. A. Kuperman and F. Ingenito}, {\em Attenuation of the coherent component of sound propagating in shallow water with rough boundaries}, J. Acoust. Soc. Am., 61 (1977), pp. 1178--1187.

\bibitem{kuperman3}
{\sc W.A. Kuperman and J.F. Lynch}, \emph{Shallow-Water acoustics}, Phys. Today, 57 (October 2004), pp. 55--61.


\bibitem{kushner}
{\sc H.-J.~Kushner}, {\em Approximation and weak convergence methods for random processes}, MIT press, Cambridge, 1984.


\bibitem{kuttler}
{\sc J. R. Kuttler and J. D. Huffaker}, {\em Solving the parabolic wave equation with a rough surface boundary condition}, J. Acoust. Soc. Am., 94 (1993), pp. 2451--2453.


\bibitem{magnanini}
{\sc R.~Magnanini and F.~Santosa}, {\em Wave propagation in a 2-D optical waveguide}, SIAM J. Appl. Math, 61 (2000),
pp.~1237--1252.


\bibitem{marcuse}
{\sc D.~Marcuse}, {\em Theory of dielectric optical waveguides}, $2^{nd}$ ed., Academic press, New York, 1991.


\bibitem{mitoma}
{\sc I.~Mitoma}, {\em On the sample continuity of $\mathcal{S}'$-processes}, J. Math. Soc. Japan, 35 (1983), pp.~629--636.

\bibitem{nehari}
{\sc Z. Nehari}, \emph{Conformal mapping}, Dover Publications, New York, 2001.

\bibitem{pekeris}
{\sc C.I.~Pekeris}, {\em Theory of propagation of explosive sound in shallow water}, propagation of sound in the ocean, geological society of America, memoir 27 (1948), pp.~1--117. 


\bibitem{reed}
{\sc M. Reed and B. Simon}, \emph{Methods of Modern Mathematical Physics, Fourier Analysis, Self-Adjointness}, Academic Press, New York, 1980.

\bibitem{rowe}
{\sc H.E.~Rowe}, {\em Electromagnetic propagation in multi-mode random media}, Wiley, New York, 1999.


\bibitem{rudin}
{\sc W. Rudin}, {\em Real and complex analysis}, 3rd ed., McGraw-Hill, New York, 1987.

\bibitem{stroock}
{\sc D.W. Stroock and S.R.S. Varadhan}, \emph{Multidimensional diffusion processes}, Springer-Verlag,
Berlin, 1979.

\bibitem{titch}
{\sc E. C. Titchmarsh}, \emph{Introduction to the theory of Fourier integral}, The Clarendon Press, Oxford, 1948.

\bibitem{wilcox}
{\sc C.~Wilcox}, {\em Spectral analysis of the Pekeris operator in the theory of acoustic wave propagation in shallow water},  Arch. Rational Mech. Anal, 60 (1975/76), no. 3, pp.~259--300.

\bibitem{wilcox2}
{\sc C.~Wilcox}, {\em Transient electromagnetic wave propagation in a dielectric waveguide}, in Symposia Mathematica, vol. XVIII (Convegno sulla Teoria Matematica dell'Elettromagnetismo, INDAM, Rome, 1974), Academic Press, London, 1976, pp.~239--277. 

\bibitem{inf}
{\sc M.~Yor}, {\em Existence et unicit{\'e} de diffusion {\`a} valeurs dans un espace de Hilbert}, Ann. Inst. Henri Poincar{\'e}, 10 (1974), pp.~55--88.


\end{thebibliography}
\end{document}